# SPECIAL CLASSES OF SET CODES AND THEIR APPLICATIONS

**W. B. Vasantha Kandasamy**
e-mail: vasanthakandasamy@gmail.com
web: http://mat.iitm.ac.in/~wbv
www.vasantha.in

**Florentin Smarandache**
e-mail: smarand@unm.edu

**2008**

# SPECIAL CLASSES OF SET CODES AND THEIR APPLICATIONS

**W. B. Vasantha Kandasamy**

**Florentin Smarandache**

# CONTENTS

# PREFACE

In this book the authors introduce the notion of set codes, set bicodes and set n-codes. These are the most generalized notions of semigroup n-codes and group n-codes. Several types of set n-codes are defined. Several examples are given to enable the reader to understand the concept. These new classes of codes will find applications in cryptography, computer networking (where fragmenting of codes is to be carried out) and data storage (where confidentiality is to be maintained). We also describe the error detection and error correction of these codes. The authors feel that these codes would be appropriate to the computer-dominated world.

This book has three chapters. Chapter One gives basic concepts to make the book a self-contained one. In Chapter Two, the notion of set codes is introduced. The set bicodes and their generalization to set n-codes ($n \geq 3$) is carried out in Chapter Three. This chapter also gives the applications of these codes in the fields mentioned above. Illustrations of how these codes are applied are also given. The authors deeply acknowledge the unflinching support of Dr.K.Kandasamy, Meena and Kama.

This book is dedicated to the memory of the first author's father Mr.W.Balasubramanian, on his 100th birth anniversary. A prominent educationalist, who also severed in Ethiopia, he has been a tremendous influence in her life and the primary reason why she chose a career in mathematics.

W.B.VASANTHA KANDASAMY
FLORENTIN SMARANDACHE

## Chapter One

# BASIC CONCEPTS

This chapter has two sections. In section one we introduce basic concepts about set vector spaces, semigroup vector space, group vector spaces and set n-vector space. In section two we recall the basic definition and properties about linear codes and other special linear codes like Hamming codes, parity check codes, repetition codes etc.

### 1.1 Definition of Linear Algebra and its Properties

In this section we just recall the definition of linear algebra and enumerate some of its basic properties. We expect the reader to be well versed with the concepts of groups, rings, fields and matrices. For these concepts will not be recalled in this section.

Throughout this section, V will denote the vector space over F where F is any field of characteristic zero.

**DEFINITION 1.1.1:** *A vector space or a linear space consists of the following:*

i. *a field F of scalars.*
ii. *a set V of objects called vectors.*
iii. *a rule (or operation) called vector addition; which associates with each pair of vectors* $\alpha, \beta \in V$; $\alpha + \beta$ *in V, called the sum of* $\alpha$ *and* $\beta$ *in such a way that*
   a. *addition is commutative* $\alpha + \beta = \beta + \alpha$.
   b. *addition is associative* $\alpha + (\beta + \gamma) = (\alpha + \beta) + \gamma$.
   c. *there is a unique vector 0 in V, called the zero vector, such that*
   $$\alpha + 0 = \alpha$$
   *for all* $\alpha$ *in V.*
   d. *for each vector* $\alpha$ *in V there is a unique vector* $-\alpha$ *in V such that*
   $$\alpha + (-\alpha) = 0.$$
   e. *a rule (or operation), called scalar multiplication, which associates with each scalar c in F and a vector* $\alpha$ *in V, a vector* $c \bullet \alpha$ *in V, called the product of c and* $\alpha$, *in such a way that*

      1. $1 \bullet \alpha = \alpha$ *for every* $\alpha$ *in V.*
      2. $(c_1 \bullet c_2) \bullet \alpha = c_1 \bullet (c_2 \bullet \alpha)$.
      3. $c \bullet (\alpha + \beta) = c \bullet \alpha + c \bullet \beta$.
      4. $(c_1 + c_2) \bullet \alpha = c_1 \bullet \alpha + c_2 \bullet \alpha$.

   *for* $\alpha, \beta \in V$ *and* $c, c_1 \in F$.

It is important to note as the definition states that a vector space is a composite object consisting of a field, a set of 'vectors' and two operations with certain special properties. The same set of vectors may be part of a number of distinct vectors.

We simply by default of notation just say V a vector space over the field F and call elements of V as vectors only as matter of convenience for the vectors in V may not bear much resemblance to any pre-assigned concept of vector, which the reader has.

***Example 1.1.1:*** Let R be the field of reals. R[x] the ring of polynomials. R[x] is a vector space over R. R[x] is also a vector space over the field of rationals Q.

***Example 1.1.2:*** Let Q[x] be the ring of polynomials over the rational field Q. Q[x] is a vector space over Q, but Q[x] is clearly not a vector space over the field of reals R or the complex field **C**.

***Example 1.1.3:*** Consider the set V = R × R × R. V is a vector space over R. V is also a vector space over Q but V is not a vector space over **C**.

***Example 1.1.4:*** Let $M_{m \times n} = \{(a_{ij}) \mid a_{ij} \in Q\}$ be the collection of all m × n matrices with entries from Q. $M_{m \times n}$ is a vector space over Q but $M_{m \times n}$ is not a vector space over R or **C**.

***Example 1.1.5:*** Let

$$P_{3 \times 3} = \left\{ \begin{pmatrix} a_{11} & a_{12} & a_{13} \\ a_{21} & a_{22} & a_{23} \\ a_{31} & a_{32} & a_{33} \end{pmatrix} \middle| \; a_{ij} \in Q, 1 \le i \le 3, \; 1 \le j \le 3 \right\}.$$

$P_{3 \times 3}$ is a vector space over Q.

***Example 1.1.6:*** Let Q be the field of rationals and G any group. The group ring, QG is a vector space over Q.

*Remark:* All group rings KG of any group G over any field K are vector spaces over the field K.

We just recall the notions of linear combination of vectors in a vector space V over a field F. A vector β in V is said to be a linear combination of vectors $\nu_1, \ldots, \nu_n$ in V provided there exists scalars $c_1, \ldots, c_n$ in F such that

$$\beta = c_1\nu_1 + \ldots + c_n\nu_n = \sum_{i=1}^{n} c_i \ \nu_i \ .$$

Now we proceed on to recall the definition of subspace of a vector space and illustrate it with examples.

**DEFINITION 1.1.2:** *Let V be a vector space over the field F. A subspace of V is a subset W of V which is itself a vector space over F with the operations of vector addition and scalar multiplication on V.*

**DEFINITION 1.1.3:** *Let S be a set. V another set. We say V is a set vector space over the set S if for all v $\in$ V and for all s $\in$ S; vs and sv $\in$ V.*

***Example 1.1.7:*** Let V = {1, 2, ..., ∞} be the set of positive integers. S = {2, 4, 6, ..., ∞} the set of positive even integers. V is a set vector space over S. This is clear for sv = vs $\in$ V for all s $\in$ S and v $\in$ V.

It is interesting to note that any two sets in general may not be a set vector space over the other. Further even if V is a set vector space over S then S in general need not be a set vector space over V.

For from the above example 2.1.1 we see V is a set vector space over S but S is also a set vector space over V for we see for every s $\in$ S and v $\in$ V, v.s = s.v $\in$ S. Hence the above example is both important and interesting as one set V is a set vector space another set S and vice versa also hold good inspite of the fact S $\neq$ V.

Now we illustrate the situation when the set V is a set vector space over the set S. We see V is a set vector space over the set S and S is not a set vector space over V.

***Example 1.1.8:*** Let V = {$Q^+$ the set of all positive rationals} and S = {2, 4, 6, 8, …, ∞}, the set of all even integers. It is easily verified that V is a set vector space over S but S is not a set vector space over V, for $\frac{7}{3} \in$ V and 2 $\in$ S but $\frac{7}{3}.2 \notin$ S. Hence the claim.

**DEFINITION 1.1.4:** *Let V be a set with zero, which is non empty. Let G be a group under addition. We call V to be a group vector space over G if the following condition are true.*

1. *For every v $\in$ V and g $\in$ G gv and vg $\in$ V.*
2. *0.v = 0 for every v $\in$ V, 0 the additive identify of G.*

We illustrate this by the following examples.

***Example 1.1.9:*** Let V = {0, 1, 2, …, 15} integers modulo 15. G = {0, 5, 10} group under addition modulo 15. Clearly V is a group vector space over G, for gv ≡ $v_1$ (mod 15), for g $\in$ G and v, $v_1 \in$ V.

***Example 1.1.10:*** Let V = {0, 2, 4, …, 10} integers 12. Take G = {0, 6}, G is a group under addition modulo 12. V is a group vector space over G, for gv ≡ $v_1$ (mod 12) for g $\in$ G and v, $v_1 \in$ V.

***Example 1.1.11:*** Let

$$M_{2\times 3} = \left\{ \begin{pmatrix} a_1 & a_2 & a_3 \\ a_4 & a_5 & a_6 \end{pmatrix} \middle| a_i \in \{-\infty, \ldots, -4, -2, 0, 2, 4, \ldots, \infty\} \right\}.$$

Take G = Z be the group under addition. $M_{2\times 3}$ is a group vector space over G = Z.

***Example 1.1.12:*** Let V = Z × Z × Z = {(a, b, c) / a, b, c $\in$ Z}. V is a group vector space over Z.

***Example 1.1.13:*** Let V = {0, 1} be the set. Take G = {0, 1} the group under addition modulo two. V is a group vector space over G.

***Example 1.1.14:*** Let

$$V = \left\{\begin{pmatrix} 0 & 1 \\ 0 & 0 \end{pmatrix}, \begin{pmatrix} 1 & 1 \\ 0 & 0 \end{pmatrix}, \begin{pmatrix} 1 & 0 \\ 1 & 0 \end{pmatrix}, \begin{pmatrix} 0 & 1 \\ 1 & 0 \end{pmatrix}, \begin{pmatrix} 1 & 0 \\ 0 & 1 \end{pmatrix}, \begin{pmatrix} 0 & 0 \\ 1 & 0 \end{pmatrix}, \begin{pmatrix} 0 & 0 \\ 0 & 1 \end{pmatrix}, \begin{pmatrix} 0 & 0 \\ 0 & 0 \end{pmatrix}\right\}$$

be set. Take G = {0, 1} group under addition modulo 2. V is a group vector space over G.

***Example 1.1.15:*** Let

$$V = \left\{\begin{pmatrix} a_1 & a_2 & \ldots & a_n \\ 0 & 0 & \ldots & 0 \end{pmatrix}, \begin{pmatrix} 0 & 0 & \ldots & 0 \\ 0 & 0 & \ldots & 0 \end{pmatrix}, \begin{pmatrix} 0 & 0 & \ldots & 0 \\ b_1 & b_2 & \ldots & b_n \end{pmatrix} \middle| \ a_i, b_i \in Z; 1 \leq i \leq n\right\}$$

be the non empty set. Take G = Z the group of integers under addition. V is the group vector space over Z.

***Example 1.1.16:*** Let

$$V = \left\{\begin{pmatrix} a_1 & 0 & \ldots & 0 \\ a_2 & 0 & \ldots & 0 \end{pmatrix}, \begin{pmatrix} 0 & 0 & \ldots & 0 \\ 0 & 0 & \ldots & 0 \end{pmatrix}, \begin{pmatrix} 0 & b_1 & 0 & \ldots & 0 \\ 0 & b_2 & 0 & \ldots & 0 \end{pmatrix}, \ldots, \begin{pmatrix} 0 & 0 & \ldots & t_1 \\ 0 & 0 & \ldots & t_2 \end{pmatrix} \middle| \ a_i, b_i, \ldots, t_i \in Z;\ 1 \leq i \leq 2\right\}$$

be the set of 2 × n matrices of this special form. Let G = Z be the group of integers under addition. V is a group vector space over Z.

***Example 1.1.17:*** Let

$$V = \left\{\begin{pmatrix} a_1 & 0 \\ 0 & 0 \end{pmatrix}, \begin{pmatrix} 0 & a_2 \\ 0 & 0 \end{pmatrix}, \begin{pmatrix} 0 & 0 \\ a_3 & 0 \end{pmatrix}, \right.$$

$$\begin{pmatrix} 0 & 0 \\ 0 & 0 \end{pmatrix}, \begin{pmatrix} 0 & 0 \\ 0 & a_4 \end{pmatrix} \Bigg| a_1, a_2, a_3, a_4 \in Z \Bigg\}$$

be the set. Z = G the group of integers V is a group vector space over Z.

Now having seen examples of group vector spaces which are only set defined over an additive group.

***Example 1.1.18:*** Let V = {(0 1 0 0), (1 1 1), (0 0 0), (0 0 0 0), (1 1 0 0), (0 0 0 0 0), (1 1 0 0 1), (1 0 1 1 0)} be the set. Take $Z_2$ = G = {0, 1} group under addition modulo 2. V is a group vector space over $Z_2$.

***Example 1.1.19:*** Let

$$V = \left\{ \begin{pmatrix} a_1 & a_2 & a_2 \\ 0 & 0 & 0 \\ 0 & 0 & 0 \end{pmatrix}, \begin{pmatrix} b_1 & 0 & 0 \\ b_2 & 0 & 0 \\ b_3 & 0 & 0 \end{pmatrix}, \begin{pmatrix} 0 & c_1 & 0 \\ 0 & c_2 & 0 \\ 0 & c_3 & 0 \end{pmatrix}, \begin{pmatrix} 0 & 0 & a'_1 \\ 0 & 0 & a'_2 \end{pmatrix}, \right.$$

$$\left. \begin{pmatrix} 0 & 0 & 0 \\ 0 & 0 & 0 \\ 0 & 0 & 0 \end{pmatrix}, \begin{pmatrix} 0 & 0 & 0 \\ 0 & 0 & 0 \end{pmatrix} \Bigg| a_i\, b_i\, c_i \in Z;\; a'_1, a'_2 \in Z; 1 \le i \le 3 \right\}$$

be the set, Z = G the group under addition. V is just a set but V is a group vector space over Z.

It is important and interesting to note that this group vector spaces will be finding their applications in coding theory.

Now we proceed onto define the notion of substructures of group vector spaces.

**DEFINITION 1.1.5:** *Let V be the set which is a group vector space over the group G. Let P $\subseteq$ V be a proper subset of V. We say P is a group vector subspace of V if P is itself a group vector space over G.*

**DEFINITION 1.1.6:** *Let $V = V_1 \cup \ldots \cup V_n$, each $V_i$ is a distinct set with $V_i \not\subseteq V_j$ or $V_j \not\subseteq V_i$ if $i \neq j$; $1 \leq i, j \leq n$. Let each $V_i$ be a set vector space over the set S, i = 1, 2, …, n, then we call $V = V_1 \cup V_2 \cup \ldots \cup V_n$ to be the set n-vector space over the set S.*

We illustrate this by the following examples.

***Example 1.1.20:*** Let

$$
\begin{aligned}
V &= V_1 \cup V_2 \cup V_3 \cup V_4 \\
&= \{(1\ 1\ 1), (0\ 0\ 0), (1\ 0\ 0), (0\ 1\ 0), (1\ 1), (0\ 0), (1\ \ 1\ \ 1\ \ 1), \\
&\quad (1\ 0\ 0\ 0), (0\ 0\ 0)\} \cup \\
&\quad \left\{ \begin{pmatrix} a & a \\ a & a \end{pmatrix} \middle| a \in Z_2 \right\} \cup \left\{ \begin{pmatrix} a & a & a & a \\ a & a & a & a \end{pmatrix} \middle| a \in Z_2 = \{0, 1\} \right\} \\
&\quad \cup \{Z_2\,[x]\}.
\end{aligned}
$$

V is a set 4 vector space over the set S = {0, 1}.

***Example 1.1.21:*** Let

$$
\begin{aligned}
V &= V_1 \cup V_2 \cup V_3 \cup V_4 \cup V_5 \cup V_6 \\
&= \left\{ \begin{pmatrix} a & a \\ a & a \end{pmatrix} \middle| a \in Z^+ \right\} \cup \{Z^+ \times Z^+ \times Z^+\} \cup \{(a,\ a,\ a), \\
&\quad (a, a, a, a, a) \mid a \in Z^+\} \cup \left\{ \begin{pmatrix} a \\ a \\ a \\ a \end{pmatrix} \middle| a \in Z^+ \right\} \cup \{Z^+[x]\} \\
&\quad \cup \left\{ \begin{pmatrix} a_1 & a_2 & a_3 \\ a_4 & a_5 & a_6 \end{pmatrix} \middle| a_i \in Z^+; 1 \leq i \leq 6 \right\}
\end{aligned}
$$

be the set 6-vector space over the set $S = Z^+$.

***Example 1.1.22:*** Let

$$V = V_1 \cup V_2 \cup V_3$$

$$= \{Z_6 [x]\} \cup \{Z_6 \times Z_6 \times Z_6\} \cup \left\{\begin{pmatrix} a & a & a \\ a & a & a \end{pmatrix} \middle| a \in Z_6\right\}$$

be the set 3 vector space over the set S = {0, 2, 4}. We call this also as set trivector space over the set S. Thus when n = 3 we call the set n vector space as set trivector space.

We define set n-vector subspace of a set n-vector space V.

**DEFINITION 1.1.7:** *Let $V = V_1 \cup \ldots \cup V_n$ be a set n-vector space over the set S. If $W = W_1 \cup \ldots \cup W_n$ with $W_i \neq W_j$; $i \neq j$, $W_i \not\subseteq W_j$ and $W_j \not\subseteq W_i$, $1 \leq i, j \leq n$ and $W = W_1 \cup W_2 \cup \ldots \cup W_n \subseteq V_1 \cup V_2 \cup \ldots \cup V_n$ and W itself is a set n-vector space over the set S then we call W to be the set n vector subspace of V over the set S.*

We illustrate this by a simple example.

***Example 1.1.23:*** Let

$$\begin{aligned} V &= V_1 \cup V_2 \cup V_3 \cup V_4 \\ &= \{(a, a, a), (a, a) \mid a \in Z^+\} \cup \left\{\begin{pmatrix} a & a \\ b & b \end{pmatrix} \middle| a, b \in Z^+\right\} \cup \\ &\quad \{Z^+[x]\} \cup \left\{\begin{pmatrix} a \\ a \\ a \end{pmatrix} \middle| a \in Z^+\right\}, \end{aligned}$$

V is a set 4-vector space over the set $S = Z^+$. Take

$$\begin{aligned} W &= W_1 \cup W_2 \cup W_3 \cup W_4 \\ &= \{(a, a, a) \mid a \in Z^+\} \cup \left\{\begin{pmatrix} a & a \\ 0 & 0 \end{pmatrix} \middle| a \in Z^+\right\} \cup \end{aligned}$$

$$\{\text{all polynomial of even degree}\} \cup \left\{ \begin{pmatrix} a \\ a \\ a \end{pmatrix} \middle| \; a \in 2Z^+ \right\}$$

$$\subseteq \quad V_1 \cup V_2 \cup V_3 \cup V_4 \quad = \quad V,$$

is a set 4-vector subspace of V over the set S = $Z^+$.

We can find several set 4-vector subspaces of V.

Now we proceed on to define the n-generating set of a set n-vector space over the set S.

**DEFINITION 1.1.8:** *Let V = $V_1 \cup \ldots \cup V_n$ be a set n-vector space over set S. Let X = $X_1 \cup \ldots \cup X_n \subset V_1 \cup V_2 \cup \ldots \cup V_n$ = V. If each set $X_i$ generates $V_i$ over the set S, i = 1, 2, …, n then we say the set n vector space V = $V_1 \cup \ldots \cup V_n$ is generated by the n-set X = $X_1 \cup X_2 \cup \ldots \cup X_n$ and X is called the n-generator of V. If each of $X_i$ is of cardinality $n_i$, i = 1, 2, …, n then we say the n-cardinality of the set n vector space V is given by $|X_1| \cup \ldots \cup |X_n| = \{|X_1|, |X_2|, \ldots, |X_n|\} = \{(n_1, n_2, \ldots, n_n)\}$. If even one of the $X_i$ is of infinite cardinality we say the n-cardinality of V is infinite. Thus if all the sets $X_1, \ldots, X_n$ have finite cardinality then we say the n-cardinality of V is finite.*

## 1.2 Introduction to linear codes

In this section we give the definition and some basic properties about linear codes. In this section we recall the basic concepts about linear codes. Here we recall the definition of Repetition code, Parity check code, Hamming code, cyclic code, dual code and illustrate them with examples.

In this section we just recall the definition of linear code and enumerate a few important properties about them. We begin by describing a simple model of a communication transmission system given by the figure 1.2.1.

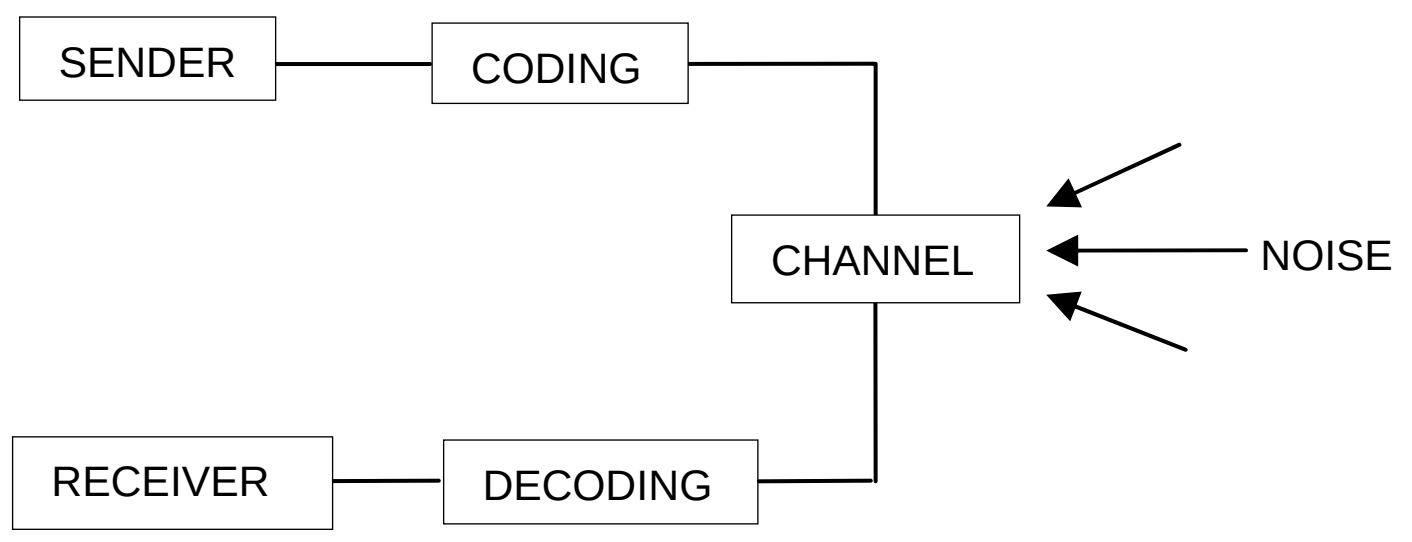

**Figure 1.2.1**

Messages go through the system starting from the source (sender). We shall only consider senders with a finite number of discrete signals (eg. Telegraph) in contrast to continuous sources (eg. Radio). In most systems the signals emanating from the source cannot be transmitted directly by the channel. For instance, a binary channel cannot transmit words in the usual Latin alphabet. Therefore an encoder performs the important task of data reduction and suitably transforms the message into usable form. Accordingly one distinguishes between source encoding the channel encoding. The former reduces the message to its essential(recognizable) parts, the latter adds redundant information to enable detection and correction of possible errors in the transmission. Similarly on the receiving end one distinguishes between channel decoding and source decoding, which invert the corresponding channel and source encoding besides detecting and correcting errors.

One of the main aims of coding theory is to design methods for transmitting messages error free cheap and as fast as possible. There is of course the possibility of repeating the message. However this is time consuming, inefficient and crude. We also note that the possibility of errors increases with an increase in the length of messages. We want to find efficient algebraic methods (codes) to improve the reliability of the transmission of messages. There are many types of algebraic codes; here we give a few of them.

Throughout this book we assume that only finite fields represent the underlying alphabet for coding. Coding consists of transforming a block of k message symbols $a_1, a_2, \ldots, a_k$; $a_i \in F_q$ into a code word $x = x_1 x_2 \ldots x_n$; $x_i \in F_q$, where $n \geq k$. Here the first $k_i$ symbols are the message symbols i.e., $x_i = a_i$; $1 \leq i \leq k$; the remaining $n - k$ elements $x_{k+1}, x_{k+2}, \ldots, x_n$ are check symbols or control symbols. Code words will be written in one of the forms x; $x_1, x_2, \ldots, x_n$ or $(x_1 x_2 \ldots x_n)$ or $x_1 x_2 \ldots x_n$. The check symbols can be obtained from the message symbols in such a way that the code words x satisfy a system of linear equations; $Hx^T = (0)$ where H is the given $(n - k) \times n$ matrix with elements in $F_q = Z_{p^n}$ $(q = p^n)$. A standard form for H is $(A, I_{n-k})$ with $n - k \times k$ matrix and $I_{n-k}$, the $n - k \times n - k$ identity matrix.

We illustrate this by the following example.

***Example 1.2.1:*** Let us consider $Z_2 = \{0, 1\}$. Take $n = 7$, $k = 3$. The message $a_1 a_2 a_3$ is encoded as the code word $x = a_1 a_2 a_3 x_4 x_5 x_6 x_7$. Here the check symbols $x_4 x_5 x_6 x_7$ are such that for this given matrix

$$H = \begin{bmatrix} 0 & 1 & 0 & 1 & 0 & 0 & 0 \\ 1 & 0 & 1 & 0 & 1 & 0 & 0 \\ 0 & 0 & 1 & 0 & 0 & 1 & 0 \\ 0 & 0 & 1 & 0 & 0 & 0 & 1 \end{bmatrix} = (A; I_4);$$

we have $Hx^T = (0)$ where

$x = a_1 a_2 a_3 x_4 x_5 x_6 x_7$.

$$a_2 + x_4 = 0$$
$$a_1 + a_3 + x_5 = 0$$
$$a_3 + x_6 = 0$$
$$a_3 + x_7 = 0.$$

Thus the check symbols $x_4 x_5 x_6 x_7$ are determined by $a_1 a_2 a_3$. The equation $Hx^T = (0)$ are also called check equations. If the message a = 1 0 0 then, $x_4 = 0$, $x_5 = 1$, $x_6 = 0$ and $x_7 = 0$. The code word x is 1 0 0 0 1 0 0. If the message a = 1 1 0 then $x_4 = 1$, $x_5 = 1$, $x_6 = 1 = x_7$. Thus the code word x = 1 1 0 1 1 0 0.

We will have altogether $2^3$ code words given by

| | |
|---|---|
| 0 0 0 0 0 0 0 | 1 1 0 1 1 0 0 |
| 1 0 0 0 1 0 0 | 1 0 1 0 0 1 1 |
| 0 1 0 1 1 0 0 | 0 1 1 1 1 1 1 |
| 0 0 1 0 1 1 1 | 1 1 1 1 0 1 1 |

**DEFINITION 1.2.1:** *Let H be an n – k × n matrix with elements in $Z_q$. The set of all n-dimensional vectors satisfying $Hx^T = (0)$ over $Z_q$ is called a linear code(block code) C over $Z_q$ of block length n. The matrix H is called the parity check matrix of the code C. C is also called a linear(n, k) code.*

*If H is of the form(A, $I_{n-k}$) then the k-symbols of the code word x is called massage(or information) symbols and the last n – k symbols in x are the check symbols. C is then also called a systematic linear(n, k) code. If q = 2, then C is a binary code. k/n is called transmission (or information) rate.*

*The set C of solutions of x of $Hx^T = (0)$. i.e., the solution space of this system of equations, forms a subspace of this system of equations, forms a subspace of $Z_q^n$ of dimension k. Since the code words form an additive group, C is also called a group code. C can also be regarded as the null space of the matrix H.*

***Example 1.2.2:*** (Repetition Code) If each codeword of a code consists of only one message symbol $a_1 \in Z_2$ and (n – 1) check symbols $x_2 = x_3 = \ldots = x_n$ are all equal to $a_1$ ($a_1$ is repeated n – 1 times) then we obtain a binary (n, 1) code with parity check matrix

$$H = \begin{bmatrix} 1 & 1 & 0 & 0 & \ldots & 1 \\ 0 & 0 & 1 & 0 & \ldots & 0 \\ 0 & 0 & 0 & 1 & \ldots & 0 \\ \vdots & \vdots & \vdots & \vdots & & \vdots \\ 1 & 0 & 0 & 0 & \ldots & 1 \end{bmatrix}.$$

There are only two code words in this code namely 0 0 … 0 and 1 1 …1.

If is often impracticable, impossible or too expensive to send the original message more than once. Especially in the transmission of information from satellite or other spacecraft, it is impossible to repeat such messages owing to severe time limitations. One such cases is the photograph from spacecraft as it is moving it may not be in a position to retrace its path. In such cases it is impossible to send the original message more than once. In repetition codes we can of course also consider code words with more than one message symbol.

***Example 1.2.3: (Parity-Check Code):*** This is a binary (n, n – 1) code with parity-check matrix to be H = (1 1 … 1). Each code word has one check symbol and all code words are given by all binary vectors of length n with an even number of ones. Thus if sum of the ones of a code word which is received is odd then atleast one error must have occurred in the transmission.

Such codes find its use in banking. The last digit of the account number usually is a control digit.

**DEFINITION 1.2.2:** *The matrix $G = (I_k, -A^T)$ is called a canonical generator matrix (or canonical basic matrix or encoding matrix) of a linear (n, k) code with parity check matrix $H = (A, I_{n-k})$. In this case we have $GH^T = (0)$.*

**DEFINITION 1.2.3:** *The Hamming distance d(x, y) between two vectors $x = x_1 x_2 \ldots x_n$ and $y = y_1 y_2 \ldots y_n$ in $F_q^n$ is the number of coordinates in which x and y differ. The Hamming weight $\omega(x)$ of a vector $x = x_1 x_2 \ldots x_n$ in $F_q^n$ is the number of non zero co ordinates in $x_i$. In short $\omega(x) = d(x, 0)$.*

We just illustrate this by a simple example.
Suppose x = [1 0 1 1 1 1 0] and y $\in$ [0 1 1 1 1 0 1 ] belong to $F_2^7$ then D(x, y) = (x ~ y) = (1 0 1 1 1 1 0) ~ (0 1 1 1 1 0 1) = (1~0, 0~1, 1~1, 1~1, 1~1, 1~0, 0~1) = (1 1 0 0 0 1 1) = 4. Now Hamming weight ω of x is ω(x) = d(x, 0) = 5 and ω(y) = d(y, 0) = 5.

**DEFINITION 1.2.4:** *Let C be any linear code then the minimum distance $d_{min}$ of a linear code C is given as*

$$d_{\min} = \min_{\substack{u,v \in C \\ u \neq v}} d(u,v).$$

*For linear codes we have*

$$d(u, v) = d(u - v, 0) = \omega(u - v).$$

Thus it is easily seen minimum distance of C is equal to the least weight of all non zero code words. A general code C of length n with k message symbols is denoted by C(n, k) or by a binary (n, k) code. Thus a parity check code is a binary (n, n – 1) code and a repetition code is a binary (n, 1) code.

If H = (A, $I_{n-k}$) be a parity check matrix in the standard form then G = ($I_k$, $-A^T$) is the canonical generator matrix of the linear (n, k) code.

The check equations (A, $I_{n-k}$) $x^T$ = (0) yield

$$\begin{bmatrix} x_{k+1} \\ x_{k+2} \\ \vdots \\ x_n \end{bmatrix} = -A \begin{bmatrix} x_1 \\ x_2 \\ \vdots \\ x_k \end{bmatrix} = -A \begin{bmatrix} a_1 \\ a_2 \\ \vdots \\ a_k \end{bmatrix}.$$

Thus we obtain

$$\begin{bmatrix} x_1 \\ x_2 \\ \vdots \\ x_n \end{bmatrix} = \begin{bmatrix} I_k \\ -A \end{bmatrix} \begin{bmatrix} a_1 \\ a_2 \\ \vdots \\ a_k \end{bmatrix}.$$

We transpose and denote this equation as

$$(x_1\ x_2\ \ldots\ x_n) = (a_1\ a_2\ \ldots\ a_k)\ (I_k, -A^7)$$
$$= (a_1\ a_2\ \ldots\ a_k)\ G.$$

We have just seen that minimum distance

$$d_{\min} = \min_{\substack{u,v \in C \\ u \neq v}} d(u,v).$$

If d is the minimum distance of a linear code C then the linear code of length n, dimension k and minimum distance d is called an (n, k, d) code.

Now having sent a message or vector x and if y is the received message or vector a simple decoding rule is to find the code word closest to x with respect to Hamming distance, i.e., one chooses an error vector e with the least weight. The decoding method is called "nearest neighbour decoding" and amounts to comparing y with all $q^k$ code words and choosing the closest among them. The nearest neighbour decoding is the maximum likelihood decoding if the probability p for correct transmission is > ½.

Obviously before, this procedure is impossible for large k but with the advent of computers one can easily run a program in few seconds and arrive at the result.

We recall the definition of sphere of radius r. The set $S_r(x) = \{y \in F_q^n \,/\, d(x, y) \leq r\}$ is called the sphere of radius r about $x \in F_q^n$ .

In decoding we distinguish between the detection and the correction of error. We can say a code can correct t errors and can detect t + s, s ≥ 0 errors, if the structure of the code makes it possible to correct up to t errors and to detect t + j, 0 < j ≤ s errors which occurred during transmission over a channel.

A mathematical criteria for this, given in the linear code is ; A linear code C with minimum distance $d_{min}$ can correct upto t errors and can detect t + j, 0 < j ≤ s, errors if and only if zt + s ≤ $d_{min}$ or equivalently we can say "A linear code C with minimum distance d can correct t errors if and only if $t = \left[\frac{(d-1)}{2}\right]$. The real problem of coding theory is not merely to minimize errors but to do so without reducing the transmission rate unnecessarily. Errors can be corrected by lengthening the code blocks, but this reduces the number of message symbols that can be sent per second. To maximize the transmission rate we want code blocks which are numerous enough to encode a given message alphabet, but at the same time no longer than is necessary to achieve a given Hamming distance. One of the main problems of coding theory is "Given block length n and Hamming distance d, find the maximum number, A(n, d) of binary blocks of length n which are at distances ≥ d from each other".

Let u = ($u_1, u_2, \ldots, u_n$) and v = ($v_1, v_2, \ldots, v_n$) be vectors in $F_q^n$ and let u.v = $u_1v_1 + u_2v_2 + \ldots + u_nv_n$ denote the dot product of u and v over $F_q^n$. If u.v = 0 then u and v are called orthogonal.

**DEFINITION 1.2.5:** *Let C be a linear (n, k) code over $F_q$. The dual(or orthogonal)code $C^\perp = \{u \mid u.v = 0$ for all $v \in C\}$, $u \in F_q^n$. If C is a k-dimensional subspace of the n-dimensional vector space $F_q^n$ the orthogonal complement is of dimension n – k and an (n, n – k) code. It can be shown that if the code C has a generator matrix G and parity check matrix H then $C^\perp$ has generator matrix H and parity check matrix G.*

Orthogonality of two codes can be expressed by $GH^T = (0)$.

**DEFINITION 1.2.6:** *For $a \in F_q^n$ we have $a + C = \{a + x / x \in C\}$. Clearly each coset contains $q^k$ vectors. There is a partition of $F_q^n$ of the form $F_q^n = C \cup \{a^{(1)} + C\} \cup \{a^{(2)} + C\} \cup \ldots \cup \{a^t + C\}$ for $t = q^{n-k} - 1$. If y is a received vector then y must be an element of one of these cosets say $a^i + C$. If the code word $x^{(1)}$ has been transmitted then the error vector*

$$e = y - x^{(1)} \in a^{(i)} + C - x^{(1)} = a^{(i)} + C.$$

*Now we give the decoding rule which is as follows.*

*If a vector y is received then the possible error vectors e are the vectors in the coset containing y. The most likely error is the vector $\bar{e}$ with minimum weight in the coset of y. Thus y is decoded as $\bar{x} = y - \bar{e}$. [23, 4]*

*Now we show how to find the coset of y and describe the above method. The vector of minimum weight in a coset is called the coset leader.*

*If there are several such vectors then we arbitrarily choose one of them as coset leader. Let $a^{(1)}, a^{(2)}, \ldots, a^{(t)}$ be the coset leaders. We first establish the following table*

| $x^{(1)} = 0$ | $x^{(2)} = 0$ | $\ldots$ | $x^{(q^k)}$ | *code words in C* |
|---|---|---|---|---|
| $a^{(1)} + x^{(1)}$ | $a^{(1)} + x^{(2)}$ | $\ldots$ | $a^{(1)} + x^{(q^k)}$ | |
| $\vdots$ | $\vdots$ | | $\vdots$ | *other cosets* |
| $a^{(t)} + x^{(1)}$ | $a^{(t)} + x^{(2)}$ | $\ldots$ | $a^{(t)} + x^{(q^k)}$ | |
| *coset leaders* | | | | |

*If a vector y is received then we have to find y in the table. Let* $y = a^{(i)} + x^{(j)}$*; then the decoder decides that the error* $\overline{e}$ *is the coset leader* $a^{(i)}$*. Thus y is decoded as the code word* $\overline{x} = y - \overline{e} = x^{(j)}$*. The code word* $\overline{x}$ *occurs as the first element in the column of y. The coset of y can be found by evaluating the so called syndrome.*

*Let H be parity check matrix of a linear (n, k) code. Then the vector* $S(y) = Hy^T$ *of length n–k is called syndrome of y. Clearly S(y) = (0) if and only if* $y \in C$*.*
$S(y^{(1)}) = S(y^{(2)})$ *if and only if* $y^{(1)} + C = y^{(2)} + C$ *.*

*We have the decoding algorithm as follows:*

*If* $y \in F_q^n$ *is a received vector find S(y), and the coset leader* $\overline{e}$ *with syndrome S(y). Then the most likely transmitted code word is* $\overline{x} = y - \overline{e}$ *we have* $d(\overline{x}, y)$ *.= min{d(x, y)/x* $\in$ *C}.*

We illustrate this by the following example.

***Example 1.2.4:*** Let C be a (5, 3) code where the parity check matrix H is given by

$$H = \begin{bmatrix} 1 & 0 & 1 & 1 & 0 \\ 1 & 1 & 0 & 0 & 1 \end{bmatrix}$$

and

$$G = \begin{bmatrix} 1 & 0 & 0 & 1 & 1 \\ 0 & 1 & 0 & 0 & 1 \\ 0 & 0 & 1 & 1 & 0 \end{bmatrix}.$$

The code words of C are

{(0 0 0  0 0), (1 0 0 1 1), (0 1 0 0 1), (0 0 1 1 0), (1 1 0 1 0), (1 0 1 0 1), (0 1 1 1 1), (1 1 1 0 0)}.

The corresponding coset table is

| Message | 000 | 100 | 010 | 001 | 110 | 101 | 011 | 111 |
|---|---|---|---|---|---|---|---|---|
| code words | 00000 | 10011 | 01001 | 00110 | 11010 | 10101 | 01111 | 11100 |
| other cosets | 10000 | 00011 | 11001 | 10110 | 01010 | 00101 | 11111 | 01100 |
| | 01000 | 11011 | 00001 | 01110 | 10010 | 11101 | 00111 | 10100 |
| | 00100 | 10111 | 01101 | 00010 | 11110 | 10001 | 01011 | 11000 |

coset leaders

If y = (1 1 1 1 0) is received, then y is found in the coset with the coset leader (0 0 1 0 0)
y + (0 0 1 0 0) = (1 1 1 1 0) + (0 0 1 0 0 ) = (1 1 0 1 0) is the corresponding message.

Now with the advent of computers it is easy to find the real message or the sent word by using this decoding algorithm.

A binary code $C_m$ of length $n = 2^m - 1$, $m \geq 2$ with $m \times 2^m - 1$ parity check matrix H whose columns consists of all non zero binary vectors of length m is called a binary Hamming code.

We give example of them.

***Example 1.2.5:*** Let

$$H = \begin{bmatrix} 1 & 0 & 1 & 1 & 1 & 1 & 0 & 0 & 1 & 0 & 1 & 1 & 0 & 0 & 0 \\ 1 & 1 & 0 & 1 & 1 & 1 & 1 & 0 & 0 & 1 & 0 & 0 & 1 & 0 & 0 \\ 1 & 1 & 1 & 0 & 1 & 0 & 1 & 1 & 0 & 0 & 1 & 0 & 0 & 1 & 0 \\ 1 & 1 & 1 & 1 & 0 & 0 & 0 & 1 & 1 & 1 & 0 & 0 & 0 & 0 & 1 \end{bmatrix}$$

which gives a $C_4(15, 11, 4)$ Hamming code.

Cyclic codes are codes which have been studied extensively.

Let us consider the vector space $F_q^n$ over $F_q$. The mapping

$$Z: F_q^n \to F_q^n$$

where Z is a linear mapping called a "cyclic shift" if $Z(a_0, a_1, \ldots, a_{n-1}) = (a_{n-1}, a_0, \ldots, a_{n-2})$

$A = (F_q[x], +, ., .)$ is a linear algebra in a vector space over $F_q$. We define a subspace $V_n$ of this vector space by

$$\begin{aligned} V_n &= \{v \in F_q[x] \,/\, \text{degree } v < n\} \\ &= \{v_0 + v_1x + v_2x^2 + \ldots + v_{n-1}x^{n-1} \,/\, v_i \in F_q; 0 \leq i \leq n-1\}. \end{aligned}$$

We see that $V_n \cong F_q^n$ as both are vector spaces defined over the same field $F_q$. Let Γ be an isomorphism

$\Gamma(v_0, v_1, \ldots, v_{n-1}) \to \{v_0 + v_1x + v_2x^2 + \ldots + v_{n-1}x^{n-1}\}$.

w: $F_q^n$ ℧ $F_q[x] / x^n - 1$

i.e., $w(v_0, v_1, \ldots, v_{n-1}) = v_0 + v_1x + \ldots + v_{n-1}x^{n-1}$.

Now we proceed onto define the notion of a cyclic code.

**DEFINITION 1.2.7:** *A k-dimensional subspace C of $F_q^n$ is called a cyclic code if Z(v) ∈ C for all v ∈ C that is v = $v_0, v_1, \ldots, v_{n-1}$ ∈ C implies ($v_{n-1}, v_0, \ldots, v_{n-2}$) ∈ C for v ∈ $F_q^n$.*

We just give an example of a cyclic code.

***Example 1.2.6:*** Let $C \subseteq F_2^7$ be defined by the generator matrix

$$G = \begin{bmatrix} 1 & 1 & 1 & 0 & 1 & 0 & 0 \\ 0 & 1 & 1 & 1 & 0 & 1 & 0 \\ 0 & 0 & 1 & 1 & 1 & 0 & 1 \end{bmatrix} = \begin{bmatrix} g^{(1)} \\ g^{(2)} \\ g^{(3)} \end{bmatrix}.$$

The code words generated by G are {(0 0 0 0 0 0 0), (1 1 1 0 1 0 0), (0 1 1 1 0 1 0), (0 0 1 1 1 0 1), (1 0 0 1 1 1 0), (1 1 0 1 0 0 1), (0 1 0 0 1 1 1), (1 0 1 0 0 1 1)}.

Clearly one can check the collection of all code words in C satisfies the rule if $(a_0 \ldots a_5) \in C$ then $(a_5\, a_0 \ldots a_4) \in C$ i.e., the codes are cyclic. Thus we get a cyclic code.

Now we see how the code words of the Hamming codes looks like.

***Example 1.2.7:*** Let

$$H = \begin{bmatrix} 1 & 0 & 0 & 1 & 1 & 0 & 1 \\ 0 & 1 & 0 & 1 & 0 & 1 & 1 \\ 0 & 0 & 1 & 0 & 1 & 1 & 1 \end{bmatrix}$$

be the parity check matrix of the Hamming (7, 4) code.

Now we can obtain the elements of a Hamming(7,4) code.

We proceed on to define parity check matrix of a cyclic code given by a polynomial matrix equation given by defining the generator polynomial and the parity check polynomial.

**DEFINITION 1.2.8:** *A linear code C in $V_n = \{v_0 + v_1x + \ldots + v_{n-1}x^{n-1} \mid v_i \in F_q, 0 \leq i \leq n-1\}$ is cyclic if and only if C is a principal ideal generated by $g \in C$.*

*The polynomial g in C can be assumed to be monic. Suppose in addition that $g / x^n - 1$ then g is uniquely determined and is called the generator polynomial of C. The elements of C are called code words, code polynomials or code vectors.*

*Let $g = g_0 + g_1x + \ldots + g_mx^m \in V_n$, $g / x^n - 1$ and deg $g = m < n$. Let C be a linear (n, k) code, with $k = n - m$ defined by the generator matrix,*

$$G = \begin{bmatrix} g_0 & g_1 & \cdots & g_m & 0 & \cdots & 0 \\ 0 & g_0 & \cdots & g_{m-1} & g_m & \cdots & 0 \\ \vdots & \vdots & & & & & \\ 0 & 0 & & g_0 & g_1 & & g_m \end{bmatrix} = \begin{bmatrix} g \\ xg \\ \\ x^{k-1}g \end{bmatrix}.$$

*Then C is cyclic. The rows of G are linearly independent and rank G = k, the dimension of C.*

***Example 1.2.8:*** Let g = $x^3 + x^2 + 1$ be the generator polynomial having a generator matrix of the cyclic(7,4) code with generator matrix

$$G = \begin{bmatrix} 1 & 0 & 1 & 1 & 0 & 0 & 0 \\ 0 & 1 & 0 & 1 & 1 & 0 & 0 \\ 0 & 0 & 1 & 0 & 1 & 1 & 0 \\ 0 & 0 & 0 & 1 & 0 & 1 & 1 \end{bmatrix}.$$

The codes words associated with the generator matrix is

0000000, 1011000, 0101100, 0010110, 0001011, 1110100, 1001110, 1010011, 0111010, 0100111, 0011101, 1100010, 1111111, 1000101, 0110001, 1101001.

The parity check polynomial is defined to be

$$h = \frac{x^7 - 1}{g}$$

$$h = \frac{x^7 - 1}{x^3 + x^2 + 1} = x^4 + x^3 + x^2 + 1.$$

$$\text{If } \frac{x^n - 1}{g} = h_0 + h_1 x + \ldots + h_k x^k.$$

the parity check matrix H related with the generator polynomial g is given by

$$H = \begin{bmatrix} 0 & \dots & 0 & h_k & \dots & h_1 & h_0 \\ 0 & \dots & h_k & h_{k-1} & & h_0 & 0 \\ \vdots & \vdots & \vdots & \vdots & & & \vdots \\ h_k & \dots & h_1 & h_0 & & \dots & 0 \end{bmatrix}.$$

For the generator polynomial g = $x^3 + x^2 + 1$ the parity check matrix

$$H = \begin{bmatrix} 0 & 0 & 1 & 1 & 1 & 0 & 1 \\ 0 & 1 & 1 & 1 & 0 & 1 & 0 \\ 1 & 1 & 1 & 0 & 1 & 0 & 0 \end{bmatrix}$$

where the parity check polynomial is given by $x^4 + x^3 + x^2 + 1 = \frac{x^7 - 1}{x^3 + x^2 + 1}$. It is left for the reader to verify that the parity check matrix gives the same set of cyclic codes.

We now proceed on to give yet another new method of decoding procedure using the method of best approximations.

We just recall this definition given by [4, 23, 39]. We just give the basic concepts needed to define this notion. We know that $F_q^n$ is a finite dimensional vector space over $F_q$. If we take $Z_2$ = (0, 1) the finite field of characteristic two. $Z_2^5 = Z_2 \times Z_2 \times Z_2 \times Z_2 \times Z_2$ is a 5 dimensional vector space over $Z_2$. Infact {(1 0 0 0 0), (0 1 0 0 0), (0 0 1 0 0), (0 0 0 1 0), (0 0 0 0 1)} is a basis of $Z_2^5$. $Z_2^5$ has only $2^5 = 32$ elements in it. Let F be a field of real numbers and V a vector space over F. An inner product on V is a function which assigns to each ordered pair of vectors α, β in V a scalar $\langle \alpha / \beta \rangle$ in F in such a way that for all α, β, γ in V and for all scalars c in F.

(a) $\langle \alpha + \beta / \gamma \rangle = \langle \alpha/\gamma \rangle + \langle \beta/\gamma \rangle$
(b) $\langle c\alpha / \beta \rangle = c\langle \alpha/\beta \rangle$
(c) $\langle \beta/\alpha \rangle = \langle \alpha/\beta \rangle$
(d) $\langle \alpha/\alpha \rangle > 0$ if $\alpha \neq 0$.

On V there is an inner product which we call the standard inner product. Let $\alpha = (x_1, x_2, \ldots, x_n)$ and $\beta = (y_1, y_2, \ldots, y_n)$

$$\langle \alpha / \beta \rangle = \sum_i x_i y_i \ .$$

This is called as the standard inner product. $\langle \alpha/\alpha \rangle$ is defined as norm and it is denoted by $||\alpha||$. We have the Gram-Schmidt orthogonalization process which states that if V is a vector space endowed with an inner product and if $\beta_1, \beta_2, \ldots, \beta_n$ be any set of linearly independent vectors in V; then one may construct a set of orthogonal vectors $\alpha_1, \alpha_2, \ldots, \alpha_n$ in V such that for each $k = 1, 2, \ldots, n$ the set $\{\alpha_1, \ldots, \alpha_k\}$ is a basis for the subspace spanned by $\beta_1, \beta_2, \ldots, \beta_k$ where $\alpha_1 = \beta_1$.

$$\alpha_2 = \beta_2 - \frac{\langle \beta_1 / \alpha_1 \rangle}{|| \alpha_1 ||^2} \alpha_1$$

$$\alpha_3 = \beta_3 - \frac{\langle \beta_3 / \alpha_1 \rangle}{|| \alpha_1 ||^2} \alpha_1 - \frac{\langle \beta_3 / \alpha_2 \rangle}{|| \alpha_2 ||^2} \alpha_2$$

and so on.

Further it is left as an exercise for the reader to verify that if a vector $\beta$ is a linear combination of an orthogonal sequence of non-zero vectors $\alpha_1, \ldots, \alpha_m$, then $\beta$ is the particular linear combination, i.e.,

$$\beta = \sum_{k=1}^{m} \frac{\langle \beta / \alpha_k \rangle}{|| \alpha_k ||^2} \alpha_k \ .$$

In fact this property that will be made use of in the best approximations.

We just proceed on to give an example.

***Example 1.2.9:*** Let us consider the set of vectors $\beta_1 = (2, 0, 3)$, $\beta_2 = (-1, 0, 5)$ and $\beta_3 = (1, 9, 2)$ in the space $R^3$ equipped with the standard inner product.

Define $\alpha_1 = (2, 0, 3)$

$$\alpha_2 = (-1, 0, 5) - \frac{\langle(-1, 0, 5)/(2, 0, 3)\rangle}{13}(2, 0, 3)$$

$$= (-1, 0, 5) - \frac{13}{13}(2, 0, 3) = (-3, 0, 2)$$

$$\alpha_3 = (1, 9, 2) - \frac{\langle(-1, 9, 2)/(2, 0, 3)\rangle}{13}(2, 0, 3)$$

$$- \frac{\langle(1, 9, 2)/(-3, 0, 2)\rangle}{13}(-3, 0, 2)$$

$$= (1, 9, 2) - \frac{8}{13}(2, 0, 3) - \frac{1}{13}(-3, 0, 2)$$

$$= (1, 9, 2) - \left(\frac{16}{13}, 0, \frac{24}{13}\right) - \left(\frac{3}{13}, 0, \frac{2}{13}\right)$$

$$= (1, 9, 2) - \left\{\left(\frac{16-3}{13}, 0, \frac{24+2}{13}\right)\right\}$$

$$= (1, 9, 2) - (1, 0, 2)$$

$$= (0, 9, 0).$$

Clearly the set {(2, 0, 3), (–3, 0, 2), (0, 9, 0)} is an orthogonal set of vectors.

Now we proceed on to define the notion of a best approximation to a vector $\beta$ in V by vectors of a subspace W where $\beta \notin W$. Suppose W is a subspace of an inner product space V and let $\beta$ be an arbitrary vector in V. The problem is to find a best possible approximation to $\beta$ by vectors in W. This means we want to find a vector $\alpha$ for which $||\beta - \alpha||$ is as small as possible subject to the restriction that $\alpha$ should belong to W. To be precisely in mathematical terms: A best approximation to $\beta$ by vectors in W is a vector $\alpha$ in W such that $||\beta - \alpha|| \leq ||\beta - \gamma||$ for every vector $\gamma$ in W ; W a subspace of V.

By looking at this problem in $R^2$ or in $R^3$ one sees intuitively that a best approximation to $\beta$ by vectors in W ought to be a vector $\alpha$ in W such that $\beta - \alpha$ is perpendicular (orthogonal) to W and that there ought to be exactly one such $\alpha$.

These intuitive ideas are correct for some finite dimensional subspaces, but not for all infinite dimensional subspaces.

We just enumerate some of the properties related with best approximation.

Let W be a subspace of an inner product space V and let $\beta$ be a vector in V.

(i) The vector $\alpha$ in W is a best approximation to $\beta$ by vectors in W if and only if $\beta - \alpha$ is orthogonal to every vector in W.

(ii) If a best approximation to $\beta$ by vectors in W exists, it is unique.

(iii) If W is finite-dimensional and $\{\alpha_1, \alpha_2, \ldots, \alpha_n\}$ is any orthonormal basis for W, then the vector $\alpha = \sum_k \frac{\langle \beta / \alpha_k \rangle}{\|\alpha_k\|^2} \alpha_k$, where $\alpha$ is the (unique)best approximation to $\beta$ by vectors in W.

Now this notion of best approximation for the first time is used in coding theory to find the best approximated sent code after receiving a message which is not in the set of codes used. Further we use for coding theory only finite fields $F_q$. i.e., $|F_q| < \infty$ . If C is a code of length n; C is a vector space over $F_q$ and C $\cong F_q^k \subseteq F_q^n$, k the number of message symbols in the code, i.e., C is a C(n, k) code. While defining the notion of inner product on vector spaces over finite fields we see all axiom of inner product defined over fields as reals or complex in general is not true. The main property which is not true is if $0 \neq x \in V$; the inner product of x with itself i.e., $\langle x / x \rangle = \langle x, x \rangle \neq 0$ if $x \neq 0$ is not true i.e., $\langle x / x \rangle = 0$ does not imply $x = 0$ .

To overcome this problem we define for the first time the new notion of pseudo inner product in case of vector spaces defined over finite characteristic fields.

**DEFINITION 1.2.9:** *Let V be a vector space over a finite field $F_p$ of characteristic p, p a prime. Then the pseudo inner product on V is a map $\langle , \rangle_p : V \times V \to F_p$ satisfying the following conditions.*

*1. $\langle x, x\rangle_p \geq 0$ for all $x \in V$ and $\langle x, x\rangle_p = 0$ does not in general imply $x = 0$.*

*2. $\langle x, y\rangle_p = \langle y, x\rangle_p$ for all $x, y \in V$.*

*3. $\langle x + y, z\rangle_p = \langle x, z\rangle_p + \langle y, z\rangle_p$ for all $x, y, z \in V$.*

*4. $\langle x, y + z\rangle_p = \langle x, y\rangle_p + \langle x, z\rangle_p$ for all $x, y, z \in V$.*

*5. $\langle \alpha.x, y\rangle_p = \alpha \langle x, y\rangle_p$ and*

*6. $\langle x, \beta.y\rangle_p = \beta\langle x, y\rangle_p$ for all $x, y, \in V$ and $\alpha, \beta \in F_p$.*

*Let V be a vector space over a field $F_p$ of characteristic p, p is a prime; then V is said to be a pseudo inner product space if there is a pseudo inner product $\langle,\rangle_p$ defined on V. We denote the pseudo inner product space by $(V, \langle,\rangle_p)$.*

Now using this pseudo inner product space $(V, \langle,\rangle_p)$ we proceed on to define pseudo-best approximation.

**DEFINITION 1.2.10:** *Let V be a vector space defined over the finite field $F_p$ (or $Z_p$). Let W be a subspace of V. For $\beta \in V$ and for a set of basis $\{\alpha_1, \ldots, \alpha_k\}$ of the subspace W the pseudo best approximation to $\beta$, if it exists is given by $\sum_{i=1}^{k} \langle \beta, \alpha_i \rangle_p \alpha_i$. If $\sum_{i=1}^{k} \langle \beta, \alpha_i \rangle_p \alpha_i = 0$, then we say the pseudo best approximation does not exist for this set of basis $\{\alpha_1, \alpha_2, \ldots, \alpha_k\}$. In this case we choose another set of basis for W say $\{\gamma_1, \gamma_2, \ldots, \gamma_k\}$ and calculate $\sum_{i=1}^{k} \langle \beta, \gamma_i \rangle_p \gamma_i$ and $\sum_{i=1}^{k} \langle \beta, \gamma_i \rangle_p \gamma_i$ is called a pseudo best approximation to $\beta$.*

*Note:* We need to see the difference even in defining our pseudo best approximation with the definition of the best approximation. Secondly as we aim to use it in coding theory and most of our linear codes take only their values from the field of characteristic two we do not need $\langle x, x\rangle$ or the norm to be divided by the pseudo inner product in the summation of finding the pseudo best approximation.

Now first we illustrate the pseudo inner product by an example.

***Example 1.2.10:*** Let V = $Z_2 \times Z_2 \times Z_2 \times Z_2$ be a vector space over $Z_2$. Define $\langle , \rangle_p$ to be the standard pseudo inner product on V; so if x = (1 0 1 1) and y = (1 1 1 1) are in V then the pseudo inner product of

$$\langle x, y\rangle_p = \langle (1\ 0\ 1\ 1\ ), (1\ 1\ 1\ 1)\rangle_p = 1 + 0 + 1 + 1 = 1.$$

Now consider

$$\langle x, x\rangle_p = \langle (1\ 0\ 1\ 1), (1\ 0\ 1\ 1)\rangle_p = 1 + 0 + 1 + 1 \neq 0$$

but

$$\langle y, y\rangle_p = \langle (1\ 1\ 1\ 1), (1\ 1\ 1\ 1)\rangle_p = 1 + 1 + 1 + 1 = 0.$$

We see clearly y ≠ 0, yet the pseudo inner product is zero.

Now having seen an example of the pseudo inner product we proceed on to illustrate by an example the notion of pseudo best approximation.

***Example 1.2.11:*** Let

$$V = Z_2^8 = \underbrace{Z_2 \times Z_2 \times \ldots \times Z_2}_{8 \text{ times}}$$

be a vector space over $Z_2$.
Now
W = {0 0 0 0 0 0 0 0), (1 0 0 0 1 0 11), (0 1 0 0 1 1 0 0), (0 0 1 0 0 1 1 1), (0 0 0 1 1 1 0 1), (1 1 0 0 0 0 1 0), (0 1 1 0 1 1 1 0), (0 0 1 1 1 0 1 0), (0 1 0 1 0 1 0 0), (1 0 1 0 1 1 0 0), (1 0 0 1 0 1 1 0), (1 1 1 0 0 1 0 1), (0 1 1 1 0 0 1 1), (1 1 0 1 1 1 1 1), (1 0 1 1 0 0 0 1), (1 1 1 1 1 0 0 0)}

be a subspace of V. Choose a basis of W as B = {$\alpha_1, \alpha_2, \alpha_3, \alpha_4$} where

$$\alpha_1 = (0\ 1\ 0\ 0\ 1\ 0\ 0\ 1),$$
$$\alpha_2 = (1\ 1\ 0\ 0\ 0\ 0\ 1\ 0),$$
$$\alpha_3 = (1\ 1\ 1\ 0\ 0\ 1\ 0\ 1)$$

and

$$\alpha_4 = (1\ 1\ 1\ 1\ 1\ 0\ 0\ 0).$$

Suppose β = (1 1 1 1 1 1 1 1) is a vector in V using pseudo best approximations find a vector in W close to β. This is given by α relative to the basis B of W where

$$\alpha = \sum_{k=-1}^{4} \langle \beta, \alpha_k \rangle_p \alpha_k$$

= ⟨(1 1 1 1 1 1 1 1), (0 1 0 0 1 0 0 1)⟩$_p$ $\alpha_1$ +
⟨(1 1 1 1 1 1 1 1), (1 1 0 0 0 0 1 0)⟩$_p$ $\alpha_2$ +
⟨(1 1 1 1 1 1 1 1), (1 1 1 0 0 1 0 1)⟩$_p$ $\alpha_3$ +
⟨(1 1 1 1 1 1 1 1), (1 1 1 1 1 0 0 0)⟩$_p$ $\alpha_4$.

= $1.\alpha_1 + 1.\alpha_2 + 1.\alpha_3 + 1.\alpha_4$.

= (0 1 0 0 1 0 0 1) + (1 1 0 0 0 0 1 0) + (1 1 1 0 0 1 0 1) + (1 1 1 1 1 0 0 0)

= (1 0 0 1 0 1 1 0) ∈ W.

Now having illustrated how the pseudo best approximation of a vector β in V relative to a subspace W of V is determined, now we illustrate how the approximately the nearest code word is obtained.

***Example 1.2.12:*** Let C = C(4, 2) be a code obtained from the parity check matrix

$$H = \begin{bmatrix} 1 & 0 & 1 & 0 \\ 1 & 1 & 0 & 1 \end{bmatrix}.$$

The message symbols associated with the code C are {(0, 0), (1, 0), (1, 0), (1, 1)}. The code words associated with H are C = {(0 0 0 0), (1 0 1 1), (0 1 0 1), (1 1 1 0)}. The chosen basis for C is B = {$\alpha_1$, $\alpha_2$} where $\alpha_1$ = (0 1 0 1) and $\alpha_2$ = (1 0 1 1). Suppose the received message is β = (1 1 1 1), consider $H\beta^T$ = (0 1) ≠ (0) so β ∉ C. Let α be the pseudo best approximation to β relative to the basis B given as

$$\alpha = \sum_{k=1}^{2} \langle \beta, \alpha_k \rangle_p \alpha_k \quad = \quad \langle (1\ 1\ 1\ 1), (0\ 1\ 0\ 1) \rangle_p \alpha_1$$

$$+ \langle (1\ 1\ 1\ 1), (1\ 0\ 1\ 1) \rangle_p \alpha_2.$$

$$= \quad (1\ 0\ 1\ 1)\ .$$

Thus the approximated code word is (1 0 1 1).
This method could be implemented in case of algebraic linear bicodes and in general to algebraic linear n-codes; $n \geq 3$.

**Chapter Two**

# NEW CLASS OF SET CODES AND THEIR PROPERTIES

In this chapter we for the first time introduce a special class of set codes and give a few properties enjoyed by them. Throughout this chapter we assume the set codes are defined over the set S = {0, 1} = $Z_2$, i.e., all code words have only binary symbols. We now proceed on to define the notion of set codes.

**DEFINITION 2.1 :** *Let C = {($x_1 \ldots x_{r_1}$), ($x_1 \ldots x_{r_2}$), …, ($x_1 \ldots x_{r_n}$)} be a set of ($r_1$, …, $r_n$) tuples with entries from the set S = {0, 1} where each $r_i$-tuple ($x_1$, …, $x_{r_i}$) has some $k_i$ message symbols and $r_i - k_i$ check symbols $1 \leq i \leq n$. We call C the set code if C is a set vector space over the set S = {0, 1}.*

The following will help the reader to understand more about these set codes.

1. $r_i = r_j$ even if $i \neq j$ ; $1 \leq i, j \leq n$
2. $k_i = k_j$ even if $i \neq j$ ; $1 \leq i, j \leq n$
3. Atleast some $r_i \neq r_t$ when $i \neq t$; $1 \leq i, t \leq n$.

We first illustrate this by some examples before we proceed onto give more properties about the set codes.

***Example 2.1 :*** Let C = {(1 1 0 0), (0 0 0 0), (1 1 1 0 0), (0 1 1 0 1), (0 0 0 0 0), (1 1 1), (0 1 0), (0 0 0)}. C is a set code over the set S = {0, 1}.

***Example 2.2:*** Let C = {(0 1 1 0), (0 0 0 1), (1 1 0 1), (1 1 1 1 1 1), (0 0 0 0 0 0), (1 0 1 0 1 0), (0 1 0 1 0 1), (0 0 0 0)}. C is also a set code over the set S = {0, 1}.

We see the set codes will have code words of varying lengths. We call the elements of the set code as set code words. In general for any set code C we can have set code words of varying lengths. As in case of usual binary codes we do not demand the length of every code word to be of same length.

Further as in case of usual codes we do not demand the elements of the set codes to form a subgroup i.e., they do not form a group or a subspace of a finite dimensional vector space. They are just collection of $(r_1, \ldots, r_n)$ tuples with no proper usual algebraic structure.

***Example 2.3:*** Let C = {(1 1 1), (0 0 0), (1 1 1 1 1 1), (0 0 0 0 0 0), (1 1 1 1 1 1 1), (0 0 0 0 0 0 0)} be a set code over the set S = {0, 1}.

Now we give a special algebraic structure enjoyed by these set codes.

**DEFINITION 2.2:** *Let $C = \{(x_1 \ldots x_{r_1}), \ldots, (x_1 \ldots x_{r_n})\}$ be a set code over the set S = {0, 1} where $x_i \in \{0, 1\}$; $1 \leq i \leq r_i, \ldots, r_n$. We demand a set of matrices $H = \{H_1, \ldots, H_t \mid H_i$ takes its entries from the set {0,1}} and each set code word $x \in C$ is such that there is some matrix $H_i \in H$ with $H_i x^t = 0$, $1 < i < t$. We do* *not demand the same $H_i$ to satisfy this relation for every set code word from C. Further the set $H = \{H_1, \ldots, H_t\}$ will not form a set vector space over the set {0,1}. This set H is defined as the set parity check matrices for the set code C.*

We illustrate this by an example.

***Example 2.4:*** Let V = {(1 1 1 0 0 0), (1 0 1 1 0 1), (0 0 0 0 0 0), (0 1 1 0 1 1), (1 0 1 1), (0 0 0 0), (1 1 1 0)} be a set code. The set parity check matrix associated with V is given by

$$H = \{H_{1,} H_2\} =$$

$$\{H_1 = \begin{pmatrix} 0 & 1 & 1 & 1 & 0 & 0 \\ 1 & 0 & 1 & 0 & 1 & 0 \\ 1 & 1 & 0 & 0 & 0 & 1 \end{pmatrix},$$

and

$$H_2 = \begin{pmatrix} 1 & 0 & 1 & 0 \\ 1 & 1 & 0 & 1 \end{pmatrix}\}.$$

The following facts are important to be recorded.

(1) As in case of the ordinary code we don't use parity check matrix to get all the code words using the message symbols. Infact the set parity check matrix is used only to find whether the received code word is correct or error is present.

(2) Clearly V, the set code does not contain in general all the code words associated with the set parity check matrix H.

(3) The set codes are just set, it is not compatible even under addition of code words. They are just codes as they cannot be added for one set code word may be of length $r_1$ and another of length $r_2$; $r_1 \neq r_2$.

(4) The set codes are handy when one intends to send messages of varying lengths simultaneously.

***Example 2.5:*** V = {(0 0 0 0 0 0), (0 0 1 1 1 0), (1 0 0 0 1 1), (0 1 0 1 0 1), (1 1 0 0 1), (1 1 1 0 0), (0 1 0 1 1), (0 0 0 0 0 0)} is a set code with set parity check matrix H = ($H_{1,}$ $H_2$).

***Example 2.6:*** Let V = {(0 0 0 0 0 0 0), (1 1 1 1 1 1 1), (1 1 1 1 1), (0 0 0 0 0), (1 1 1 1), (0 0 0 0)} be a set code with set parity check matrix

$$H = \left\{ \begin{pmatrix} 1 & 1 & 0 & 0 & 0 & 0 & 0 & 0 \\ 1 & 0 & 1 & 0 & 0 & 0 & 0 & 0 \\ 1 & 0 & 0 & 1 & 0 & 0 & 0 & 0 \\ 1 & 0 & 0 & 0 & 1 & 0 & 0 & 0 \\ 1 & 0 & 0 & 0 & 0 & 1 & 0 & 0 \\ 1 & 0 & 0 & 0 & 0 & 0 & 1 & 0 \\ 1 & 0 & 0 & 0 & 0 & 0 & 0 & 1 \end{pmatrix}, \right.$$

$$\begin{pmatrix} 1 & 1 & 0 & 0 & 0 & 0 \\ 1 & 0 & 1 & 0 & 0 & 0 \\ 1 & 0 & 0 & 1 & 0 & 0 \\ 1 & 0 & 0 & 0 & 1 & 0 \\ 1 & 0 & 0 & 0 & 0 & 1 \end{pmatrix}$$

and

$$\left. \begin{pmatrix} 1 & 1 & 0 & 0 & 0 \\ 1 & 0 & 1 & 0 & 0 \\ 1 & 0 & 0 & 1 & 0 \\ 1 & 0 & 0 & 0 & 1 \end{pmatrix} \right\}.$$

Now having defined set codes we make the definition of special classes of set codes.

**DEFINITION 2.3:** *Let V = {($y_1$ ,…, $y_{r_t}$ ), ($x_1$ ,…, $x_{r_t}$ )| $y_i$ = 0, 1 ≤ i ≤ $r_t$; $x_i$ = 1, 1 ≤ i ≤ $r_t$ and t = 1, 2, …, n} be a set code where either each of the tuples are zeros or ones. The set parity check matrix H = {$H_1$, …, $H_n$} where $H_i$ is a ($r_i$ – 1) × $r_i$ matrix of the form;*

$$\begin{pmatrix} 1 & 1 & 0 & \cdots & 0 \\ 1 & 0 & 1 & \cdots & 0 \\ \vdots & \vdots & \vdots & & \vdots \\ 1 & 0 & 0 & \cdots & 1 \end{pmatrix}$$

*where first column has only ones and the rest is a ($r_i$ – 1) × ($r_i$ – 1) identity matrix; i = 1, 2, …, n. We call this set code as the repetition set code.*

We illustrate this by some examples.

***Example 2.7:*** Let V = {(0 0 0 0), (1 1 1 1), (1 1 1 1 1 1), (0 0 0 0 0 0), (1 1 1 1 1 1 1 1 1), (0 0 0 0 0 0 0 0 0), (1 1 1 1 1 1 1), (0 0 0 0 0 0 0)} be a set code with set parity check matrix H = {$H_1$, $H_2$, $H_3$, $H_4$} where

$$H_1 = \begin{pmatrix} 1 & 1 & 0 & 0 & 0 \\ 1 & 0 & 1 & 0 & 0 \\ 1 & 0 & 0 & 1 & 0 \\ 1 & 0 & 0 & 0 & 1 \end{pmatrix},$$

$$H_2 = \begin{pmatrix} 1 & 1 & 0 & 0 & 0 & 0 & 0 \\ 1 & 0 & 1 & 0 & 0 & 0 & 0 \\ 1 & 0 & 0 & 1 & 0 & 0 & 0 \\ 1 & 0 & 0 & 0 & 1 & 0 & 0 \\ 1 & 0 & 0 & 0 & 0 & 1 & 0 \\ 1 & 0 & 0 & 0 & 0 & 0 & 1 \end{pmatrix},$$

$$H_3 = \begin{pmatrix} 1 & 1 & 0 & 0 & 0 & 0 & 0 & 0 & 0 & 0 \\ 1 & 0 & 1 & 0 & 0 & 0 & 0 & 0 & 0 & 0 \\ 1 & 0 & 0 & 1 & 0 & 0 & 0 & 0 & 0 & 0 \\ 1 & 0 & 0 & 0 & 1 & 0 & 0 & 0 & 0 & 0 \\ 1 & 0 & 0 & 0 & 0 & 1 & 0 & 0 & 0 & 0 \\ 1 & 0 & 0 & 0 & 0 & 0 & 1 & 0 & 0 & 0 \\ 1 & 0 & 0 & 0 & 0 & 0 & 0 & 1 & 0 & 0 \\ 1 & 0 & 0 & 0 & 0 & 0 & 0 & 0 & 1 & 0 \\ 1 & 0 & 0 & 0 & 0 & 0 & 0 & 0 & 0 & 1 \end{pmatrix}$$

and

$$H_4 = \begin{pmatrix} 1 & 1 & 0 & 0 & 0 & 0 & 0 & 0 \\ 1 & 0 & 1 & 0 & 0 & 0 & 0 & 0 \\ 1 & 0 & 0 & 1 & 0 & 0 & 0 & 0 \\ 1 & 0 & 0 & 0 & 1 & 0 & 0 & 0 \\ 1 & 0 & 0 & 0 & 0 & 1 & 0 & 0 \\ 1 & 0 & 0 & 0 & 0 & 0 & 1 & 0 \\ 1 & 0 & 0 & 0 & 0 & 0 & 0 & 1 \end{pmatrix}.$$

V is a repetition set code.

***Example 2.8:*** Let V = {(0 0 0), (1 1 1), (1 1 1 1 1), (0 0 0 0 0), (1 1 1 1), (0 0 0 0), (1 1 1 1 1 1 1 1), (0 0 0 0 0 0 0 0)} be a set code with the set parity check matrix H = {$H_1$, $H_2$, $H_3$ and $H_4$}, where

$$H_1 = \begin{pmatrix} 1 & 1 & 0 & 0 \\ 1 & 0 & 1 & 0 \\ 1 & 0 & 0 & 1 \end{pmatrix}, H_2 = \begin{pmatrix} 1 & 1 & 0 & 0 & 0 & 0 \\ 1 & 0 & 1 & 0 & 0 & 0 \\ 1 & 0 & 0 & 1 & 0 & 0 \\ 1 & 0 & 0 & 0 & 1 & 0 \\ 1 & 0 & 0 & 0 & 0 & 1 \end{pmatrix},$$

$$H_3 = \begin{pmatrix} 1 & 1 & 0 & 0 & 0 \\ 1 & 0 & 1 & 0 & 0 \\ 1 & 0 & 0 & 1 & 0 \\ 1 & 0 & 0 & 0 & 1 \end{pmatrix}$$

and

$$H_4 = \begin{pmatrix} 1 & 1 & 0 & 0 & 0 & 0 & 0 & 0 & 0 \\ 1 & 0 & 1 & 0 & 0 & 0 & 0 & 0 & 0 \\ 1 & 0 & 0 & 1 & 0 & 0 & 0 & 0 & 0 \\ 1 & 0 & 0 & 0 & 1 & 0 & 0 & 0 & 0 \\ 1 & 0 & 0 & 0 & 0 & 1 & 0 & 0 & 0 \\ 1 & 0 & 0 & 0 & 0 & 0 & 1 & 0 & 0 \\ 1 & 0 & 0 & 0 & 0 & 0 & 0 & 1 & 0 \\ 1 & 0 & 0 & 0 & 0 & 0 & 0 & 0 & 1 \end{pmatrix}.$$

This set code V is again a repetition set code.

Now we proceed on to define the notion of parity check set code.

**DEFINITION 2.4:** *Let V = {Some set of code words from the binary ($n_i$, $n_i$ – 1) code; i = 1, 2, …, t} with set parity check matrix*

$$\begin{aligned} H &= \{H_1, \ldots, H_t\} \\ &= \{(1\ 1\ 1 \ldots 1), (1\ 1 \ldots 1), \ldots, (1\ 1 \ldots 1)\}; \end{aligned}$$

*t set of some $n_i$ tuples with ones; i = 1, 2, …, t. We call V the parity check set code.*

We illustrate this by some examples.

***Example 2.9:*** Let V = {(1 1 1 1), (1 1 0 0), (1 0 0 1), (1 1 1 1 0), (1 1 0 0 0), (0 1 1 0 0), (1 1 1 1 1 1), (1 1 0 0 1 1), (0 1 1 1 1 0), (0 0 0 0), (0 0 0 0 0), (0 0 0 0 0 0)} be set code with set parity check matrix H = {$H_1$, $H_2$, $H_3$} where $H_1$ = (1 1 1 1), $H_2$ = (1 1 1 1 1) and $H_3$ = (1 1 1 1 1 1). V is a parity check set code.

We give yet another example.

***Example 2.10:*** Let V = {(1 1 1 1 1 1 0), (1 1 0 0 0 1 1), (1 0 0 0 0 0 1), (0 0 0 0 0 0 0), (1 1 1 1 0 0 0 0 0), (1 1 1 1 1 1 0 0 0), (0 0 0 0 0 0 0 0 0), (1 1 0 0 1 1 0 1 1), (1 1 1 1 0 0), (1 1 0 0 1 1), (1 0 0 0 0 1), (0 0 0 0 0 0)} be the set code with set parity check matrix H = {$H_{1,}$ $H_2$, $H_3$} where $H_1$ = (1 1 1 1 1 1 1), $H_2$ = (1 1 1 1 1 1 1 1 1) and $H_3$ (1 1 1 1 1 1). V is clearly a parity check set code.

Now we proceed onto define the notion of Hamming distance in set codes.

**DEFINITION 2.5:** *Let V = {$X_1$, …, $X_n$} be a set code. The set Hamming distance between two vectors $X_i$ and $X_j$ in V is defined if and only if both $X_i$ and $X_j$ have same number of entries and d($X_i$, $X_j$) is the number of coordinates in which $X_i$ and $X_j$ differ.*

*The set Hamming weight $\omega(X_i)$ of a set vector $X_i$ is the number of non zero coordinates in $X_i$. In short $\omega(X_i) = d(X_i, 0)$.*

Thus it is important to note that as in case of codes in set codes we cannot find the distance between any set codes. The distance between any two set code words is defined if and only if the length of both the set codes is the same.

We first illustrate this situation by the following example.

***Example 2.11:*** Let V = {(1 1 1 0 0), (0 0 0 0 0), (1 1 0 0 1), (1 1 1 1), (0 0 0 0), (1 1 1 0), (0 1 0 1), (0 0 0 0 0 0 0 0), (1 1 0 0 1 1 0 0), (1 1 0 1 1 0 0 0)} be a set code. The set Hamming distance between (1 1 1 0 0) and (0 1 0 1) cannot be defined as both of them are of different lengths. Now d {(1 1 1 1), (0 1 0 1)} = 2,

| | | |
|---|---|---|
| d((1 1 1 0 0), (1 1 0 0 1)) | = | 2, |
| d{(1 1 0 0 1 1 0 0), (1 1 0 1 1 0 0 0)} | = | 2 and |
| d((1 1 1 1), (0 0 0 0)) | = | 4; where as |
| d ((1 1 0 0 1 1 0 0), (1 1 1 1)) | | |

is not defined. We see the set weight ω(1 1 0 0 1 1 0 0) = 4, ω(1 1 1 1) = 4 and ω(1 1 1 0) = 3.

How to detect errors in set codes. The main function of set parity check matrices is only to detect error and nothing more. Thus if x is a transmitted set code word and y is the received set code word then e = (y – x) is called the set error word or error set word or error.

Thus if V is a set code with set parity check matrix H = $\{H_1, \ldots, H_k\}$. If y is a received message we say y has error if $H_t y^t \neq (0)$. If $H_t y^t = (0)$; we assume the received set code word is the correct message.

Thus in set codes correcting error is little difficult so we define certain special type of set codes with preassigned weights.

Now we proceed on to define the notion of binary Hamming set code.

**DEFINITION 2.6:** *Let V = $\{X_1, \ldots, X_n\}$ be a set code with set parity check matrix H = $\{H_1, \ldots, H_p\}$ where each $H_t$ is a $m_t \times (2^{m_t} - 1)$ parity check matrix whose columns consists of all non zero binary vectors of length $m_i$. We don't demand all code words associated with each $H_t$ to be present in V, only a choosen few alone are present in V; $1 \leq t \leq p$. We call this V as the binary Hamming set code.*

We illustrate this by some simple examples.

***Example 2.12:*** Let V = {(0 0 0 0 0 0 0), (1 0 0 1 1 0 0), (1 1 1 1 1 1 1), (0 0 0 0 0 0 0 0), (1 1 1 1 1 1 1 1), (1 0 1 0 1 0 1 0)} be a set code with associated set parity matrix H = $(H_1, H_2)$ where

$$H_1 = \begin{pmatrix} 0 & 0 & 0 & 1 & 1 & 1 & 1 \\ 0 & 1 & 1 & 0 & 0 & 1 & 1 \\ 1 & 0 & 1 & 0 & 1 & 0 & 1 \end{pmatrix}$$

and

$$H_2 = \begin{pmatrix} 1 & 1 & 1 & 1 & 1 & 1 & 1 & 1 \\ 0 & 0 & 0 & 1 & 1 & 1 & 1 & 0 \\ 0 & 1 & 1 & 0 & 0 & 1 & 1 & 0 \\ 1 & 0 & 1 & 0 & 1 & 0 & 1 & 0 \end{pmatrix}.$$

H is a binary Hamming set code.

Now we give yet another example of a binary Hamming set code.

***Example 2.13:*** Let V = {(0 0 0 0 0 0 0), (1 0 0 1 1 0 1), (0 1 0 1 0 1 1), (0 0 0 0 0 0 0 0 0 0 0 0 0 0 0), (1 0 0 0 1 0 0 1 1 0 1 0 1 1 1), (0 0 0 1 0 0 1 1 0 1 0 1 1 1 1)} be a set code with set parity check matrix H = ($H_1$, $H_2$) where

$$H_1 = \begin{pmatrix} 1 & 0 & 0 & 1 & 1 & 0 & 1 \\ 0 & 1 & 0 & 1 & 0 & 1 & 1 \\ 0 & 0 & 1 & 0 & 1 & 1 & 1 \end{pmatrix}$$

and

$$H_2 = \begin{pmatrix} 1 & 0 & 0 & 0 & 1 & 0 & 0 & 1 & 1 & 0 & 1 & 0 & 1 & 1 & 1 \\ 0 & 1 & 0 & 0 & 1 & 1 & 0 & 1 & 0 & 1 & 1 & 1 & 1 & 0 & 0 \\ 0 & 0 & 1 & 0 & 0 & 1 & 1 & 0 & 1 & 0 & 1 & 1 & 1 & 1 & 0 \\ 0 & 0 & 0 & 1 & 0 & 0 & 1 & 1 & 0 & 1 & 0 & 1 & 1 & 1 & 1 \end{pmatrix}.$$

Clearly V is a Hamming set code.

Now we proceed on to define the new class of codes called m – weight set code (m ≥ 2).

**DEFINITION 2.7:** *Let V = {$X_1$, …, $X_n$} be a set code with a set parity check matrix H = {$H_1$, …, $H_t$}; t < n. If the set Hamming weight of each and every $X_i$ in V is m, m < n and m is less than the least length of set code words in V. i.e., $\omega(X_i) = m$ for i = 1, 2, …, n and $X_i \neq (0)$. Then we call V to be a m-weight set code.*

These set codes will be useful in cryptology, in computers and in channels in which a minimum number of messages is to be preserved.

We will illustrate them by some examples.

***Example 2.14:*** Let V = {(0 0 0 0 0 0), (1 1 1 1 0 0), (0 1 1 1 1 0), (1 1 1 1 0 0 0 0), (0 1 1 1 1 0 0 0), (1 1 0 0 0 0 1 1), (1 1 1 0 0 0 1), (0 0 1 1 1 1 0), (0 0 0 0 0 0 0), (1 0 1 0 1 0 1)} be a set code with set parity check matrix H = {$H_1$, $H_2$, $H_3$}. This set code is a 4 weight set code for we see every non zero code word is of weight 4.

One of the advantages of m weight set code is that error detection becomes very easy. So for these m weight set codes even if set parity check matrix is not provided we can detect the errors easily.

We give yet another example of a m-weight set code.

***Example 2.15:*** Let V = {(1 1 1 1 0 1), (1 1 1 1 1 0 0 0), (0 0 0 0 0 0 0 0), (0 0 0 0 0 0), (1 0 1 0 1 0 1 1), (1 1 0 0 1 1 1 0), (0 1 1 1 1 1), (1 1 1 0 1 1), (1 1 1 1 1 0 0), (0 0 1 1 1 1 1), (0 0 0 0 0 0 0)}. V is a set code which is a 5-weight set code.

Another usefulness of this m-weight set code is that these codes are such that the error is very quickly detected; they can be used in places where several time transmission is possible. In places like passwords in e-mails, etc; these codes is best suited.

Now having defined m-weight set codes we now proceed on to define cyclic set codes.

**DEFINITION 2.8:** *Let V = {$x_1$, …, $x_n$} be n-set code word of varying length, if $x_i = (x_{i_1}, \ldots, x_{i_n}) \in V$ then $(x_{i_n}, x_{i_1}, \ldots, x_{i_{n-1}}) \in V$ for $1 \le i \le n$. Such set codes V are defined as cyclic set codes.*

*For cyclic set codes if H = {$H_1$, …, $H_t$}; t < n is the set parity check matrix then we can also detect errors.*

Now if in a set cyclic code we have weight of each set code word is m then we call such set codes as m weight cyclic set codes.

We now illustrate this by the following example.

***Example 2.16:*** Let V = {(1 0 1 0 1 0), (0 1 0 1 0 1), (0 0 0 0 0 0), (0 1 1 1 0 0 0 0), (0 0 1 1 1 0 0 0), (0 0 0 1 1 1 0 0), (0 0 0 0 1 1 1 0), (0 0 0 0 0 1 1 1), (1 0 0 0 0 0 1 1), (1 1 0 0 0 0 0 1), (1 1 1 0 0 0 0 0), (0 0 0 0 0 0 0 0)}.V is a 3 weight cyclic set code.

We do not even need the set parity check matrix to find error during transmission as the weight of the codes takes care of it.

***Example 2.17:*** Let V = {(0 0 0 0 0 0 0), (0 0 1 0 1 1 1), (0 1 0 1 1 1 0), (1 0 1 1 1 0 0), (1 0 0 1 0 0), (0 0 0 0 0 0), (0 1 0 0 1 0), (0 0 1 0 0 1)} be a set cyclic code with associated set parity check matrix H = ($H_1$ $H_2$) where

$$H_1 = \begin{pmatrix} 0 & 0 & 1 & 0 & 1 & 1 & 1 \\ 0 & 1 & 0 & 1 & 1 & 1 & 0 \\ 1 & 0 & 1 & 1 & 1 & 0 & 0 \end{pmatrix}$$

and

$$H_2 = \begin{pmatrix} 1 & 0 & 0 & 1 & 0 & 0 \\ 0 & 1 & 0 & 0 & 1 & 0 \\ 0 & 0 & 1 & 0 & 0 & 1 \end{pmatrix}.$$

Now we proceed on to define the notion of dual set code or an orthogonal set code.

**DEFINITION 2.9:** *Let V = {$X_1$, …, $X_n$} be a set code with associated set parity check matrix H = {$H_1$, …, $H_t$ | t < n}. The dual set code (or the orthogonal set code) $V^{\perp}$ = {$Y_i$ / $Y_i$. $X_i$ = (0) for all those $X_i \in V$; such that both $X_i$ and $Y_i$ are of same length in V}.*

Now it is important to see as in case of usual codes we can define orthogonality only between two vectors of same length alone.

We just illustrate this by an example.

***Example 2.18:*** Let V = {(0 0 0 0 0), (1 1 0 0 0), (1 0 1 0 0), (1 1 1 1 1), (0 0 0 0 0 0), (1 1 1 1 0 0), (0 1 1 1 1 0)} be a set code. The dual set code $V^{\perp}$ = {(0 0 0 0 0), (1 1 1 1 1), (1 1 1 0 0), (0 0 0 1 1), (1 1 1 1 0), (1 1 1 0 1), (0 1 1 0 0 0), (0 0 1 1 0 0), (0 1 0 1 0 0) and so on}.

For error detection in dual codes we need the set parity check matrix.

We can over come the use of set parity check matrix H by defining m-weight dual set codes.

**DEFINITION 2.10:** *Let V = {$X_1$, …, $X_n$} be a set code with a set parity check matrix H = {$H_1$, $H_2$, …, $H_t$ | t < n}. We call a proper subset S $\subseteq V^{\perp}$ for which the weight of each set code word in S to be only m, as the m-weight dual set code of V. For this set code for error detection we do not require the set parity check matrix.*

We illustrate this by some simple examples.

***Example 2.19:*** Let V = {(1 1 0 0), (0 0 0 0), (0 0 1 1), (0 1 0 1 0), (0 0 0 0 0), (0 0 0 1 1), (1 1 0 0 0), (0 0 1 1 0)} be a set code over {0, 1}. Now $V^{\perp}$ the dual set code of V of weight 2 is given by $V^{\perp}$ = {(0 0 0 0 0), (0 0 0 0), (0 0 1 1), (1 1 0 0 0), (0 0 1 1 0)}.

We give yet another example.

***Example 2.20:*** Let V = {(0 0 0 0 0), (1 1 1 0 0), (0 0 1 1 1), (1 1 1 0 0 0), (0 0 0 0 0 0), (0 0 0 1 1 1), (0 1 0 1 0 1)} be a set code.

The dual set code of weight 3 is given by $V^{\perp}$ = {(0 0 0 0 0), (0 0 0 0 0 0), (0 1 1 1 0), (0 1 1 0 1), (1 0 1 0 1), (1 0 1 0 1)}. The advantage of using m-weight set code or m-weight dual set code is that they are set codes for which error can be detected very easily.

Because these codes do not involve high abstract concepts, they can be used by non-mathematicians. Further these codes can be used when retransmission is not a problem as well as one wants to do work at a very fast phase so that once a error is detected retransmission is requested. These set codes are such that only error detection is possible. Error correction is not easy for all set codes, m-weight set codes are very simple class of set codes for which error detection is very easy. For in this case one need not even go to work with the set parity check matrix for error detection.

Yet another important use of these set codes is these can be used when one needs to use many different lengths of codes with varying number of message symbols and check symbols.

These set codes are best suited for cryptologists for they can easily mislead the intruder as well as each customer can have different length of code words. So it is not easy for the introducer to break the message for even gussing the very length of the set code word is very difficult.

Thus these set codes can find their use in cryptology or in places where extreme security is to be maintained or needed.

The only problem with these codes is error detection is possible but correction is not possible and in channels where retransmission is possible it is best suited. At a very short span of time the error detection is made and retransmission is requested.

Now we proceed on to define yet another new class of set codes which we call as semigroup codes.

**DEFINITION 2.11:** *Let $V = \{X_1, \ldots, X_n\}$ be a set code over the semigroup $S = \{0, 1\}$. If the set of codes of same length form semigroups under addition with identity i.e., monoid, then $V = \{S_1, \ldots, S_r \mid r < n\}$ is a semigroup code, here each $S_i$ is a semigroup having elements as set code words of same length; $1 \leq i \leq r$. The elements of V are called semigroup code words.*

We illustrate this by some examples.

***Example 2.21:*** Let S = {(0 0 0 0), (1 1 0 0), (0 0 1 0), (1 1 1 0), (0 0 0 0 0), (1 0 0 0 0), (0 0 0 0 1), (1 0 0 0 1), (0 0 1 0 0), (1 0 1 0 0), (1 0 1 0 1), (0 0 1 0 1)} be the semigroup code over the semigroup $Z_2$ = {0, 1}. Clearly S = {$S_1$, $S_2$} where $S_1$ = {(0 0 0 0), (1 1 0 0), (0 0 1 0), (1 1 1 0)} and $S_2$ = {(0 0 0 0 0), (1 0 0 0 0), (0 0 0 0 1), (1 0 0 0 1), (0 0 1 0 0), (1 0 1 0 0), (1 0 1 0 1), (0 0 1 0 1)} are semigroups with identity or monoids under addition.

***Example 2.22:*** Let V = {(1 1 1 1 1 1), (0 0 0 0 0 0), (1 1 0 0 0), (0 0 0 0 0), (0 0 1 1 1), (1 1 1 1 1), (1 1 1 1 1 1 1), (0 0 0 0 0 0 0), (0 0 1 1 0 0 1), (1 1 0 0 1 1 0)} be a semigroup code over $Z_2$ = {0, 1}. V = {$S_1$, $S_2$, $S_3$} where $S_1$ = {(1 1 1 1 1 1), (0 0 0 0 0 0)}, $S_2$ = {(1 1 0 0 0), (0 0 0 0 0), (0 0 1 1 1), (1 1 1 1 1)} and $S_3$ = {(1 1 1 1 1 1 1), (0 0 0 0 0 0 0), (0 0 1 1 0 0 1), (1 1 0 0 1 1 0)} are monoids under addition.

***Example 2.23:*** Let V = {(0 0 0 0 0), (1 1 1 1 1), (1 1 1 1 1 1 1), (0 0 0 0 0 0 0), (1 1 1 1 1 1 1 1 1 1), (0 0 0 0 0 0 0 0 0 0)} be a set code over $Z_2$ = {0, 1}. V is a repetition set code. Infact V = {$S_1$, $S_2$, $S_3$} where $S_1$ = {(0 0 0 0 0), (1 1 1 1 1)}, $S_2$ = {(0 0 0 0 0 0 0), (1 1 1 1 1 1 1)} and $S_3$ = {(1 1 1 1 1 1 1 1 1 1), (0 0 0 0 0 0 0 0 0 0)} are semigroups. Thus V is a semigroup code.

Now we prove the following.

**THEOREM 2.1:** *Every set repetition code V is a semigroup code.*

*Proof:* We know if C is any repetition code (ordinary), it has only two elements viz (0 0 … 0) and (1 1 … 1 1) i.e., a n tuple with all its entries zeros or ones. Thus if V = {$X_1$ … $X_n$} is a set repetition code it has only its elements to be zero tuples and one tuples of different lengths. So V can be written as {$S_1$ … $S_{n/2}$} where each $S_i$ is {(0 0 … 0), (1 1 … 1)}; $1 \leq i \leq n/2$. n is always even if V is a set repetition code. Now each $S_i$ is a semigroup $1 \leq i \leq n/2$. Thus V is a semigroup code.

Hence every set repetition code is a semigroup code.

We have the following result.

**THEOREM 2.2:** *Every semigroup code is a set code and a set code in general is not a semigroup code.*

*Proof*: Suppose V = {$X_1$, …, $X_n$} is a semigroup code then clearly V is a set code. Conversely suppose V = {$X_1$, …, $X_n$} is a set code it need not in general be a semigroup code. We prove this by the following example.

***Example 2.24:*** Let V = {(1 1 1 1 0), (0 0 0 0 0), (1 1 0 0 1), (1 0 0 0 1), (1 1 1 1 1 1 1), (0 0 0 0 0 0 0), (1 1 0 0 0 0 0 0), (0 0 0 0 0 1 1)} be a set code over the set {0, 1}. Clearly V is not a semigroup code for (1 1 0 0 1), (1 0 0 0 1) $\in$ V but (1 1 0 0 1) + (1 0 0 0 1) = (0 1 0 0 0) $\notin$ V. Hence the claim.

We can with a semigroup code V = {$X_1$, …, $X_n$}, associate a set parity check matrix H = {$H_1$, $H_2$, …, $H_r$ | r < n}. Here also the set parity check matrix only serves the purpose for error detection we do not use it to get all set code words for we do not use the group of code words. We use only a section of them which forms a semigroup.

Now as in case of set codes even in case of semigroup codes we can define m-weight semigroup code.

**DEFINITION 2.12:** *Let V = {$X_1$, …, $X_n$} = {$S_1$, …, $S_r$ | r < n} be a semigroup code if Hamming weight of each semigroup code is just only m then we call V to be a m-weight semigroup code.*

We illustrate this by some examples.

***Example 2.25:*** Let V = {(1 1 1 0 0 0), (0 0 0 0 0 0), (1 1 1 0 0 0 0), (0 0 0 0 0 0 0), (1 0 0 0 0 0 1 1), (0 0 0 0 0 0 0 0)} be a semigroup code. We see weight of each semigroup code word is 3 so V is a 3-weight semigroup code.

We have the following interesting result.

**THEOREM 2.3:** *Let V = {$X_1$, …, $X_n$} = {$S_1$ , $S_2$, …, $S_r$ | r < n} be a semigroup code. If V is a semigroup repetition code then V is not a m-weight semigroup code for any positive integer m.*

*Proof:* Given V is a repetition semigroup code so every semigroup $S_i$ in the semigroup code V has either weight 0 or weight $n_i$ if $n_i$ is the weight of $X_i$. So no two semigroup codes in V can have same weight. Thus a repetition semigroup code can never be a m-weight semigroup code.

Next we proceed to give the definition of semigroup parity check code.

**DEFINITION 2.13:** *Let V = {$X_1$, …, $X_n$} = {$S_1$, …, $S_r$ | r < n} be a semigroup code where each $S_i$ is a semigroup of the parity check binary ($n_i$, $n_i$ – 1) code, i = 1, 2, …, r. Then we call V to be a semigroup parity check code and the set parity check matrix H = {$H_1$, …, $H_r$} is such that*

$$H_i = \underbrace{\begin{pmatrix} 1 & 1 & \cdots & 1 \end{pmatrix}}_{n_i - times} \; ; i = 1, 2, \ldots, r.$$

We illustrate this by the following example.

***Example 2.26:*** Let V = {(1 1 1 1 1 1), (1 1 0 0 0 0), (0 0 1 1 1 1), (0 0 0 0 0 0), (1 1 0 0 0 1 1), (1 1 1 1 0 0 0), (0 0 0 0 0), (1 1 1 1 0), (1 1 0 0 0), (0 0 1 1 0)} = {$S_1$, $S_2$, $S_3$} where $S_1$ = {(0 0 0 0 0 0), (1 1 1 1 1 1), (1 1 0 0 0 0), (0 0 1 1 1 1)}, $S_2$ = {(0 0 0 0 0 0 0), (1 1 1 1 0 0 0), (0 0 1 1 0 1 1), (1 1 0 0 0 1 1)} and $S_3$ = {(0 0 0 0 0), (1 1 0 0 0), (1 1 1 1 0 0), (0 0 1 1 0)} is a semigroup parity check code. The set parity check matrix associated with V is {(1 1 1 1 1), (1 1 1 1 1 1), (1 1 1 1 1 1 1)} = {$H_1$, $H_2$, $H_3$}.

It is interesting to note that even a semigroup parity check code need not in general be a m-weight semigroup code. The above example is one such instance of a semigroup parity check code which is not a m-weight semigroup code.

As in case of set codes even in case of semigroup codes we can define orthogonal semigroup code.

**DEFINITION 2.14:** *Let $S = \{X_1, \ldots, X_n\} = \{S_1, \ldots, S_r \mid r < n\}$ be a semigroup code over {0, 1}. The dual or orthogonal semigroup code $S^{\perp}$ of S is defined by $S^{\perp} = \left\{S_1^{\perp}, \ldots, S_r^{\perp};\ r < n\right\}$ $S_i^{\perp} = \{s \mid s_i. s = 0$ for all $s_i \in S_i\}$; $1 \le i \le r$. The first fact to study about is will $S^{\perp}$ also be a semigroup code. Clearly $S^{\perp}$ is a set code. If x, y $\in S_i^{\perp}$ then $x.s_i = 0$ for all $s_i \in S_i$ and $y.s_i = 0$ for all $s_i \in S_i$. To prove $(x + y)\, s_i = 0$ i.e., $(x + y) \in S_i^{\perp}$. At this point one cannot always predict that the closure axiom will be satisfied by $S_i^{\perp}$, $1 \le i \le r$.*

We first atleast study some examples.

***Example 2.27:*** Let S = {(0 0 0 0), (0 0 1 0), (0 0 0 1), (0 0 1 1), (0 0 0 0 0), (0 0 0 1 1), (0 0 0 1 0), (0 0 0 0 1), (1 0 0 0 0), (1 0 0 1 1), (1 0 0 1 0), (1 0 0 0 1)} be a semigroup code. S = {$S_1$, $S_2$} where $S_1$ = {(0 0 0 0), (0 0 1 0), (0 0 0 1), (0 0 1 1)} and $S_2$ = {(0 0 0 0 0), (0 0 0 1 1), (0 0 0 1 0), (0 0 0 0 1), (1 0 0 0 0), (1 0 0 1 1), (1 0 0 1 0), (1 0 0 0 0 1)}.

$S_1^{\perp}$ = {(0 0 0 0), (1 0 0 0), (0 1 0 0), (1 1 0 0)} and
$S_2^{\perp}$ = {(0 0 0 0 0), (0 1 1 0 0), (0 1 0 0 0), (0 0 1 0 0)};

we see $S^{\perp} = \{S_1^{\perp}, S_2^{\perp}\}$ = {(0 0 0 0), (1 0 0 0), (0 1 0 0), (1 1 0 0), (0 0 0 0 0), (0 1 1 0 0), (0 1 0 0 0), (0 0 1 0 0)} is a semigroup code called the orthogonal dual semigroup code.

Now we proceed on to define the notion of semigroup cyclic code.

**DEFINITION 2.15:** *Let $S = \{X_1, \ldots, X_n\} = \{S_1, \ldots, S_r \mid r < n\}$ be a semigroup code if each of code words in $S_i$ are cyclic where $S_i$ is a semigroup under addition with identity, for each i; $1 \le i \le r$, then we call S to be a semigroup cyclic code.*

We now try to find some examples to show the class of semigroup cyclic codes is non empty.

**THEOREM 2.4:** *Let V = {$X_1$, …, $X_n$} = {$S_1$, …, $S_{n/2}$} be the semigroup repetition code, V is a semigroup cyclic code.*

*Proof:* Since any $X_i \in V$ is either of the form (0 0 … 0) or (1 1 … 1) so they happen to be cyclic. Hence the claim.

We leave it as an exercise for the reader to give examples of semigroup cyclic codes other than the semigroup repetition code.

Now we see this classes of semigroup codes form a sub class of set codes i.e., set codes happen to be the most generalized one with no algebraic operation on them. Now we give some approximate error correcting means to these new classes of codes. We know that for these set codes one can easily detect the error. Now how to correct error in them.

This is carried out by the following procedure.

Let V = {$X_1$, …, $X_n$} be a set code with set parity check matrix H = {$H_1$, …, $H_r$ | r < n}. Suppose $Y_i$ happens to be the received message by finding $H_iY_i^t = (0)$; no error, only $Y_i$ is the received message. If $H_iY_i^t \neq (0)$ then $Y_i$ is not the real message some error has occurred during transmission. To find approximately the correct message. Suppose $Y_i$ is of length $n_i$ then find the Hamming distance between all set code words $X_i$ in V of length $n_i$ and choose that $X_i$ which has least distance from $Y_i$ as the approximately correct message. If more than one set code has same minimal value then take that set code word which has least value between the message symbol weights. For instance if $Y_i$ = (1 1 1 0 0 1 1) is the received set code word and $Y_i \notin V$. This $Y_i$ has first four coordinates to be the message symbols and the last three coordinates are the check symbols. Suppose $X_i$ = (1 1 0 1 1 0 0) and

$X_j$ = (0 1 1 1 1 0 1) $\in$ V;

$d(Y_i, X_i) = 5$
$d(Y_i, X_j) = 4;$

we choose $X_j$ to be the approximately correct word. Suppose $X_t$ = (1 1 1 1 1 0 0) $\in$V and $d(Y_i, X_t) = 4$.

Now $X_j$ and $X_t$ are at the same distance from $Y_i$, which is to be choosen $X_t$ or $X_j$ since the difference in the message symbols between $Y_i$ and $X_i$ is 2 where as the difference between the message symbols between $Y_i$ and $X_t$ is one we choose $X_t$ to $Y_i$. This is the way the choice is made without ambiguity as same message symbol cannot be associated with two distinct check symbols. Thus the approximately correct message is determined.

This method is used in set codes which are not m-weight set codes or m-weight semigroup codes.

Next we proceed on to define the new notion of group codes over group under addition.

**DEFINITION 2.16:** *Let $V = \{X_1, \ldots, X_n\}$ be a set code if $V = \{G_1, \ldots, G_k \mid k < n\}$ where each $G_i$ is a collection of code words which forms a group under addition then we call V to be a group code over the group $Z_2 = \{0,1\}$.*

We illustrate this situation by the following examples.

***Example 2.28:*** Let V = {(0 0 0 0 0), (1 1 1 1 1), (0 0 0 0 0 0), (1 1 1 0 0 0), (1 1 1 1 1 1), (0 0 0 1 1 1), (1 1 1 1 1 1 1 1), (0 0 0 0 0 0 0 0), (1 1 0 0 1 1 0 0), (0 0 1 1 0 0 1 1)} = $\{G_1, G_2, G_3\}$ where $G_1$ = {(0 0 0 0 0), (1 1 1 1 1)}, $G_2$ = {(0 0 0 0 0 0), (0 0 0 1 1 1), (1 1 1 0 0 0)} and $G_3$ = {(1 1 1 1 1 1 1 1), (0 0 0 0 0 0 0 0), (1 1 0 0 1 1 0 0), (0 0 1 1 0 0 1 1)} where $G_1$, $G_2$ and $G_3$ are group under addition. Clearly $G_1 \subseteq Z_2^5$; $G_2 \subseteq Z_2^6$ and $G_3 \subseteq Z_2^8$; i.e., $G_i$ are subgroups of the groups $Z_2^n$ (n = 5, 6, 8). Thus we can verify whether a received word X is in V or not using the set parity check matrix.

We give yet another example.

***Example 2.29:*** Let V = {(1 1 1 1 1 1 1), (0 0 0 0 0 0 0), (1 1 0 0 0 0), (1 0 0 0 0 0), (0 1 0 0 0 0), (0 0 0 0 0 0), (1 1 0 0 1), (0 0 1 1 0), (1 1 1 1 1), (0 0 0 0 0)} be a group code over $Z_2$.

Clearly V = {$G_1$, $G_2$, $G_3$} where $G_1$ = {(1 1 1 1 1 1 1), (0 0 0 0 0 0 0)} $\subseteq Z_2^7$. $G_2$ = {(1 1 0 0 0 0), (1 0 0 0 0 0), (0 1 0 0 0 0), (0 0 0 0 0 0)} $\subseteq Z_2^6$ and $G_3$ = {(1 1 1 1 1), (0 0 0 0 0), (1 1 0 0 1), (0 0 1 1 0)} $\subseteq Z_2^5$ are groups and the respective subgroups of $Z_2^7$, $Z_2^6$ and $Z_2^5$ respectively.

This subgroup property helps the user to adopt coset leader method to correct the errors. However the errors are detected using set parity check matrices.

All group codes are semigroup codes and semigroup codes are set codes.

But in general all the set codes need not be a semigroup code or a group code. In view of this we prove the following theorem.

**THEOREM 2.5:** *Let V be a set code V in general need not be a group code.*

*Proof:* We prove this by a counter example. Take V = {(0 0 0 0 0), (1 1 0 0 0), (0 0 1 0 0), (1 1 1 1 1 1), (0 0 0 0 0 0), (0 1 0 1 0 1), (1 1 1 0 0 0 0), (0 0 1 1 1 1 0), (0 0 0 0 0 0 0), (0 1 0 1 0 1 0), (1 1 0 0 1 1 0)} to be a set code over $Z_2$ = {0, 1}.

Now take X = (1 1 0 0 0) and Y = (0 0 1 0 0) $\in$ V. Clearly X + Y = (1 1 1 0 0) is not in V so V is not even closed under addition hence V cannot be a group code over $Z_2$ = {0, 1}.

Next we prove every repetition set code is a group code.

**THEOREM 2.6:** *Let V be a set repetition code then V is a group code.*

*Proof:* Given V = {$X_1$, …, $X_n$} is a set repetition code. So if $X_i$ $\in$ V then

$$X_i = \underbrace{\begin{pmatrix} 0 & 0 & \cdots & 0 \end{pmatrix}}_{n_i - \text{tuples}}$$

or

$$X_j = \underbrace{\begin{pmatrix} 1 & 1 & \cdots & 1 \end{pmatrix}}_{n_i - \text{tuples}}.$$

Thus V = $\{G_1, \ldots, G_{n/2}\}$ and the order of a set repetition code is always even. Each $G_i$ is a group under addition; $1 \le i \le n/2$. $G_i = \{(0\ 0 \ldots 0), (1\ 1 \ldots 1)\}$ is a group. Thus V is a group code. Every set repetition code is a group code.

We give some examples of these codes.

***Example 2.30:*** Let V = {(1 1 1 1 1 1), (0 0 0 0 0 0), (1 1 1 1 1), (0 0 0 0 0), (1 1 1 1 1 1 1 1 1), (0 0 0 0 0 0 0 0 0), (1 1 1 1 1 1 1), (0 0 0 0 0 0 0)} be a set code. We see order of V is eight and V = $\{G_1, G_2, G_3, G_4\}$ where $G_1$ = {(0 0 0 0 0 0), (1 1 1 1 1 1)}, $G_2$ = {(1 1 1 1 1) (0 0 0 0 0)}, $G_3$ = {(0 0 0 0 0 0 0 0 0), (1 1 1 1 1 1 1 1 1)} and $G_4$ = {(1 1 1 1 1 1 1), (0 0 0 0 0 0 0)} are groups. So V is a group code.

We now define these codes as repetition code.

**DEFINITION 2.17:** *Let V = $\{X_1, \ldots, X_n\} = \{G_1, \ldots, G_{n/2}\}$ where V is a set repetition code. Clearly V is a group code. We call V to be the group repetition code.*

The following facts are interesting about these group repetition code.

1. Every set repetition code is a group repetition code.
2. Every group repetition code is of only even order.
3. The error detection and correction is very easily carried out.

If $X_i \in V$ is a length $n_i$, if the number of ones in $X_i$ is less than $n_i/2$ then we accept $X_i$ = (0 0 … 0 0).

If the number of ones in $X_i$ is greater than $n_i/2$ then we accept $X_i = (1\ 1\ \ldots\ 1)$. Thus the easy way for both error correction and error detection is possible.

Now we proceed on to describe group parity check code and group Hamming code.

**DEFINITION 2.18:** *Let $V = \{X_1, \ldots, X_n\}$ be a group code if $H = \{H_1, \ldots, H_r \mid r < n\}$ be the associated set parity check matrix where each $H_i$ is of the form*

$$\underbrace{(1 \quad 1 \quad \cdots \quad 1)}_{n_i - times}$$

*$1 < i < r$. i.e., $V = \{G_1, \ldots, G_r \mid r < n\}$ and each $G_i$ is a set of code words of same length $n_i$ and forms a group under addition modulo 2; then V is defined as the group parity check code.*

We now illustrate this by the following example.

***Example 2.31:*** Let $V = \{X_1, \ldots, X_n\}$ = {(0 0 0 0 0 0), (1 1 0 0 0 0), (0 0 1 1 1 1), (1 1 1 1 1 1), (1 1 1 1 0 0 0), (0 0 0 0 0 0 0), (1 1 1 1 0), (1 1 0 0 0), (0 0 1 1 0), (0 0 0 0 0)} = $\{G_1, G_2, G_3\}$ where

$G_1$ = {(0 0 0 0 0 0), (1 1 0 0 0 0), (0 0 1 1 1 1)},

$G_2$ = {(0 0 0 0 0 0 0), (1 1 1 1 0 0 0)}

and

$G_3$ = {(1 1 0 0 0), (0 0 1 1 0), (0 0 0 0 0)}

are groups. So V is a group code which is a parity check group code.

The associated set parity check matrix is $H = \{H_1, H_2, H_3\}$ where

$H_1$ = (1 1 1 1 1 1),

$H_2$ = (1 1 1 1 1)

and

$H_3$ = (1 1 1 1 1 1 1).

We give yet another example of parity check group code.

***Example 2.32:*** Let V = {(0 0 0 0), (1 1 0 0), (0 0 1 1), (1 1 1 1), (1 1 1 1 0 0 0 0), (1 1 0 0 0 0 0 0), (0 0 1 1 0 0 0 0), (0 0 0 0 0 0 0 0), (1 1 1 0 0 0 1), (0 0 1 0 0 0 1), (1 1 0 0 0 0 0), (0 0 0 0 0 0 0)} = {$G_1$, $G_2$, $G_3$} where

$$G_1 = \{(0\ 0\ 0\ 0), (1\ 1\ 0\ 0), (0\ 0\ 1\ 1), (1\ 1\ 1\ 1)\},$$

$$G_2 = \{(1\ 1\ 1\ 1\ 0\ 0\ 0\ 0), (1\ 1\ 0\ 0\ 0\ 0\ 0\ 0), (0\ 0\ 1\ 1\ 0\ 0\ 0\ 0), (0\ 0\ 0\ 0\ 0\ 0\ 0\ 0)\}$$

and

$$G_3 = \{(1\ 1\ 1\ 0\ 0\ 0\ 1), (0\ 0\ 1\ 0\ 0\ 0\ 1), (1\ 1\ 0\ 0\ 0\ 0\ 0)\}$$

are groups and V is a parity check group code with set parity check matrix H = {$H_1$, $H_2$, $H_3$} where

$$H_1 = (1\ 1\ 1\ 1),$$

$$H_2 = (1\ 1\ 1\ 1\ 1\ 1\ 1\ 1)$$

and

$$H_3 = (1\ 1\ 1\ 1\ 1\ 1\ 1).$$

Now we proceed on to define the new notion of binary Hamming group code.

**DEFINITION 2.19:** *Let V = {$X_1$, …, $X_n$} = {$G_1$, …, $G_r$ | r < n} be a group code. If each $G_i$ is a Hamming code with parity check matrix $H_i$, i = 1, 2, …, r i.e., of the set parity check matrix H = {$H_1$, $H_2$, …, $H_r$}, then we call V to be a group Hamming code or Hamming group code.*

We now illustrate this situation by few examples.

***Example 2.33:*** Let V = {(0 0 0 0 0 0 0), (0 0 0 1 1 1 1), (0 1 1 0 0 1 1), (0 1 1 1 1 0 0), (0 0 0 0 0 0 0 0), (1 1 1 1 1 1 1 1), (0 0 0 1 1 1 1 0), (1 1 1 0 0 0 0 1)} = {$G_1$, $G_2$} where

$G_1$ = {(0 0 0 0 0 0 0), (0 0 0 1 1 1 1), (0 1 1 0 0 1 1), (0 1 1 1 1 0 0)}

and

$G_2$ = {(0 0 0 0 0 0 0 0), (1 1 1 1 1 1 1 1), (0 0 0 1 1 1 1 0), (1 1 1 0 0 0 0 1)}

is a Hamming group code with set parity check matrix H = ($H_1$, $H_2$) where

$$H_1 = \begin{pmatrix} 0 & 0 & 0 & 1 & 1 & 1 & 1 \\ 0 & 1 & 1 & 0 & 0 & 1 & 1 \\ 1 & 0 & 1 & 0 & 1 & 0 & 1 \end{pmatrix}$$

and

$$H_2 = \begin{pmatrix} 1 & 1 & 1 & 1 & 1 & 1 & 1 & 1 \\ 0 & 0 & 0 & 1 & 1 & 1 & 1 & 0 \\ 0 & 1 & 1 & 0 & 0 & 1 & 1 & 0 \\ 1 & 0 & 1 & 0 & 1 & 0 & 1 & 0 \end{pmatrix}.$$

$H_2$ is the parity check matrix of the (8, 4, 4) extended Hamming code.

We give yet another example of a Hamming group code before we proceed on to define group cyclic codes or cyclic group codes.

***Example 2.34:*** Let V = {(0 0 0 0 0 0 0 0 0 0 0 0 0 0 0), (1 0 0 0 1 0 0 1 1 0 1 0 1 1 1), (0 1 0 0 1 1 0 1 0 1 1 1 1 0 0), (1 1 0 0 0 0 1 0 0 1 1 0 1 0 1 1), (0 0 0 0 0 0 0), (1 0 0 1 1 0 1), (0 1 0 1 0 1 1), (1 1 0 0 1 1 0), (0 0 1 0 1 1 1), (0 1 1 1 1 0 0), (1 1 1 0 0 0 1), (1 0 1 1 0 1 0)} be a Hamming group code with set parity check matrix H = {$H_1$, $H_2$} where

$$H_1 = \begin{pmatrix} 1 & 0 & 0 & 0 & 1 & 0 & 0 & 1 & 1 & 0 & 1 & 0 & 1 & 1 & 1 \\ 0 & 1 & 0 & 0 & 1 & 1 & 0 & 1 & 0 & 1 & 1 & 1 & 1 & 0 & 0 \\ 0 & 0 & 1 & 0 & 0 & 1 & 1 & 0 & 1 & 0 & 1 & 1 & 1 & 1 & 0 \\ 0 & 0 & 0 & 1 & 0 & 0 & 1 & 1 & 0 & 1 & 0 & 1 & 1 & 1 & 1 \end{pmatrix}$$

and

$$H_2 = \begin{pmatrix} 1 & 0 & 0 & 1 & 1 & 0 & 1 \\ 0 & 1 & 0 & 1 & 0 & 1 & 1 \\ 0 & 0 & 1 & 0 & 1 & 1 & 1 \end{pmatrix}.$$

Thus V = $\{G_1, G_2\}$ where $G_1$ and $G_2$ are groups given by $G_1$ = {(0 0 0 0 0 0 0 0 0 0 0 0 0 0), (1 0 0 0 1 0 0 1 1 0 1 0 1 1 1), (0 1 0 0 1 1 0 1 0 1 1 1 1 0 0), (1 1 0 0 0 1 0 0 1 1 0 1 0 1 1)} and $G_2$ = {(0 0 0 0 0 0 0), (1 0 0 1 1 0 1), (0 1 0 1 0 1 1), (1 1 0 0 1 1 0), (0 0 1 0 1 1 1), (0 1 1 1 1 0 0), (1 1 1 0 0 0 1), (1 0 1 1 0 1 0)}.

Now having seen few examples of Hamming group code we now proceed on to define cyclic group code or group cyclic code.

**DEFINITION 2.20:** *Let V = $\{X_1, \ldots, X_n\} = \{G_1, G_2, \ldots, G_t \mid t < n\}$ be a group code. If each of the $G_i$ is a cyclic code which forms a subgroup of some cyclic code $C_i$; i = 1, 2, …, t with set parity check matrix H = $\{H_1, \ldots, H_t\}$.Then we call V to be a group cyclic code or cyclic group code.*

It is important to mention that we need not take the complete cyclic code $C_i$ got using the parity check matrix $H_i$; i = 1, 2, …, t.

We illustrate this situation by a few examples before we proceed on to give more properties about group codes.

***Example 2.35:*** Let V = {(0 0 0 0 0 0), (1 1 1 1 1 1), (0 0 0 0 0 0 0 0), (1 1 1 1 1 1 1 1), (0 0 0 0 0), (1 1 1 1 1)} be a cyclic group code with V = $\{G_1, G_2, G_3\}$ where $G_1$ = {(0 0 0 0 0 0), (1 1 1 1 1 1)}, $G_2$ = {(0 0 0 0 0 0 0 0), (1 1 1 1 1 1 1 1)} and $G_3$ = {(0 0 0 0 0), (1 1 1 1 1)} are group codes.

One of the important features about this code in this example leads to the following result.

**THEOREM 2.7:** *Every group repetition code $V = \{X_1, \ldots, X_n\}$ is a cyclic group code with even number of code words in V. However the converse i.e., every cyclic group code is a repetition code is not true.*

*Proof:* Consider the repetition group code $V = \{X_1, \ldots, X_n\}$; clearly n is even for if $X_i = (0, \ldots, 0)$ is of length $n_i$ then their exists a $X_j$ in V such that $X_j = (1\ 1, \ldots, 1)$ which is of length $n_i$ and $\{X_i, X_j\}$ forms a group. Thus $V = \{G_1, G_2, \ldots, G_{n/2}\}$, where each $G_i$ is a repetition code and order of each $G_i$ is two. Clearly each $G_i$ is a cyclic code as every repetition code is a cyclic code. Hence V is a cyclic group code with even number of elements in it.

The proof of the converse is established using a counter example. Consider the cyclic group code V = {(0 0 0 0 0 0), (0 0 1 0 0 1), (0 1 0 0 1 0), (0 1 1 0 1 1), (1 0 0 1 0 0), (1 0 1 1 0 1), (1 1 0 1 1 0), (1 1 1 1 1 1), (0 0 0 0 0 0 0), (0 0 1 0 1 1 1), (0 1 0 1 1 1 0), (1 0 1 1 1 0 0), (0 1 1 1 0 0 1), (1 0 0 1 0 1 1), (1 1 0 0 1 0 1), (1 1 1 0 0 1 0)} = $\{G_1, G_2\}$ where $G_1$ = {(0 0 0 0 0 0) (0 0 1 0 0 1), (0 1 0 0 1 0), (0 1 1 0 1 1), (1 0 0 1 0 0), (1 0 1 1 0 1), (1 1 0 1 1 0), (1 1 1 1 1 1)} and $G_2$ = {(0 0 0 0 0 0 0) (0 0 1 0 1 1 1) (0 1 0 1 1 1 0), (1 0 1 1 1 0 0), (0 1 1 1 0 0 1), (1 0 0 1 0 1 1), (1 1 0 0 1 0 1), (1 1 1 0 0 1 0)} are the group codes. It is easily observed both the codes are cyclic so V is a group cyclic code however both the codes $G_1$ and $G_2$ are not repetition codes. Hence the claim.

Now we proceed on to define the dual (orthogonal) group code.

**DEFINITION 2.21:** *Let $V = \{X_1, \ldots, X_n\} = \{G_1, \ldots, G_t \mid t < n\}$ be a group code. The dual (or orthogonal) code $V^{\perp}$ of V is defined by $G_i^{\perp} = \{Y \mid Y. X_i = 0$ for every $X_i \in G_i\}$, $1 \leq i \leq t$; where $V^{\perp} = \{G_1^{\perp}, G_2^{\perp}, \ldots, G_t^{\perp}\}$*

A natural question would be will $V^{\perp}$ be a group code? We leave the answer to this question as exercise. We now illustrate this by some simple examples.

***Example 2.36:*** Let V = {(0 0 0 0 0 0), (0 0 0 0), (1 0 1 1), (0 1 0 1), (1 1 1 0), (1 1 1 1 1 1)} = {$G_1$, $G_2$} where $G_1$ = {(0 0 0 0), (1 0 1 1), (0 1 0 1), (1 1 1 0)} and $G_2$ = {(0 0 0 0 0 0) (1 1 1 1 1 1)} be group codes so that V is a group code. To find $V^{\perp}$ = {(0 0 0 0 0 0), (1 1 1 1 1 1), (1 1 0 0 0 0), (0 0 1 1 0 0), (0 0 0 0 1 1), (1 0 1 0 0 0), (1 0 0 1 0 0), (1 0 0 0 1 0), (0 1 1 0 0 0), (0 1 0 1 0 0), (0 1 0 0 1 0), (0 1 0 0 0 1), (0 0 1 0 1 0), (0 0 1 0 0 1), (0 0 0 1 0 1), (0 0 0 1 1 0), (1 1 1 1 0 0), (0 1 1 1 1 0), (0 0 1 1 1 1), (1 1 0 1 1 0), (1 1 0 0 1 1), (0 1 1 0 1 1), (1 0 1 1 0 1), (0 1 1 1 0 1), (0 1 0 1 1 1), (1 1 1 0 1 0), (1 1 1 0 0 1), (0 1 0 1 1 1), (1 0 0 1 1 1), (1 0 1 1 1 0), (1 0 1 0 1 1)}, {(0 0 0 0), (1 0 1 0), (1 1 0 1), (0 1 1 1)} = ($G_1^{\perp}$, $G_2^{\perp}$).

We see the dual group code of a group code is also a group code.

Now we describe how we can do the error correction and error detection in group codes. Suppose V = {$X_1$, …, $X_n$} = {$G_1$, …, $G_r$ | r < n} be a group code with associated set parity check matrix H = {$H_1$, …, $H_r$}. Suppose Y is a received message knowing the length of Y we find the appropriate parity check matrix $H_i$ of H and find $H_iY^t$; if $H_iY^t$ = (0) then the received message is the correct one. If $H_iY^t \neq 0$ then we say some error has occured during transmission.

Now how to correct the error once we have detected it. The error correction can be carried out in two ways.

**Method I:** In this method we obtain only approximately the close word to the sent word, once error has been detected. If y $\notin$ V, by monitoring the length of y we can say the code word is only from the group code say $G_i$. Now y $\notin G_i$ we want to find a x $\in G_i$ such that the Hamming distance between x and y is the least. We calculate d (x, y) for every x in $G_i$. If more than one x has the minimal distance with y then we observe the Hamming distance between those x and y only taking into account the message symbols and accept the least of them. The least distance code word in $G_i \subset$ V is taken as the approximately correct message.

Now we describe this by the following example so that the reader can know how to detect and correct errors in group codes.

***Example 2.37:*** Let V = $\{X_1, \ldots, X_n\} = \{G_1, \ldots, G_r \mid r < n\}$ be a group code. Suppose H = $\{H_1, \ldots, H_r\}$ be the set parity check matrix associated with V. Let y = (1 1 1 0 0 1) be the received word. Suppose $y \notin G_i$ = {(0 0 0 0 0 0) (0 0 1 0 0 1), (0 1 0 0 1 0), (0 1 1 0 1 1), (1 0 0 1 0 0), (1 0 1 1 0 1), (1 1 0 1 1 0), (1 1 1 1 1 1)} with the associated parity check matrix $H_i \in H$ where

$$H_i = \begin{pmatrix} 1 & 0 & 0 & 1 & 0 & 0 \\ 0 & 1 & 0 & 0 & 1 & 0 \\ 0 & 0 & 1 & 0 & 0 & 1 \end{pmatrix}.$$

Clearly $G_i$ are codes with 3 message symbols and 3 check symbols and length of the code words in $G_i$ is 6. We first find

$$H_i\, y^T = \begin{pmatrix} 1 & 0 & 0 & 1 & 0 & 0 \\ 0 & 1 & 0 & 0 & 1 & 0 \\ 0 & 0 & 1 & 0 & 0 & 1 \end{pmatrix} \begin{pmatrix} 1 \\ 1 \\ 1 \\ 0 \\ 0 \\ 1 \end{pmatrix} = \begin{pmatrix} 1 & 1 & 0 \end{pmatrix} \neq 0 .$$

So we confirm an error has occurred as $H_i\, y^T \neq 0$.

Now how to find approximately a code word close to y. This we do by using method I i.e., we find the Hamming distance between y and every $x \in G_i$.

d ((0 0 1 0 0 1), (1 1 1 0 0 1)) = 2 --- (1)
d ((0 1 0 0 1 0), (1 1 1 0 0 1)) = 4
d ((0 1 1 0 1 1), (1 1 1 0 0 1)) = 2 --- (2)
d ((1 0 0 1 0 0), (1 1 1 0 0 1)) = 4
d ((1 0 1 1 0 1), (1 1 1 0 0 1)) = 2 --- (3)

d ((1 1 0 1 1 0), (1 1 1 0 0 1)) = 4
d ((1 1 1 1 1 1), (1 1 1 0 0 1)) = 2 --- (4)

Now we find the Hamming distance between the sent messages and the codes given in (1) (2) (3) and (4). We see the difference among message symbols in (1) is 2 where as in (2) it is one, in (3) it is also one but the difference between the message symbols given by 4 is zero. So we accept (1 1 1 1 1 1) to be the approximately the correct message as d (message symbols of (1 1 1 1 1 1) and message symbols of (1 1 1 0 0 1)) is 0 as well as the same code also gives the minimal number of differences.

Now we describe the method II. The way of finding correct message or the method of correcting the error that has occurred during the transmission.

**Method II** This method is as in case of the coset leaders described in chapter one of this book.

Suppose V = $\{X_1, \ldots, X_n\}$ is the group code, i.e., V = $\{X_1, \ldots, X_n\} = \{G_1, \ldots, G_r \mid r < n\}$. Suppose y is the received message from knowing the length of the code word we using the parity check matrix from the set parity check matrix H = $\{H_1, \ldots, H_r\}$ find whether the received word has an error or not by finding the syndrome of y. If S(y) = 0 then we accept y to be the correct word if S(y) $\neq$ 0 we confirm that some error has occurred during transmission. How to correct the error which has occurred during transmission.

If

$$H_i = \begin{pmatrix} 1 & 0 & 0 & 1 & 0 & 0 \\ 0 & 1 & 0 & 0 & 1 & 0 \\ 0 & 0 & 1 & 0 & 0 & 1 \end{pmatrix}$$

is the associated parity check matrix of the received code word y and $H_i\, y^T \neq (0)$. Now we use the coset leader method. Suppose y = (1 1 1 0 0 1) is the received word.

The corresponding coset table is

| Message words | 0 0 0 | 0 1 0 | 0 0 1 | 1 0 0 |
|---|---|---|---|---|
| Code words | 0 0 0 0 0 0 | 0 1 0 0 1 0 | 0 0 1 0 0 1 | 1 0 0 1 0 0 |
| Other cosets | 1 0 0 0 0 0<br>0 1 0 0 0 0<br>0 0 1 0 0 0<br>1 1 0 0 0 0<br>0 1 1 0 0 0<br>1 0 0 0 0 1<br>1 1 1 0 0 0 | 1 1 0 0 1 0<br>0 0 0 0 1 0<br>0 1 1 0 1 0<br>1 0 0 0 1 0<br>0 0 1 0 1 0<br>1 1 0 0 1 1<br>1 0 1 0 1 0 | 1 0 1 0 0 1<br>0 1 1 0 0 1<br>0 0 0 0 0 1<br>1 1 1 0 0 1<br>0 1 0 0 0 1<br>1 0 1 1 0 1<br>1 1 0 0 0 1 | 0 0 0 1 0 0<br>1 1 0 1 0 0<br>1 0 1 1 0 0<br>0 1 0 1 0 0<br>1 1 1 1 0 0<br>0 0 0 1 0 1<br>0 1 1 1 0 0 |
| Message words | 1 1 0 | 0 1 1 | 1 0 1 | 1 1 1 |
| Code words | 1 1 0 1 1 0 | 0 1 1 0 1 1 | 1 0 1 1 0 1 | 1 1 1 1 1 1 |
| Other cosets | 0 1 0 1 1 0<br>1 0 0 1 1 0<br>1 1 1 0 1 1<br>0 0 0 1 1 0<br>1 0 1 1 1 0<br>0 1 0 1 1 1<br>0 0 1 1 1 0 | 1 1 1 0 1 1<br>0 0 1 0 1 1<br>0 1 0 0 1 1<br>1 0 1 0 1 1<br>0 0 0 0 1 1<br>1 1 1 0 1 0<br>1 0 0 0 1 1 | 0 0 1 1 0 1<br>1 1 1 1 0 1<br>1 0 0 1 0 1<br>0 1 1 1 0 1<br>1 1 0 1 0 1<br>0 0 1 1 0 0<br>0 1 0 1 0 1 | 0 1 1 1 1 1<br>1 0 1 1 1 1<br>1 1 0 1 1 1<br>0 0 1 1 1 1<br>1 0 0 1 1 1<br>0 1 1 1 1 0<br>0 0 0 1 1 1 |

Now we see (1 1 1 0 0 1) has occurred with the coset leader 1 1 0 0 0 0 so the sent code word is (0 0 1 0 0 1) i.e., 1 1 1 0 0 1 + 1 1 0 0 0 0 = 0 0 1 0 0 1.

It is very important to note that the hamming difference between the received code and the corrected message happens to give a minimal difference for d(0 0 1 0 0 1), (1 1 1 0 0 1)) is 2. But however the difference between the message symbols are very large. It is appropriate to use any method as per the wishes of the designer who use these codes.

It is pertinent to mention here that as in case of usual codes we can adopt in group codes the method of coset leader for error correction. The main feature is that in general codes we cannot use codes of different length from the same channel. All codes would be of same length n having only the same number of check symbols say (n – k) and same number of message symbols k such that the length of the code is n – k + k = n. But in group codes we have codes of different length and also we have different number of message symbols and check symbols depending from which group code $G_i$ we are sending the message. The advantage of these codes is that in the same machine we can transform codes of different lengths different sets of messages with different sets of check symbols. With the advent of computer it is not very difficult to process the correctness of the received message!

Chapter Three

# SET BICODES AND THEIR GENERALIZATIONS

In this chapter we for the first time define the notion of set bicodes, group bicodes and semigroup bicodes and generalize them. We give or indicate the applications of these codes then and there. This chapter has two sections.

## 3.1 Set bicodes and their properties

In this section we proceed on to define the new notion of set bicodes and enumerate a few of their properties.

**DEFINITION 3.1.1:** *Let* $C = C_1 \cup C_2$

$$= \left\{\left(x_1^1 \ldots x_{r_1}^1\right), \left(x_1^1 \ldots x_{r_2}^1\right), \ldots, \left(x_1^1 \ldots x_{r_{n_1}}^1\right)\right\}$$

$$\cup \left\{\left(x_1^2 \ldots x_{s_1}^2\right)\left(x_1^2 \ldots x_{s_2}^2\right), \ldots, \left(x_1^2 \ldots x_{s_{n_2}}^2\right)\right\}$$

*be a biset of* $\left(r_1,...,r_{n_1}\right) \cup \left(s_1,...,s_{n_2}\right)$ *bituples with entries from the set S = {0, 1} where each* $\left(r_{i_1}, s_{j_2}\right)$ *tuple* $\left(x_1^1 \ldots x_{r_{i_1}}^1\right) \cup \left(x_1^2 \ldots x_{s_{j_2}}^2\right)$ *has some* $\left(k_{i_1}^1, k_{j_2}^2\right)$ *message symbols and* $\left(r_{i_1} - k_{i_1}^1, s_{j_2} - k_{j_2}^2\right)$ *bicheck symbols* $1 \le i_1 < r_i$, $1 \le j_2 < s_{j_2}$, $1 \le i_1 \le n_1$ *and* $1 \le j_2 \le n_2$.

*We call* $C = C_1 \cup C_2$ *the set bicode if C is a set bivector space over the set S = {0, 1}.*

The following observations will throw more light about these set bicodes.

(1) $r_{i_1} = r_{j_1}, s_{i_2} = s_{j_2}$ even if $i_t \neq j_t$; $1 \le i_t, j_t \le n_t$, t = 1, 2.

(2) $k_{i_1}^1 = k_{j_1}^2; k_{i_2}^1 = k_{j_2}^2$ even if $i_t \neq j_t$; $1 \le i_t, j_t \le n_t$, t = 1, 2.

(3) At least some $r_{i_1} \neq r_{t_1}$ when $i_1 \neq t_1$; $1 \le i_1, t_1 \le n_1$, $s_{j_2} = s_{t_2}$ where $j_2 \neq t_2$, $1 \le j_2, t_2 \le n_2$.

(4) $C_1 \not\subseteq C_2$, $C_2 \not\subseteq C_1$ i.e., $C_1$ and $C_2$ are distinct, however $C_1 \cap C_2$ need not be empty.

We illustrate this situation by the following example.

***Example 3.1.1:*** Let $C = C_1 \cup C_2$ = {(1 1 1), (0 0 0), (1 1 0), (1 0 0), (0 1 0), (0 0 0 0 0), (1 1 1 0 0), (1 1 0 0 0), (0 0 0 0 1), (0 0 0 0 0), (1 1 1 0 0 0), (0 1 1 1 0 0), (0 0 1 1 1 0), (1 0 1 1 0 0), (1 0 0 0 1 1)} ∪ {(0 1 1 0), (0 0 0 0), (1 0 1 0), (1 1 1 1), (0 0 0 0 0 0 0), (0 0 0 1 0 0 0), (1 1 1 0 0 0 1), (1 1 0 0 1 1 0), (1 0 1 0 1 0 1), (0 0 0 0 0), (1 1 1 1 1)}. C is a set bicode over the set S = {0, 1}.

The main advantage of these set bicodes is that they have two set codes of varying lengths. We call the elements of the set bicodes as set bicode words. Any $x \in C = C_1 \cup C_2$ is denoted by

$$x = \left(x_1^1, \ldots, x_{r_1}^1\right) \cup \left(x_1^2, \ldots, x_{s_2}^2\right); \left(x_1^1, x_2^1, \ldots, x_{r_1}^1\right) \in C_1$$

and

$$\left(x_1^2, x_2^2, \ldots, x_{s_2}^2\right) \in C_2.$$

Another advantage of these codes over the binary codes or binary bicodes in which we have code words of same length in $C_1$ and like wise the code words in $C_2$ are of same length though the length of the code words in $C_1$ and $C_2$ need not be of the same length in general where as the set codes and set bicodes are such that the set code words in them can be of different lengths. So in a set bicode we can send two messages of different lengths which is not possible in case of ordinary bicodes. Further we do not demand the set bicodes to form a bigroup or subbigroup or a subbispace of any finite dimensional vector bispace. They are just a collection of $\left(r_1, r_2, \ldots, r_{n_1}\right) \cup \left(s_1, s_2, \ldots, s_{n_2}\right)$-bituples with no proper algebraic structure defined on it.

***Example 3.1.2:*** Let C = $C_1 \cup C_2$ = {(1 1 1 1), (0 0 0 0), (1 1 1 1 1), (0 0 0 0 0), (1 1 0 0 0), (0 0 0 1 1), (1 1 1 1 1 1 1), (0 0 0 0 0 0 0), (1 1 0 0 0 0 0), (1 1 1 0), (1 1 0 0)} $\cup$ {(1 1 1), (0 0 0), (1 1 0), (0 1 0), (1 1 1 1 1 1), (0 0 0 0 0 0), (1 1 0 0 0 0), (0 0 1 1 0 0), (0 0 1 0 1 0), (1 1 1 1 1 1 1 1), (0 0 0 0 0 0 0 0), (1 1 0 0 1 1 0 0), (1 0 1 1 0 1 0 1)} is a set bicode over the set S = {0, 1}.

However these set bicodes submit to certain algebraic conditions which is expressed in the form of a definition.

**DEFINITION 3.1.2**: *Let*

$$C = C_1 \cup C_2$$

$$= \left\{\left(x_1^1, x_2^1, \ldots, x_{r_1}^1\right), \left(x_1^1, \ldots, x_{r_2}^1\right), \ldots, \left(x_1^1, x_2^1, \ldots, x_{r_{n_1}}^1\right)\right\} \cup \left\{\left(x_1^2, x_2^2, \ldots, x_{s_1}^2\right), \left(x_1^1, x_2^2, \ldots, x_{s_2}^2\right), \ldots, \left(x_1^2, x_2^2, \ldots, x_{s_{n_2}}^2\right)\right\}$$

*be a set bicode over the set S = {0, 1} when each* $x_i^1, x_j^2 \in$ *{0, 1},* $1 \leq i \leq (r_1, r_2, \ldots, r_{n_1})$, $1 \leq j \leq (s_1, s_2, \ldots, s_{n_2})$. *We demand a set of set bimatrices H =* $H_1 \cup H_2 = \left\{H_1^1 \ldots H_{t_1}^1\right\} \cup \left\{H_1^2 \ldots H_{t_2}^2\right\}$ *such that for every set bicode word* $x = x_i^1 \cup x_j^2 \in C_1 \cup C_2$ *there exists a unique bimatrix* $H_i^1 \cup H_j^2 \in H_1 \cup H_2 = H$ *satisfying the condition*

$$Hx^t = (H_i^1 \cup H_j^2), \quad (x_i^1 \cup x_j^2),$$
$$= H_i^1 \; x_i^1 \cup H_j^2 \; x_j^2 = 0 \cup 0.$$

*We call H the set parity check set bimatrices associated with the set bicode C.*

However it is very important to note that the set bicode does not contain all x such that $Hx^T = (H_i^1 \cup H_j^2)(x_i^1 \cup x_j^2) = (0) \cup (0)$, $1 \leq i \leq n_1$ and $1 \leq j \leq n_2$. The only criteria which is essential is that the set bicodes x which is in C satisfies $Hx^t = (0) \cup (0)$ and nothing more. That is if $x = x_1 \cup x_2$ is such that $Hx^T = H_1x_1^t \cup H_2x_2^t = 0 \cup 0$ it is not necessary $x \in C = C_1 \cup C_2$.

We will illustrate this by a simple example.

***Example 3.1.3:*** Let

$C = C_1 \cup C_2$

= {(1 1 1 1), (0 0 0 0), (1 0 1 0), (1 1 1 1 1 1), (0 0 0 0 0 0), (1 0 0 1 0 0), (0 0 1 0 0 1), (1 1 1 1 1 1), (0 0 0 0 0 0), (0 0 0 0 0 0 0), (1 0 0 0 1 0 0), (0 1 0 0 0 1 0)} ∪ {(1 1 1 1 1), (0 0 0 0 0), (1 1 1 1 1 1), (0 0 0 0 0 0), (0 1 0 0 1 0), (1 0 0 1 0 0), (1 1 1 1), (0 0 0 0), (0 1 0 1)}

$= \{C_1^1, C_2^1, C_3^1\} \cup \{C_1^2, C_2^2, C_3^2\}$

where

$C_1^1$ = {(1 1 1 1), (0 0 0 0), (1 0 1 0)},

$C_2^1$ = {(1 1 1 1 1 1), (0 0 0 0 0 0), (1 0 0 1 0 0), (0 0 1 0 0 1)},

$C_3^1$ = {(1 1 1 1 1 1 1), (0 0 0 0 0 0 0), (1 0 0 0 1 0 0), (0 1 0 0 0 1 0)},

$C_1^2$ = {(1 1 1 1 1), (0 0 0 0 0), (1 0 1 0 0)},

$C_2^2 = \{(0\ 0\ 0\ 0\ 0\ 0), (1\ 1\ 1\ 1\ 1\ 1), (0\ 0\ 1\ 0\ 0\ 1)\}$ and
$C_3^2 = \{(1\ 1\ 1\ 1), (0\ 0\ 0\ 0), (0\ 1\ 0\ 1)\}$

be a set bicode with an associated set bimatrix

$$H = \{H_1 \cup H_2\}$$
$$= \left\{H_1^1, H_2^1, H_3^1\right\} \cup \left\{H_1^2, H_2^2, H_3^2\right\} =$$

$$\left\{\begin{pmatrix} 1 & 0 & 1 & 0 \\ 0 & 1 & 0 & 1 \end{pmatrix}, \begin{pmatrix} 1 & 0 & 0 & 1 & 0 & 0 \\ 0 & 1 & 0 & 0 & 1 & 0 \\ 0 & 0 & 1 & 0 & 0 & 1 \end{pmatrix}, \begin{pmatrix} 1 & 0 & 0 & 0 & 1 & 0 & 0 \\ 0 & 1 & 1 & 0 & 0 & 1 & 0 \\ 0 & 0 & 1 & 1 & 0 & 0 & 1 \end{pmatrix}\right\}$$

$$\cup \left\{\begin{pmatrix} 1 & 0 & 1 & 0 & 0 \\ 0 & 1 & 0 & 1 & 0 \\ 0 & 1 & 0 & 0 & 1 \end{pmatrix}, \begin{pmatrix} 1 & 0 & 0 & 1 & 0 & 0 \\ 0 & 1 & 0 & 0 & 1 & 0 \\ 0 & 0 & 1 & 0 & 0 & 1 \end{pmatrix}, \begin{pmatrix} 1 & 0 & 1 & 0 \\ 0 & 1 & 0 & 1 \end{pmatrix}\right\}$$

which is the set parity check bimatrix of $C = C_1 \cup C_2$.

It is easily verified that all bicode words $x = x_1 \cup x_2$ such that $Hx^T = (0 \cup 0) = H_1x_1^t \cup H_2x_2^t$ is not included in $C = C_1 \cup C_2$.

Further suppose $x = (1\ 0\ 1\ 0) \cup (0\ 0\ 1\ 1\ 1\ 0) = x_1 \cup x_2$ is any received bimessage. We now check the bimessage is the correct one by finding

$$Hx^t = (H_1^1 \cup H_2^2)\ \{(1\ 0\ 1\ 0) \cup (0\ 0\ 1\ 1\ 1\ 0)\}$$
$$= (0\ 0) \cup (1\ 1\ 1)$$
$$\neq (0\ 0) \cup (0\ 0\ 0).$$

So $x \notin C = C_1 \cup C_2$.

In view of this we define the new notion of set bisyndrome of a set bicode word. This is mainly used to detect the errors. However for correction we need to adopt different methods by resending the message or correcting the error using Hamming bidistance. The main use of set parity check bimatrices is mainly for error detection.

**DEFINITION 3.1.3:** *Let C =* $\left(C_1^1, C_2^1, ..., C_{n_1}^1\right) \cup \left(C_1^2, C_2^2, ..., C_{n_2}^2\right)$ = $C_1 \cup C_2$ *be a set bicode over S = {0, 1} where* $C_{j_i}^i$ *are set codes i.e.,* $C_{j_i}^i$ *consists of set of code words of same length but in general* $C_{j_i}^i \neq C_{j_t}^i$*, if* $i \neq t$*;* $1 \leq j_t, j_i \leq n_i$*; i = 1, 2, so* $C_{j_i}^i$ *and* $C_{j_t}^i$ *are code words of different lengths* $1 \leq j_i, j_t \leq n_i$*; i = 1, 2. Suppose H =* $H_1 \cup H_2 = \left\{H_1^1, H_2^1, \ldots, H_{n_1}^1\right\} \cup \left\{H_1^2, H_2^2, \ldots, H_{n_2}^2\right\}$ *be a set parity check bimatrix associated with the set bicode C =* $C_1 \cup C_2$*. Let* $y \in C_1 \cup C_2 = C$ *i.e., y =* $y_{j_1}^1 \cup y_{j_2}^2$*. We define S(y), =* $S(y_{j_1}^1) \cup S(y_{j_2}^2)$*, =* $H_{j_1}^1 \left(y_{j_1}^1\right)^t \cup H_{j_2}^2 \left(y_{j_2}^2\right)^t$ *to be the set bisyndrome of the set bicode C =* $C_1 \cup C_2$*.*

The set bisyndrome helps one to detect the error in the following way.

We find for any y = $y_1 \cup y_2$ a set bicode word, the set bisyndrome S(y) = $S(y_{j_1}^1) \cup S(y_{j_2}^2)$; if S(y) = $(0 \cup 0)$ then we say $y \in C = C_1 \cup C_2$. If $S(y) \neq 0 \cup 0$ and even if one of $S\left(y_{i_t}^t\right) \neq 0$; t = 1 or 2 then also we understand that the error has occurred and it is detected. Thus the error S(y) = $e_1^1 \cup e_2^2$ is called as the set bierror. This technique is useful only in set bierror detection.

We will illustrate this by a simple example.

***Example 3.1.4:*** Let C = $C_1 \cup C_2$ be a set bicode over the set S = {0, 1} where

$C_1^1$ = {(1 1 1 1 1 1), (0 0 0 0 0 0), (1 1 0 0 0 0), (0 1 1 0 0 0), (0 0 0 0 1 1), (0 0 0 1 1 0)},

$C_2^1$ = {(1 1 1 0), (0 1 1 1), (0 0 0 0)},

$C_3^1$ = {(1 1 1 0 1 0 0), (0 1 1 1 0 1 0), (1 1 0 1 0 0 1), (0 0 0 0 0 0 0), (1 1 1 1 1 1 1)},

$C_1^2 =$ {(0 0 0 1 1 1 1), (0 1 1 0 0 1 1), (1 1 1 1 1 1 1), (0 0 0 0 0 0 0), (1 0 1 0 1 0 1), (1 0 0 1 1 0 0)} ,

$C_2^2 =$ {(1 1 0 0), (0 1 0 1), (0 0 0 0)},

$C_3^2 =$ {(0 0 1 0 1 1 1), (0 1 0 1 1 1 0), (1 0 1 1 1 0 0), (0 0 0 0 0 0 0)} and

$C_4^2 =$ {(0 0 0 0 0 0), (0 0 1 1 1 0), (1 0 0 0 1 1), (1 1 1 0 0 0), (0 1 0 1 0 1)}.

Let the set parity check bimatrix

$$H = H_1 \cup H_2$$

$$= \left\{ H_1^1 = \begin{pmatrix} 1 & 1 & 1 & 1 & 0 & 0 \\ 0 & 0 & 0 & 1 & 1 & 1 \\ 1 & 1 & 1 & 1 & 1 & 1 \end{pmatrix}, H_2^1 = \begin{pmatrix} 1 & 0 & 1 & 0 \\ 0 & 1 & 0 & 1 \end{pmatrix}, \right.$$

$$\left. H_3^1 = \begin{pmatrix} 1 & 1 & 1 & 0 & 1 & 0 & 0 \\ 0 & 1 & 1 & 1 & 0 & 1 & 0 \\ 1 & 1 & 0 & 1 & 0 & 0 & 1 \end{pmatrix} \right\} \cup$$

$$\left\{ H_2^1 = \begin{pmatrix} 0 & 0 & 0 & 1 & 1 & 1 & 1 \\ 0 & 1 & 1 & 0 & 0 & 1 & 1 \\ 1 & 0 & 1 & 0 & 1 & 0 & 1 \end{pmatrix}, H_2^2 = \begin{pmatrix} 1 & 1 & 0 & 1 \\ 0 & 1 & 0 & 1 \end{pmatrix}, \right.$$

$$H_3^2 = \begin{pmatrix} 0 & 0 & 1 & 0 & 1 & 1 & 1 \\ 0 & 1 & 0 & 1 & 1 & 1 & 0 \\ 1 & 0 & 1 & 1 & 1 & 0 & 0 \end{pmatrix}$$

and

$$\left. H_4^2 = \begin{pmatrix} 0 & 1 & 1 & 1 & 0 & 0 \\ 1 & 0 & 1 & 0 & 1 & 0 \\ 1 & 1 & 0 & 0 & 0 & 1 \end{pmatrix} \right\}.$$

We observe the following from this example.

(1) Suppose y = (1 1 1 0 0 0) $\cup$ (1 1 1 1 1 0 0) = $y_1 \cup y_2$ is the received set bicode word. To find out whether the received set bicode word is correct we use the method of detecting the set bi error using the set bisyndrome technique.

$$S(y) = S(y_1) \cup S(y_2)$$
$$= H_1^1 y_1^t \cup H_3^2 y_2^t$$

$$= \begin{pmatrix} 1 & 1 & 1 & 1 & 0 & 0 \\ 0 & 0 & 0 & 1 & 1 & 1 \\ 1 & 1 & 1 & 1 & 1 & 1 \end{pmatrix} \begin{bmatrix} 1 \\ 1 \\ 1 \\ 0 \\ 0 \\ 0 \end{bmatrix}$$

$$\cup \begin{pmatrix} 0 & 0 & 1 & 0 & 1 & 1 & 1 \\ 0 & 1 & 0 & 1 & 1 & 1 & 0 \\ 1 & 0 & 1 & 1 & 1 & 0 & 0 \end{pmatrix} \begin{bmatrix} 1 \\ 1 \\ 1 \\ 1 \\ 1 \\ 0 \\ 0 \end{bmatrix}$$

$$= (1\ 0\ 1)^t \cup (0\ 1\ 0)^t$$
$$\neq 0 \cup 0.$$

Hence the received set bicode word has error. Thus we have to adopt to some error correcting techniques to retrieve the correct message.

(2) Further in a bicode C we can send only a bicode word of length $(n_{1,} n_2)$ which is always fixed once a bicode C is choosen. But in case of set bicodes we can send bicode words of different lengths say if the set bicode has biset code words of length $\left(n_1^1, n_2^1, ..., n_{t_1}^1\right) \cup \left(n_1^2, n_2^2, ..., n_{t_2}^2\right)$ then a set bicode word can be of $t_1$, $t_2$ different lengths i.e., of length $\left(n_i^1, n_j^2\right)$; $1 \leq i \leq t_1$, $1 \leq j \leq t_2$. The length of set bicode is called the set bilength denoted generally by $\left(n_1^1, n_2^1, ..., n_{t_1}^1\right) \cup \left(n_1^2, n_2^2, ..., n_{t_2}^2\right)$; $t_i > t_2$ or $t_1 < t_2$ or $t_1 = t_2$. Also some $n_j^1$ can be equal to $n_i^2$.

Thus the use of these set bicodes is handy when one wants to send set bicode words of different lengths varying from time to time.

Thus x = (0 0 0 0 1 1) ∪ (0 0 0 1 1 1 1) is a set bicode word of length (6, 7). x = (1 1 1 0) ∪ (1 1 0 0) is a set bicode word of length (4, 4); x = (1 1 1 0) ∪ (0 1 0 1 0 1) ∈ $C = C_1 \cup C_2$ is set bicode word of length (4, 6) and so on. The set bilength of the given set bicode is (6, 4, 7) ∪ (7, 4, 7, 6).

The only means to correct the errors in these set bicodes is use the Hamming bidistance once we know error has been detected. Suppose $x = x_1 \cup x_2$ and $y = y_1 \cup y_2$ the Hamming bidistance $d_b$ (x, y) = $(d(x_1, y_1), d(x_2, y_2))$ where $d(x_i, y_i)$ is the usual Hamming distance between $x_i$ and $y_i$; $1 \leq i \leq 2$. If $y = y_1 \cup y_2$ is the received set bicode word which has some error then we find d(x, y), for all $x = x_1 \cup x_2$ where $x_1 \in C_i^1$ and $x_2 \in C_j^2$. The set bicode word x such that d(x, y), is minimum is taken as the approximately correct set bicode word.

We shall illustrate this by a simple example.

***Example 3.1.5:*** Let $C = \left\{C_1^1, C_2^1, C_3^1\right\} \cup \left\{C_1^2, C_2^2, C_3^2, C_4^2\right\} = C_1 \cup C_2$ be a set bicode over the set S = {0, 1}. $C_1^1$ = {(1 1 1 1 1), (1 1 0 0 0), (0 0 0 0 0), (0 0 0 1 1)} be a set code with the set parity check matrix.

$$H_1^1 = \begin{pmatrix} 1 & 1 & 0 & 0 & 0 \\ 0 & 0 & 0 & 1 & 1 \\ 1 & 0 & 1 & 0 & 0 \end{pmatrix}.$$

$C_2^1$ = {(1 1 1 1 0 1), (0 0 0 0 0 0), (1 0 1 0 0 0), (0 1 0 0 0 1)} be a set code with the set parity check matrix,

$$H_2^1 = \begin{pmatrix} 1 & 0 & 1 & 0 & 1 & 0 \\ 0 & 1 & 0 & 0 & 0 & 1 \\ 0 & 1 & 0 & 1 & 0 & 1 \end{pmatrix}.$$

$C_3^1$ = {(0 0 0 0 0 0), (1 1 1 0 0 0), (1 0 1 1 0 1), (0 0 1 1 1 0), (1 1 0 1 1 0)} be a set code with set parity check matrix;

$$H_3^1 = \begin{pmatrix} 0 & 1 & 1 & 1 & 0 & 0 \\ 1 & 0 & 1 & 0 & 1 & 0 \\ 1 & 1 & 0 & 0 & 0 & 1 \end{pmatrix}.$$

$C_1^2$ = {(0 0 0 0 0 0), (1 0 0 0 1 1), (0 1 1 0 1 1), (1 1 1 0 0 0)} be a set code with set parity check matrix,

$$H_1^2 = \begin{pmatrix} 0 & 1 & 1 & 1 & 0 & 0 \\ 1 & 0 & 1 & 0 & 1 & 0 \\ 1 & 1 & 0 & 0 & 0 & 1 \end{pmatrix}.$$

$C_2^2$ = {(0 0 0 0 0 0 0), (0 1 1 0 0 1 1), (0 0 0 1 1 1 1), (1 0 1 0 1 0 1), (0 1 0 1 0 1 0)} be a set bicode with set parity check matrix

$$H_2^2 = \begin{pmatrix} 0 & 0 & 0 & 1 & 1 & 1 & 1 \\ 0 & 1 & 1 & 0 & 0 & 1 & 1 \\ 1 & 0 & 1 & 0 & 1 & 0 & 1 \end{pmatrix}.$$

$C_3^2$ = {(0 0 0 0), (1 0 1 1), (1 1 1 0)} be a set bicode with set parity check matrix,

$$H_3^2 = \begin{pmatrix} 1 & 0 & 1 & 0 \\ 1 & 1 & 0 & 1 \end{pmatrix}.$$

$C_4^2$ = {(1 1 1 1 1 1), (1 0 1 1 0 1), (0 0 0 0 0 0), (0 1 0 0 1 0), (1 1 0 1 1 0)} is a set bicode with set parity check matrix,

$$H_4^2 = \begin{pmatrix} 1 & 0 & 0 & 1 & 0 & 0 \\ 0 & 1 & 0 & 0 & 1 & 0 \\ 0 & 0 & 1 & 0 & 0 & 1 \end{pmatrix}.$$

Thus C = $C_1 \cup C_2$ = $\{C_1^1, C_2^1, C_3^1\} \cup \{C_1^2, C_2^2, C_3^2, C_4^2\}$ is a set bicode over S = {0, 1} with the set parity check bimatrix H = $H_1 \cup H_2$ = $\{H_1^1, H_2^1, H_3^1\} \cup \{H_1^2, H_2^2, H_3^2, H_4^2\}$. Suppose y = (1 1 0 1 0 1) ∪ (1 0 1 0 1 1 1), is the received set bicode word, we first check whether y ∈ $C_1 \cup C_2$. We find the set bisyndrome

$$S_b(y) = S(y_1) \cup S(y_2)$$

$$= H_2^1 y_1^t \cup H_2^2 y_2^t$$

$$= \begin{pmatrix} 1 & 0 & 1 & 0 & 1 & 0 \\ 0 & 1 & 0 & 0 & 0 & 1 \\ 0 & 1 & 0 & 1 & 0 & 1 \end{pmatrix} \begin{bmatrix} 1 \\ 1 \\ 0 \\ 1 \\ 0 \\ 1 \end{bmatrix} \cup$$

$$\begin{pmatrix} 0 & 0 & 0 & 1 & 1 & 1 & 1 \\ 0 & 1 & 1 & 0 & 0 & 1 & 1 \\ 1 & 0 & 1 & 0 & 1 & 0 & 1 \end{pmatrix} \begin{bmatrix} 1 \\ 0 \\ 1 \\ 0 \\ 1 \\ 1 \\ 1 \end{bmatrix}$$

$$= \begin{bmatrix} 1 \\ 0 \\ 1 \end{bmatrix} \cup \begin{bmatrix} 1 \\ 1 \\ 0 \end{bmatrix}$$

$$\neq \quad (0) \cup (0).$$

Thus the received set bicode word does not belong to $C_1 \cup C_2$. Now we find the approximately a close set bicode word to y using the set Hamming bidistance. $d(x, y) = d(x_1, y_1) \cup d(x_2, y_2)$ where $x = x_1 \cup x_2 = (1\ 1\ 1\ 1\ 0\ 1) \cup (0\ 1\ 1\ 0\ 0\ 1\ 1)$. $d(x, y) = (1, 3)$ we choose $x_1 = (1\ 1\ 1\ 1\ 0\ 1)$ to be the approximately close to the code word (1 1 0 1 0 1). Now we find d (x, y) where

$$x = x_1 \cup (0\ 0\ 0\ 1\ 1\ 1\ 1)$$

with

$$x_1 = (1\ 1\ 1\ 1\ 0\ 1)$$
$$d(x, y) = (1, 3)$$

$d(x, y)$ where $x = x_1 \cup (1\ 0\ 1\ 0\ 1\ 0\ 1) = (1, 1)$.

Thus we take (1 0 1 0 1 0 1) to be the received probable set code word in the place of the received set code word (1 0 1 0 1 1 1). Thus the nearest set bicode word is $y_1 = (1\ 1\ 1\ 1\ 0\ 1) \cup (1\ 0\ 1\ 0\ 1\ 1\ 1)$.

Now we proceed on to define the notion of special types of set bicodes.

**DEFINITION 3.1.4**: *Let*

$$C = C_1 \cup C_2$$

$$= \left\{\left(y_1^1 \ldots y_{r_{t_1}}^1\right), \left(x_1^1 \ldots x_{r_{t_1}}^1\right) \Big/ y_i^1 = 0, 1 \leq i \leq r_{t_1}, x_i^1 = 1, 1 \leq i \leq r_{t_1} ;\ t_1 = \{1, 2, \ldots, n_1\}\right\}$$

$$\cup \left\{\left(y_1^2 \ldots y_{s_{t_2}}^2\right), \left(x_1^2 \ldots x_{s_{t_2}}^2\right) \Big/ y_i^2 = 0, 1 \leq i \leq r_{t_2}, x_i^2 = 1, 1 \leq i \leq r_{t_2}\right.$$

*and $t_2 = 1, 2, \ldots, n_2\}$ be a set bicode where either each of the tuples are zeros or ones. The set parity check bimatrix $H = H_1 \cup H_2$* $= \left\{H_1^1, H_2^1, \ldots, H_{n_1}^1\right\} \cup \left\{H_1^2, H_2^2, \ldots, H_{n_2}^2\right\}$ *where* $H_i^1$ *is a* $(r_i - 1) \times r_i$ *matrix of the form*

$$\begin{pmatrix} 1 & 1 & 0 & 0 & \cdots & 0 \\ 1 & 0 & 1 & \cdots & \cdots & 0 \\ \vdots & \vdots & \vdots & & & \vdots \\ 1 & 0 & 0 & \cdots & \cdots & 1 \end{pmatrix}$$

*where first column has only ones and the rest is a* $(r_i - 1) \times (r_i - 1)$ *identity matrix,* $i = 1, 2, \ldots, n_1$ *and* $H_j^2$ *is a* $(s_j - 1) \times (s_j - 1)$ *matrix of the form*

$$\begin{pmatrix} 1 & 1 & 0 & \cdots & 0 \\ 1 & 0 & 1 & \cdots & 0 \\ \vdots & \vdots & \vdots & & \vdots \\ 1 & 0 & 0 & \cdots & 1 \end{pmatrix}$$

*where first column has only ones and the rest is a* $(s_j - 1) \times (s_j - 1)$ *identity matrix;* $j = 1, 2, \ldots, n_2$. *We call this the set bicode* $C = C_1 \cup C_2$ *as the repetition set bicode.*

We shall illustrate this by an example.

***Example 3.1.6:*** Let $C = C_1 \cup C_2$ = {(0 0 0 0), (1 1 1 1), (0 0 0 0 0), (1 1 1 1 1), (0 0 0), (1 1 1), (1 1 1 1 1 1 1), (0 0 0 0 0 0 0)} ∪ {(0 0 0 0), (1 1 1 1), (0 0 0 0 0 0), (1 1 1 1 1 1), (0 0 0 0 0 0 0 0), (1 1 1 1 1 1 1 1), (0 0 0 0 0), (1 1 1 1 1)} be the repetition set bicode

$C = C_1 \cup C_2 = C_1^1$ = {(0 0 0 0), (1 1 1 1)},

$C_2^1$ = {(0 0 0 0 0), (1 1 1 1 1)},

{(1 1 1 1 1 1 1), (0 0 0 0 0 0 0)} = $C_3^1$,

$C_4^1$ = {(0 0 0), (1 1 1)} and $C_1^2$ = {(0 0 0 0), (1 1 1 1)},

$C_2^2$ = {(0 0 0 0 0 0), (1 1 1 1 1 1)},

$$C_3^2 = \{(0\ 0\ 0\ 0\ 0\ 0\ 0\ 0), (1\ 1\ 1\ 1\ 1\ 1\ 1\ 1)\}$$

and

$$C_4^2 = \{(0\ 0\ 0\ 0\ 0), (1\ 1\ 1\ 1\ 1)\}$$

with the associated set parity check bimatrix

$H = H_1 \cup H_2 =$

$$\left\{ \left\{ \begin{pmatrix} 1 & 1 & 0 & 0 \\ 1 & 0 & 1 & 0 \\ 1 & 0 & 0 & 1 \end{pmatrix} \right\} = H_1^1 \right.$$

$$H_2^1 = \left\{ \begin{pmatrix} 1 & 1 & 0 & 0 & 0 \\ 1 & 0 & 1 & 0 & 0 \\ 1 & 0 & 0 & 1 & 0 \\ 1 & 0 & 0 & 0 & 1 \end{pmatrix} \right\}$$

$$H_3^1 = \left\{ \begin{pmatrix} 1 & 1 & 0 & 0 & 0 & 0 & 0 \\ 1 & 0 & 1 & 0 & 0 & 0 & 0 \\ 1 & 0 & 0 & 1 & 0 & 0 & 0 \\ 1 & 0 & 0 & 0 & 1 & 0 & 0 \\ 1 & 0 & 0 & 0 & 0 & 1 & 0 \\ 1 & 0 & 0 & 0 & 0 & 0 & 1 \end{pmatrix} \right\}$$

$$\left. H_4^1 = \left\{ \begin{pmatrix} 1 & 1 & 0 \\ 0 & 0 & 1 \end{pmatrix} \right\} \right\} \cup$$

$$\left\{ H_1^2 = \left\{ \begin{pmatrix} 1 & 1 & 0 & 0 \\ 1 & 0 & 1 & 0 \\ 1 & 0 & 0 & 1 \end{pmatrix} \right\} \right.$$

$$H_2^2 = \left\{ \begin{pmatrix} 1 & 1 & 0 & 0 & 0 & 0 \\ 1 & 0 & 1 & 0 & 0 & 0 \\ 1 & 0 & 0 & 1 & 0 & 0 \\ 1 & 0 & 0 & 0 & 1 & 0 \\ 1 & 0 & 0 & 0 & 0 & 1 \end{pmatrix} \right\},$$

$$H_3^2 = \left\{ \begin{pmatrix} 1 & 1 & 0 & 0 & 0 & 0 & 0 & 0 \\ 1 & 0 & 1 & 0 & 0 & 0 & 0 & 0 \\ 1 & 0 & 0 & 1 & 0 & 0 & 0 & 0 \\ 1 & 0 & 0 & 0 & 1 & 0 & 0 & 0 \\ 1 & 0 & 0 & 0 & 0 & 1 & 0 & 0 \\ 1 & 0 & 0 & 0 & 0 & 0 & 1 & 0 \\ 1 & 0 & 0 & 0 & 0 & 0 & 0 & 1 \end{pmatrix} \right\}$$

$$H_4^2 = \left\{ \begin{pmatrix} 1 & 1 & 0 & 0 & 0 \\ 1 & 0 & 1 & 0 & 0 \\ 1 & 0 & 0 & 1 & 0 \\ 1 & 0 & 0 & 0 & 1 \end{pmatrix} \right\}.$$

Thus we can send any repetition set bicode of varying or desired length. This is the advantage of the repetition bicodes.

Now we proceed on to define the notion of set parity check bicode.

**DEFINITION 3.1.5**: *Let*

$$C = C_1 \cup C_2 = \left\{C_1^1, C_2^1, \ldots, C_{n_1}^1\right\} \cup \left\{C_1^2, C_2^2, \ldots, C_{n_2}^2\right\}$$

*be a set bicode over S = {0, 1} with set parity check bimatrices* $H_1 \cup H_2$ *= {* $H_1^1$ *= (1 1 …1),* $H_2^1$ *= (1 1 … 1), …,* $H_{n_1}^1$ *= (1 1 1 … 1), vectors of lengths say* $(r_1, r_2, \ldots, r_{n_1})$ *respectively}* $\cup$ *{* $H_1^2$ *= (1 1 …1),* $H_2^2$ *= (1 1 …1), …,* $H_{n_2}^1$ *= (1 1 …1)} be vectors of* *length say {$s_1, s_2, \ldots, s_{n_2}$} respectively. We call $C = C_1 \cup C_2$ to be a set parity check bicode.*

We illustrate this by an example.

***Example 3.1.7:*** Let C = $C_1 \cup C_2$ = {(1 1 0 0 0), (0 1 1 0 0), (0 0 0 0 0), (1 1 0 1 1), (1 0 0 0 1), (1 1 1 1), (1 1 0 0), (0 0 0 0), (1 1 1 1 1 0), (1 0 1 0 0 0 0), (0 1 0 0 0 1 0), (1 1 0 0 0 1 1), (0 0 1 1 1 1 0)} ∪ {(1 1 1 1 1 1), (1 1 0 0 0 0), (0 0 0 0 1 1), (0 0 1 1 1 0 0), (0 0 0 0 0 0), (1 1 1 1 1 1 1 1), (1 1 1 1 0 0 0 0), (1 1 0 0 1 1 0 0), (1 1 1 0 0 1 1 1), (0 0 0 0 0 0 0 0)} be a set bicode with the parity check set bimatrix H = $H_1 \cup H_2$ ={(1 1 1 1 1) = $H_1^1$, (1 1 1 1) = $H_2^1$, (1 1 1 1 1 1 1) = $H_3^1$} ∪ {(1 1 1 1 1 1) = $H_1^2$, $H_2^2$ = (1 1 1 1 1 1 1 1)}. Clearly C is a parity check set bicode.

Next we proceed on to define the notion of binary Hamming set bicode.

**DEFINITION 3.1.6:** *Let*

$$C = C_1 \cup C_2 = \left\{C_1^1, C_2^1, \ldots, C_{n_1}^1\right\} \cup \left\{C_1^2, C_2^2, \ldots, C_{n_2}^2\right\}$$

*be a set bicode with a set parity check bimatrix*

$$H = H_1 \cup H_2 = \left\{H_1^1, H_2^1, \ldots, H_{n_1}^1\right\} \cup \left\{H_1^2, H_2^2, \ldots, H_{n_2}^2\right\}$$

*where each $H_{t_i}^i$ is a $m_{t_i}^i \times (2^{m_{t_i}^i} - 1)$ parity check matrix i = 1, 2; $1 \leq t_i \leq n_i$; whose columns consists of all non zero binary vectors of length $m_{t_i}^i$; $1 \leq t_i \leq n_i$, i = 1, 2. We first state each $C_{t_i}^i$ is a set of code words associated with $H_{t_i}^i$ but we do not demand $C_{t_i}^i$ to contain all x such that $H_{t_i}^i x^t = (0)$, only $C_{t_i}^i$ is a subset of all the code words associated with $H_{t_i}^i$; $1 \leq t_i \leq n_i$; $1 \leq i \leq 2$. We call this $C = C_1 \cup C_2$ to be a binary Hamming set bicode.*

We illustrate this by a simple example.

***Example 3.1.8:*** Let C = $C_1 \cup C_2$ = {(0 0 0 0 0 0 0), (1 0 0 1 1 0 0), (1 1 1 1 1 1 1), (0 0 0 0 0 0 0 0), (1 1 1 1 1 1 1 1), (1 0 1 0 1 0 1 0)} $\cup$ {(0 0 0 0 0 0 0), (1 0 0 1 1 0 1), (0 1 0 1 0 1 1), (0 0 0 0 0 0 0 0 0 0 0 0 0 0 0), (1 0 0 0 1 0 0 1 1 0 1 0 1 1 1), (0 0 0 1 0 0 1 1 0 1 0 1 1 1 1)} be a set bicode associated with a set parity check bimatrix

$$H = H_1 \cup H_2 = \left\{H_1^1, H_2^1\right\} \cup \left\{H_1^2,\ H_2^2\right\}$$

where

$$H_1^1 = \begin{pmatrix} 0 & 0 & 0 & 1 & 1 & 1 & 1 \\ 0 & 1 & 1 & 0 & 0 & 1 & 1 \\ 1 & 0 & 1 & 0 & 1 & 0 & 1 \end{pmatrix},$$

$$H_2^1 = \begin{pmatrix} 1 & 1 & 1 & 1 & 1 & 1 & 1 & 1 \\ 0 & 0 & 0 & 1 & 1 & 1 & 1 & 0 \\ 0 & 1 & 1 & 0 & 0 & 1 & 1 & 0 \\ 1 & 0 & 1 & 0 & 1 & 0 & 1 & 0 \end{pmatrix}$$

$$H_1^2 = \begin{pmatrix} 1 & 0 & 0 & 1 & 1 & 0 & 1 \\ 0 & 1 & 0 & 1 & 0 & 1 & 1 \\ 0 & 0 & 1 & 0 & 1 & 1 & 1 \end{pmatrix}$$

and

$$H_2^2 = \begin{pmatrix} 1 & 0 & 0 & 0 & 1 & 0 & 0 & 1 & 1 & 0 & 1 & 0 & 1 & 1 & 1 \\ 0 & 1 & 0 & 0 & 1 & 1 & 0 & 1 & 0 & 1 & 1 & 1 & 1 & 0 & 0 \\ 0 & 0 & 1 & 0 & 0 & 1 & 1 & 0 & 1 & 0 & 1 & 1 & 1 & 1 & 0 \\ 0 & 0 & 0 & 1 & 0 & 0 & 1 & 1 & 0 & 1 & 0 & 1 & 1 & 1 & 1 \end{pmatrix}.$$

Clearly C = $C_1 \cup C_2$ is a Hamming set bicode.

Now we proceed on to define two types of weight m set bicodes.

**DEFINITION 3.1.7:** *Let*

$$V = \left\{V_1^1, \ldots, V_{n_1}^1\right\} \cup \left\{V_1^2, V_1^2, \ldots, V_{n_2}^2\right\}$$

*be a set bicode with the set parity check bimatrix.*

$$H = H_1 \cup H_2 = \left\{H_1^1, H_2^1, \ldots, H_{n_1}^1\right\} \cup \left\{H_1^2, H_2^2, \ldots, H_{n_2}^2\right\}$$

*where $V_j^i$ are set codes with associated set parity check matrix $H_j^i$, $1 \le j \le n_i$, and i = 1, 2. If the set Hamming biweight of each and every bicode word $C_i^1 \cup C_j^2$ in $V = V_1 \cup V_2$ is (m, m), m < length of code word $C_i^1$ and $C_j^2$ in $V_i^1$ and $V_j^2$ respectively; $1 \le j \le n_2$ and $1 \le i \le n_1$. Then we call $V = V_1 \cup V_2$ to be a (m, m), biweight set bicode.*

These types of codes will be more useful in cryptography and computers in which minimum number of bits is fixed.

We first illustrate this by a simple example.

***Example 3.1.9:*** Let V = $V_1 \cup V_2$ = { $V_1^1$ = {(0 0 0 0 0 0), (1 1 1 1 0 0), (0 0 1 1 1 1), (1 1 1 0 1 1), (0 1 1 1 1 0), (0 1 0 1 1 1)}, $V_2^1$ = {(1 1 1 1 0 0 0), (0 0 0 0 0 0 0), (1 1 1 0 1 0 0), (1 0 1 0 1 0 1), (1 1 0 1 1 0 0), (0 0 1 1 1 1 0), (1 1 0 0 0 1 1)}; $V_3^1$ = {(1 1 1 1), (0 0 0 0)}} ∪ {{ $V_1^2$ = (0 0 0 0 0 0 0 0), (1 1 1 1 0 0 0 0), (1 1 0 0 1 1 0 0), (0 0 0 1 1 1 1 0), (1 1 0 0 0 0 1 1), (1 0 1 0 1 0 1 0), (1 1 0 1 1 0 0 0), (0 1 1 0 1 0 1 0)}, $V_2^2$ = {(1 1 1 1 0 0 0), (1 1 0 0 1 1 0), (0 0 0 0 0 0 0), (1 0 1 1 0 1 0), (0 1 1 0 1 0 1)}, $V_3^2$ = {(1 1 1 1 0), (0 0 0 0 0), (0 1 1 1 1), (1 1 0 1 1)}} be a set bicode with 4 biweight. Every non zero bicode is of 4 weight viz. (4, 4) or (0, 4) or (4, 0). One of the main advantages of the set bicode of set biweight (4, 4); the error detection is easy.

We show yet another example of a set bicode of same biweight (m, m).

***Example 3.1.10:*** Let V = $V_1 \cup V_2$ = {{ $V_1^1$ = {(1 1 1 1 1 1 0 0 0), (1 1 1 1 1 0 1 0 0), (0 0 0 1 1 1 1 1 1), (0 0 0 0 0 0 0 0 0), (0 1 0 1 0 1 1 1 1), (1 0 1 1 1 1 0 1 0), (1 1 0 1 1 0 1 1 0)}, $V_2^1$ = {(1 1 1 1 1 1), (0 0 0 0 0 0)}, $V_3^1$ = {(0 1 1 1 1 1 1), (1 1 1 0 1 1

1), (0 0 0 0 0 0 0), (1 1 1 1 0 1 1), (1 1 0 1 1 1 1), (1 1 1 1 1 0 1)}, $V_4^1$ = {(1 1 1 1 0 0 1 1), (0 0 1 1 1 1 1 1), (0 1 1 1 1 1 1 0), (1 1 1 0 0 1 1 1), (0 0 0 0 0 0 0 0), (1 1 0 1 1 1 1 0), (0 1 1 1 1 0 1 1)}} $\cup$ {{ $V_1^2$ = (1 1 1 1 1 1 0 0 0 0), (0 0 0 0 0 0 0 0 0 0), (1 1 1 0 0 0 0 1 1 1), (1 1 0 0 1 1 1 1 0 0), (0 0 0 0 1 1 1 1 1 1), (1 1 1 0 0 1 1 1 1 0 0), (1 1 0 0 0 0 1 1 1 1)}, $V_2^2$ = {(1 1 1 1 1 1 0), (0 1 1 1 1 1 1), (0 0 0 0 0 0 0), (1 1 1 0 1 1 1), (1 1 1 1 0 1 1)}, $V_3^2$ = {(1 1 1 1 1 1), (0 0 0 0 0 0)}} be a set bicode of (6, 6) biweight over the set S = {0, 1}.

Now we proceed on to define a (m, n) biweight set bicode or a set bicode of (m, n) biweight m $\neq$ n.

**DEFINITION 3.1.8:** *Let*

$$V = V_1 \cup V_2 = \left\{V_1^1, \ldots, V_{n_1}^1\right\} \cup \left\{V_1^2, V_1^2, \ldots, V_{n_2}^2\right\}$$

*be a set bicode over the set S = {0, 1}. Suppose $V_1$ is a set code of m-weight or m-weighted set code and $V_2$ is a n-weighted set code (or a set code of n weight) m $\neq$ n. Then we call V = $V_1 \cup V_2$ to be a (m, n), biweight set bicode.*

We illustrate this by an example.

***Example 3.1.11:*** Let V = $V_1 \cup V_2$ = { $V_1^1$ = {(0 0 1 1 1), (0 0 0 0 0), (1 1 0 1 0), (1 1 1 0 0), (0 1 1 1 0)}, $V_2^1$ = {(1 1 1 0 0 0), (0 0 0 0 0 0), (0 0 0 1 1 1), (1 1 0 1 0 0), (1 0 0 1 1 0), (0 0 1 1 1 0), (0 1 0 1 0 1)}, $V_3^1$ = {(1 1 1 0 0 0 0), (0 0 0 0 0 0 0), (0 1 1 1 0 0 0), (0 0 1 1 1 0 0), (0 0 0 1 1 1 0), (1 0 1 0 1 0 0), (0 1 0 1 0 1 0), (1 1 0 0 0 0 1)}} $\cup$ { $V_1^2$ = {(0 0 0 0), (1 1 1 0), (0 1 1 1), (1 0 1 1), (1 1 0 1)}, $V_2^2$ = {(1 1 1 0 0), (0 0 0 0 0), (1 1 0 1 0), (0 0 1 1 1), (1 0 1 1 0), (0 1 1 1 0)}, $V_3^2$ = {(0 0 0 0 0 1 1 1), (0 0 0 0 0 0 0 0), (1 1 1 0 0 0 0 0), (1 0 0 0 0 1 1 0)}, $V_4^2$ = {(1 1 1 0 0 0), (0 0 0 0 0 0), (1 1 0 1 0 0), (0 0 1 1 1 0), (1 0 1 0 1 0), (0 1 0 1 0 1)}} is a (3, 3) biweighted set bicode over the set (0, 1) = S.

These bicodes which are biweighted is very useful for by the sheer observation one can detect the errors.

Next we give an example of a (m, n) biweighted set bicode ($m \neq n$).

***Example 3.1.12:*** Let $V = V_1 \cup V_2$ = { $V_1^1$ = {(1 1 1 1 0 0), (0 0 0 0 0 0), (1 1 0 0 1 1), (0 1 1 1 1 0)}, $V_2^1$ = {(1 1 1 1), (0 0 0 0)} $V_3^1$ = {(0 0 0 1 1 1 1), (0 0 0 0 0 0 0), (1 1 1 0 1 0 0), (0 1 1 1 1 0 0)}} $\cup$ { $V_1^2$ = {(1 1 1 0 0 0), (0 0 0 1 1 1), (0 0 0 0 0 0), (1 1 0 0 0 1), (0 1 1 0 1 0), (0 1 1 1 0 0), (1 0 1 0 1 0)}, $V_2^2$ = {(0 0 0 0), (1 1 1 0), (0 1 1 1), (1 1 0 1)}, $V_3^2$ = {(1 1 0 0 1), (0 0 0 0 0), (1 1 1 0 0), (0 1 1 1 0), (0 0 1 1 1)}, $V_4^2$ = {(0 0 0 0 0 1 1 1 1), (1 1 0 0 0 0 0 1), (0 0 0 0 0 0 0 0), (0 1 1 1 0 0 0 0), (1 0 1 0 0 0 1 0), (1 0 0 0 1 0 0 1)}} be a (4, 3) biweighted set bicode over {0, 1}.

Now we proceed on to define set cyclic bicode.

**DEFINITION 3.1.9:** *Let $V = V_1 \cup V_2$, if each of the $V_i$'s is a set code such that every code word is cyclic of varying lengths; $1 \leq i \leq 2$; then we call $V = V_1 \cup V_2$ to be a set cyclic bicode.*

We denote this by a simple example.

***Example 3.1.13:*** Let $V = V_1 \cup V_2$ where $V_1$ = { $V_1^1$ = {(0 0 0 0 0), (1 1 1 1 1)}, $V_2^1$ = {(1 1 0 1 1 0), (0 1 1 0 1 1), (1 0 1 1 0 1), (0 0 0 0 0 0), (1 1 1 1 1 1)}, $V_3^1$ = {(1 1 1 1 1 0 0 0), (0 1 1 1 1 1 0 0), (0 0 1 1 1 1 1 0), (0 0 0 1 1 1 1 1), (1 0 0 0 1 1 1 1), (1 1 0 0 0 1 1 1), (1 1 1 0 0 0 1 1), (1 1 1 1 0 0 0 1), (0 0 0 0 0 0 0 0)}} is a cyclic set code and $V_2$ = { $V_1^2$ = {(1 1 1 1 1 1), (0 0 0 0 0 0)}, $V_2^2$ = {(0 0 0 0 0), (1 1 1 1 0), (0 1 1 1 1), (1 0 1 1 1), (1 1 0 1 1), (1 1 1 0 1)}, $V_3^2$ = {(1 1 1 1 0 0 0), (0 0 0 0 0 0 0), (0 1 1 1 1 0 0), (0 0 1 1 1 1 0), (0 0 0 1 1 1 1), (1 1 0 0 0 1 1), (1 1 1 0 0 0

1)}} is again a cyclic set code. Thus V = $V_1 \cup V_2$ is a cyclic set bicode.

Now we define the notion of cyclic set bicode of (m, m), biweight.

**DEFINITION 3.1.10:** *Let V = $V_1 \cup V_2$ be a cyclic set bicode over the set S. If V is also a (m, m) biweighted set bicode then we call V to be a cyclic set biweighted bicode over the set S.*

We illustrate this by some examples.

***Example 3.1.14:*** Let V = $V_1 \cup V_2$ = { $V_1^1$ = {(0 0 0 0 0 0), (1 1 1 1 0 0), (0 1 1 1 1 0), (0 0 1 1 1 1), (1 0 0 1 1 1), (1 1 0 0 1 1), (1 1 1 0 0 1)}, $V_2^1$ = {(0 1 0 1 0 1 1), (0 0 0 0 0 0 0), (1 0 1 0 1 0 1), (1 1 0 1 0 1 0), (0 1 1 0 1 0 1), (1 0 1 1 0 1 0), (0 1 0 1 1 0 1), (1 0 1 0 1 1 0)}, $V_3^1$ = {(1 1 0 0 0 0 1 1), (0 0 0 0 0 0 0 0), (1 1 1 0 0 0 0 1), (1 1 1 1 0 0 0 0), (0 1 1 1 1 0 0 0), (0 0 1 1 1 1 0 0), (0 0 0 1 1 1 1 0), (0 0 0 0 1 1 1 1), (1 0 0 0 0 1 1 1)}} ∪ { $V_1^2$ = {(1 1 1 1 0 0 0), (0 0 0 0 0 0 0), (0 1 1 1 1 0 0), (0 0 1 1 1 1 0), (0 0 0 1 1 1 1), (1 0 0 0 1 1 1), (1 1 0 0 0 1 1), (1 1 1 0 0 0 1)}, $V_2^2$ = {(1 0 1 0 1 1), (1 1 0 1 0 1), (1 1 1 0 1 0), (0 1 1 1 0 1), (1 0 1 1 1 0), (0 1 0 1 1 1), (0 0 0 0 0 0)}, $V_3^2$ = {(0 0 0 0), (1 1 1 1)}, $V_4^2$ = {(1 0 1 0 1 0 1 0), (0 1 0 1 0 1 0 1), (0 0 0 0 0 0 0 0)}} is a (4, 4) beweighted set cyclic bicode.

*Note:* If V = $V_1 \cup V_2$ is a (m, n) biweighted set bicode which is also cyclic then we call V to be a cyclic biweighted set bicode.

We illustrate this also by an example.

***Example 3.1.15:*** Let V = $V_1 \cup V_2$ = { $V_1^1$ = {(0 0 0 0 0), (1 1 1 0 0), (0 1 1 1 0), (0 0 1 1 1), (1 0 0 1 1), (1 1 0 0 1)}, $V_2^1$ = {(0 1 1 0 0 0 1), (1 0 1 1 0 0 0), (0 1 0 1 1 0 0), (0 0 1 0 1 1 0), (0 0 0 1 0 1 1), (1 0 0 0 1 0 1), (1 1 0 0 0 1 0), (0 0 0 0 0 0 0)}, $V_3^1$ = {(1

0 1 0 1 0), (0 1 0 1 0 1), (0 0 0 0 0 0)}} ∪ { $V_1^2$ = {(1 1 0 1), (1 1 1 0), (0 1 1 1), (1 0 1 1), (0 0 0 0)}, $V_2^2$ = {(1 0 1 0 1), (1 1 0 1 0), (0 1 1 0 1), (1 0 1 1 0), (0 1 0 1 1), (0 0 0 0 0)}, $V_3^2$ = {(1 0 1 0 1 0 0), (0 1 0 1 0 1 0), (0 0 1 0 1 0 1), (1 0 0 1 0 1 0), (0 1 0 0 1 0 1), (1 0 1 0 0 1 0), (0 1 0 1 0 0 1), (0 0 0 0 0 0 0)}, $V_4^2$ = {(1 1 0 0 0 1), (1 1 1 0 0 0), (0 1 1 1 0 0), (0 0 1 1 1 0), (0 0 0 1 1 1), (1 0 0 0 1 1), (0 0 0 0 0 0)}} is again a (3, 3) beweighted cyclic set bicode.

***Example 3.1.16:*** Let V = $V_1 \cup V_2$ = {{(1 1 1 0 0 0 1), (1 1 1 1 0 0 0), (0 1 1 1 1 0 0), (0 0 1 1 1 1 0), (0 0 0 1 1 1 1), (1 0 0 0 1 1 1), (1 1 0 0 0 1 1), (0 0 0 0 0 0 0)} = $V_1^1$, $V_2^1$ = {(0 1 1 1 0 1), (1 0 1 1 1 0), (0 1 0 1 1 1), (1 0 1 0 1 1), (1 1 0 1 0 1), (1 1 1 0 1 0), (0 0 0 0 0 0)}, $V_3^1$ = {(1 0 1 0 1 0 1), (1 1 0 1 0 1 0), (0 1 1 0 1 0 1), (1 0 1 1 0 1 0), (0 1 0 1 1 0 1), (1 0 1 0 1 1 0), (0 1 0 1 0 1 1), (0 0 0 0 0 0 0)}} ∪ { $V_1^2$ = {(1 0 1 0 1 0 0), (0 1 0 1 0 1 0), (0 0 1 0 1 0 1), (1 0 0 1 0 1 0), (0 1 0 0 1 0 1), (1 0 1 0 0 1 0), (0 1 0 1 0 0 1), (0 0 0 0 0 0 0)}, $V_2^2$ = {(1 1 1 0 0 0 0), (0 1 1 1 0 0 0), (0 0 1 1 1 0 0), (0 0 0 1 1 1 0), (0 0 0 0 1 1 1), (1 0 0 0 0 1 1), (1 1 0 0 0 0 1), (0 0 0 0 0 0 0)}, $V_3^2$ = {(1 0 1 0 1 0), (0 1 0 1 0 1), (0 0 0 0 0 0)}, $V_4^2$ = {(1 1 1 0 0), (0 1 1 1 0), (0 0 1 1 1), (1 0 0 1 1), (1 1 0 0 1), (0 0 0 0 0)}}. V = $V_1 \cup V_2$ is a (4, 3) weighted cyclic set bicode.

The main advantage of these codes is that both error detection is very easy. Now we proceed on to define the notion of dual set bicode.

**DEFINITION 3.1.11:** *Let*

$$C = C_1 \cup C_2 = \left\{C_1^1, C_2^1, \ldots, C_{n_1}^1\right\} \cup \left\{C_1^2, C_2^2, \ldots, C_{n_2}^2\right\}$$

*be a set bicode. The dual set bicode of C or the perpendicular set bicode of C denoted by*

$$C^{\perp} = (C_1 \cup C_2)^{\perp} = (C_1^{\perp} \cup C_2^{\perp})$$

$$= \left\{\left(C_1^1\right)^{\perp}, \ldots, \left(C_{n_1}^1\right)^{\perp}\right\} \cup \left\{\left(C_1^2\right)^{\perp}, \left(C_2^2\right)^{\perp}, \ldots, \left(C_{n_2}^2\right)^{\perp}\right\}$$

*where* $\left(C_{t_i}^{i}\right)^{\perp} = \left\{x^i / x^i \cdot y_{t_i}^{i} = (0)\ \text{for desired}\ y_{t_i}^{i} \in C_{t_i}^{i}\right\}$*; true for* $1 \le t_i \le n_i$*; i = 1, 2. The desired needed code words may depend on weight or cyclic nature etc.*

We illustrate this first by an example.

***Example 3.1.17:*** Let C = $C_1 \cup C_2$ = { $C_1^1$ = {(1 1 1 0 0), (1 1 1 1 1), (0 0 0 0 0), (1 0 1 0 1)}, $C_2^1$ = {(1 1 1 1 1 1), (0 0 0 0 0 0), (1 1 0 0 1 1)}, $C_3^1$ = {(0 0 0 0 0 0 0), (1 1 0 0 0 0 1), (0 0 1 1 1 0 0), (0 1 0 0 1 0 1)}} $\cup$ { $C_1^2$ = {(1 1 0 0 1), (0 1 1 0 0), (0 0 0 0 0), (0 1 1 1 1)}, $C_2^2$ = {(0 0 0 0 0 0), (1 1 1 0 0 0), (1 1 0 0 1 1), (0 0 1 1 0 1)}, $C_3^2$ = {(1 1 1 1 1 1), (0 0 0 0 0 0 0)}, $C_4^2$ = {(1 1 1 0 0 0 0 0), (0 0 0 0 0 0 0 0), (0 0 1 0 0 1 1 1), (1 1 0 0 0 1 1 0), (1 1 0 0 0 1 1 1)}} be a set bicode over the set S = {0, 1}.

The dual set bicode of C denoted by $C^{\perp} = (C_1^{\perp} \cup C_2^{\perp})$ = { $\left(C_1^1\right)^{\perp}$ = {(0 0 0 0 0), (1 0 1 0 0)}, $\left(C_2^1\right)^{\perp}$ = {(0 0 0 0 0 0), (1 1 1 1 1 1), (1 1 0 0 1 1), (0 0 1 1 0 0), (1 1 1 1 0 0), (0 0 1 1 1 1), (0 1 1 1 1 0), (1 0 0 0 1 0), (0 1 0 0 0 1), (0 1 0 0 1 0), (1 0 0 0 0 1)}, $\left(C_3^1\right)^{\perp}$ = {(0 0 0 0 0 0 0), (0 1 0 0 0 0 1), (0 0 1 1 0 0 0)}} $\cup$ { $\left(C_1^2\right)^{\perp}$ = {(0 0 0 0 0), (0 1 1 1 1), (1 1 1 0 0)}, $\left(C_2^2\right)^{\perp}$ = {(0 0 0 0 0 0), (1 1 0 0 0 0)}, $\left(C_3^2\right)^{\perp}$ = {(0 0 0 0 0 0), (1 1 1 1 1 1), (1 1 0 0 0 0), (0 0 1 1 1 1), (1 1 1 1 0 0), (0 0 1 1 0 0), (1 0 0 0 0 1), (0 1 1 1 1 0), (1 0 1 0 1 1), (1 0 1 0 0 0), (1 0 1 0 1 1), (1 0 1 0 1 1)}, $\left(C_4^2\right)^{\perp}$ = {(0 0 0 0 0 0 0 0), (1 1 0 0 0 0 0 0), (0 0 0 0 0 1 1 0), (1 1 0 0 0 1 1 0)}}.

Thus we see $C^{\perp} = (C_1^{\perp} \cup C_2^{\perp})$, is the set dual bicode of the given set code C = $(C_1 \cup C_2)$, .

*Note:* We have just given a set such that $CC^{\perp} = (C_1 \cap C_1^{\perp}) \cup (C_2 \cap C_2^{\perp}) = 0 \cup 0$.

The error detection and correction would be easy if we use the set dual codes.

Now we give a new class of set bicodes called complementing set bicodes.

**DEFINITION 3.1.11:** *Let C = ($C_1 \cup C_1^{\perp}$), be a set bicode i.e., if C =* $\left\{C_1^1, C_2^1, \ldots, C_{n_1}^1\right\} \cup \left\{C_1^2, C_2^2, \ldots, C_n^2\right\}$ *then* $\left(C_i^2\right) = \left(C_i^1\right)^{\perp}$ *i.e., C =* $\left\{C_1^1, C_2^1, \ldots, C_{n_1}^1\right\} \cup \left(\left(C_1^1\right)^{\perp}, \left(C_2^1\right)^{\perp}, \ldots, \left(C_{n_1}^1\right)^{\perp}\right)$ *is called the complementing set bicodes.*

*Note:* It is important to note that in C = ($C_1 \cup C_1^{\perp}$) we don't take all the set codes which are dual with $C_1$.

We illustrate this by the following example.

***Example 3.1.18:*** Let C = ($C_1 \cup C_1^{\perp}$) = {{ $C_1^1$ = {(0 0 0 0 0), (1 1 0 0 0), (1 0 0 0 1), (0 0 1 1 0)}, $C_2^1$ = {(0 0 0 0 0 0), (1 1 1 1 1 1)}, $C_3^1$ = {(0 0 0 0 0 0 0), (0 1 1 1 1 1 1), (0 0 1 1 1 1 0)}, $C_4^1$ = {(0 0 0 0 0 0 0 0), (1 1 0 1 1 0 1 1), (1 1 0 1 1 0 0 0), (0 0 0 0 0 0 1 1), (1 1 0 0 0 0 0 0)}} $\cup$ { $\left(C_1^1\right)^{\perp}$ ={(0 0 0 0 0), (0 0 1 1 0), (1 1 0 0 1), (0 0 1 1 0)}, $\left(C_2^1\right)^{\perp}$ = {(0 0 0 0 0), (1 1 0 0 0 0), (0 1 1 0 0 0), (0 0 1 1 0 0), (0 0 0 0 1 1), (1 1 1 1 0 0), (1 1 1 1 1 1), (1 1 1 1 0 0), (1 1 0 0 1 1)}, $\left(C_3^1\right)^{\perp}$ = {(0 0 0 0 0 0 0), (0 1 1 1 1 1 1), (0 0 1 1 1 1 0), (0 1 1 0 0 0 0), (0 0 0 0 1 1 0), (0 1 1 1 1 1 1)}, $\left(C_4^1\right)^{\perp}$ = {(0 0 0 0 0 0 0 0), (1 1 0 1 1 0 1 1), (0 0 0 0 0 0 1 1), (1 1 0 0 0 0 0 0), (0 0 0 1 1 0 0 0), (1 1 0 1 1 0 0 0), (0 0 0 0 0 1 1 1), (0 0 0 1 1 1 1 1)}} is the complementing set bicode.

It is pertinent to mention here that for a given set code $C_1$ the complementing set bicode got as $C_1^{\perp}$ may be finitely many i.e., it is not unique in general.

We give yet another complementing part of $C_1^{\perp}$ given in the above example.

***Example 3.1.19:*** Let C = ($C_1 \cup C_1^{\perp}$) where {$C_1$ = { $C_1^1$ = {(0 0 0 0 0), (1 1 0 0 0), (1 0 0 0 1), (0 0 1 1 0)}, $C_2^1$ = {(0 0 0 0 0 0), (1 1 1 1 1 1)}, $C_3^1$ = {(0 0 0 0 0 0 0), (0 1 1 1 1 1 1), (0 0 1 1 1 1 0)}, $C_4^1$ = {(0 0 0 0 0 0 0 0), (1 1 0 1 1 0 1 1), (1 1 0 1 1 0 0 0), (0 0 0 0 0 0 1 1), (1 1 0 0 0 0 0 0)}} ∪ { $\left(C_1^1\right)^{\perp}$ = {(0 0 0 0 0), (1 1 0 0 1)}, $\left(C_2^1\right)^{\perp}$ = {(0 0 0 0 0 0), (1 1 0 0 0 0), (0 0 1 1 1 1), (1 1 1 1 1 1), (0 0 1 1 0 0)}, $\left(C_3^1\right)^{\perp}$ = {(0 0 0 0 0 0 0), (0 0 1 1 0 0 0), (0 0 0 1 1 0 0), (0 0 0 0 1 1 0)}, $\left(C_4^1\right)^{\perp}$ = {(0 0 0 0 0 0 0 0), (1 1 0 0 0 0 1 1), (0 0 0 1 1 0 0 0)}}.

C = ($C_1 \cup C_1^{\perp}$) is also a complementing set bicode. Clearly C = ($C_1 \cup C_1^{\perp}$) given in example 3.1.18 is different from the complementing set bicode given here.

The main advantage of these classes of bicodes is that they are very much useful in error correction and error detection. Now to have more advantage than these bicodes we define (m, m), weighted complementary set bicodes and (m, n), (m ≠ n), weighted complementary set bicodes.

**DEFINITION 3.1.12:** *Let*

$$V = V_1 \cup V_2 = \left\{V_1^1, V_2^1, \ldots, V_n^1\right\} \cup \left\{\left(V_1^1\right)^{\perp}, \left(V_2^1\right)^{\perp}, \ldots, \left(V_n^1\right)^{\perp}\right\}$$

*be a complementary set bicode were the weight of each code word in every $V_j^1$ and $\left(V_j^1\right)^{\perp}$ are m; $1 \le j \le n$. Then we call V to be a (m, m) weighted complementary set bicode.*

The main advantage of these new classes of codes is that they are useful for error detection and error correction; further these codes can be used by cryptologists so that it cannot be easily broken by an intruder. These also can be used in channels with varying lengths but with same weight.

Now we illustrate this new classes of codes by the following examples.

***Example 3.1.20:*** Let

$$V = V_1 \cup V_1^{\perp} = \{V_1^1, V_2^1, V_3^1, V_4^1\} \cup \{(V_1^1)^{\perp}, (V_2^1)^{\perp}, (V_3^1)^{\perp}, (V_4^1)^{\perp}\}$$

where $V_1^1$ = {(0 0 0 0 0 0), (1 1 0 1 1 0 0), (1 1 1 1 0 0), (1 1 0 0 1 1), (0 1 1 1 1 0)}, $V_2^1$ = {(0 0 0 0 0 0 0), (1 0 1 0 1 0 1), (1 1 0 0 1 1 0), (0 1 1 0 0 1 1), (1 1 1 1 0 0 0), (0 1 1 1 1 0 0), (1 1 1 0 0 0 1), (0 0 0 1 1 1 1), (0 0 1 1 1 1 0)}, $V_3^1$ = {(0 0 0 0 0), (1 1 1 1 0), (0 1 1 1 1), (1 1 0 1 1)} and $V_4^1$ = {(1 1 1 1 0 0 0 0), (0 1 1 1 1 0 0 0), (0 0 0 0 0 0 0 0), (0 0 1 1 1 1 0 0), (1 1 0 0 0 0 1 1), (1 1 0 0 1 1 0 0)}, $(V_1^1)^{\perp}$ = {(0 0 0 0 0 0), (1 1 0 0 1 1)}, $(V_2^1)^{\perp}$ = {(0 0 0 0 0 0 0), (1 1 0 0 1 1 0)}, $(V_3^1)^{\perp}$ = {(0 0 0 0 0)}, $(V_4^1)^{\perp}$ = {(0 0 0 0 0 0 0 0), (0 0 1 1 0 0 1 1)}. We see V is a (4, 4) weighted complementary set bicode.

These codes can be used in cryptography as well as in channels were retransmission is not possible as the corrected bicode word can easily be obtained based on both weight and duality.
Now we proceed on to define (m, n) weighted complementary set bicodes.

**DEFINITION 3.1.13:** *Let V = $V_1 \cup V^{\perp}$ be a complementary set bicode given by V = $V_1 \cup V_2$ = $\{V_1^1, \ldots, V_n^1\} \cup \{(V_1^1)^{\perp}, \ldots, (V_n^1)^{\perp}\}$ where weight of each code word in each $V_i^1$ in $V_1$ is of weight m, m $\neq$ n; 1 $\leq$ i $\leq$ n and that is each code word in $(V_i^1)^{\perp}$ is of weight n (m $\neq$ n), 1 $\leq$ i $\leq$ n. Then we call this new set bicode as (m, n) weighted complementary set bicode.*

We illustrate this by a simple example.

***Example 3.1.21:*** Let V = $V_1 \cup V_1^{\perp}$ where $\{V_1^1, V_2^1, V_3^1, V_4^1\}$ with $V_1^1$ = {(0 0 1 1 1), (1 1 1 0 0), (0 0 0 0 0)}, $V_2^1$ = {(0 0 0 1 1 1), (0 1 1 1 0 0), (0 0 1 1 1 0), (0 0 0 0 0 0)}, $V_3^1$ = {(0 1 1 1 0 0 0), (0 0 1 1 1 0 0), (0 0 0 0 0 0 0), (0 0 0 0 1 1 1), (1 1 0 0 0 1 0)}, and $V_4^1$ = {(0 0 0 0 0 1 1 1), (0 0 0 0 0 0 0 0), (0 0 1 1 1 0 0 0), (1 1 0 0 0 0 0 1)} are codes in $V_1$ of weight 3. Now we wish to take codes in $V_1^{\perp}$ to be of weight 4 which is as follows; $V_1^{\perp}$ = $\{(V_1^1)^{\perp}, (V_2^1)^{\perp}, (V_3^1)^{\perp}, (V_4^1)^{\perp}\}$ where $(V_1^1)^{\perp}$ = {(0 0 0 0 0), (1 1 0 1 1)}, $(V_2^1)^{\perp}$ = {(0 0 0 0 0 0), (0 1 1 0 1 1), (1 1 0 1 1 0), (1 0 1 1 0 1)}, $(V_3^1)^{\perp}$ = {(0 0 0 0 0 0 0), (0 1 1 0 1 1 0)} and $(V_4^1)^{\perp}$ = {(0 0 0 0 0 0 0 0), (1 1 1 1 0 0 0 0), (0 0 0 1 1 1 1 0), (1 1 0 1 1 0 0 0), (0 0 1 1 0 1 1 0)}. Clearly V= $V_1 \cup V_1^{\perp}$ is a (3, 4) weighted complimentary set bicode.

The main advantage of these codes are easy error detection and error correction and difficult to be hacked by an intruder when any message is transmitted for it is not easy to guess, the real code which carries the message; only the receiver at the other end and the sender exactly know which of the code that carries the message. Another advantage is the user need not be very well versed in coding theory. Further these codes can be used when one need different lengths of messages to be sent across with varying number of check symbols.

Error detection is easy after error detection error correction or guessing can be easily made when we use (m, n) or (m, m) weighted complementary set bicodes.

Next we proceed on to define yet another new notion namely semigroup bicodes.

**DEFINITION 3.1.14:** *Let V = $V_1 \cup V_2$ = {$X_1, X_2, \ldots, X_{n_1}$ } $\cup${$Y_1, Y_2, \ldots, Y_{n_2}$ } be set bicode over the semigroup S = {0, 1}. If each of the codes of the same length of the set codes in $V_1$ and $V_2$ form a semigroup under addition i.e., a monoid under addition i.e., V = $\left\{S_1^1,\ldots,S_{r_1}^1\right\} \cup \left\{\left(S_1^2\right),\ldots,\left(S_{r_2}^2\right)\right\}$ such that $r_1 < n_1$ and $r_2 < n_2$ then we call V to be a semigroup bicode. The elements of V are called semigroup bicode words.*

We illustrate this situation by some examples.

***Example 3.1.22:*** Let S = {$V_1 \cup V_2$} = { $V_1^1$ = {(0 0 0 0 0), (1 1 0 1 0), (0 1 0 1 1), (1 0 0 0 1)}, $V_2^1$ = (0 0 0 0 0 0), (1 1 1 0 0 0), (0 0 0 1 1 1), (1 1 1 1 1 1)}, $V_3^1$ = {(1 1 0 0 0 1 1), (0 0 1 1 1 0 0), (1 1 1 1 1 1 1), (0 0 0 0 0 0 0)}, $V_4^1$ = {(1 1 1 1), (0 0 0 0), (1 0 0 0), (0 1 1 1)}} $\cup$ { $V_1^2$ = {(0 0 0 0 0 0), (1 1 0 0 1 1), (0 0 1 1 0 0), (1 1 1 1 0 0), (0 0 0 0 1 1), (0 0 1 1 1 1), (1 1 1 1 1 1), (1 1 0 0 0 0)}, $V_2^2$ = {(0 0 0 0 0 0 0 0), (1 1 0 0 1 1 0 0), (1 1 0 0 1 1 0 0), (0 0 1 1 0 0 1 1), (1 1 1 1 0 0 0 0), (0 0 0 0 1 1 1 1), (0 0 1 1 1 1 0 0), (1 1 1 1 1 1 1 1)}, $V_3^2$ = {(1 1 1 0 0 0 0), (0 0 0 0 0 0 0), (0 0 0 1 1 1 0), (1 1 1 1 1 1 0)}}. Clearly V = $V_1 \cup V_2$ is a semigroup bicode.

***Example 3.1.23:*** Let

$$V = V_1 \cup V_2 = \left\{V_1^1, V_2^1, V_3^1, V_4^1\right\} \cup \left\{V_1^2, V_2^2, V_3^2, V_4^2\right\}$$

where $V_1^1$ = {(1 1 1 1 1 1 1), (0 0 0 0 0 0 0), (1 1 0 0 1 1 0), (0 0 1 1 0 0 1)}, $V_2^1$ = {(1 1 1 1 1), (0 0 0 0 0), (1 1 0 0 1), (0 0 1 1 0)}, $V_3^1$ = {(1 1 1 1 1 1), (0 0 0 0 0 0)}, $V_4^1$ = {(1 1 1 1 1 1 1 1), (1 1 0 0 0 0 0 0), (0 0 1 1 1 0 0 0), (0 0 0 0 0 1 1 1), (1 1 1 1 1 0 0 0), (0 0 1 1 1 1 1 1), (0 0 0 0 0 0 0 0)}, $V_1^2$ = {(1 1 1 1 1), (0 0 0 0 0), (1 1 1 0 0), (1 0 0 0 0), (0 1 1 0 0), (1 0 0 1 1), (0 1 1 1 1), (0 0 0 1 1)}, $V_2^2$ = {(0 0 0 0 0 0), (1 0 1 0 1 0), (0 1 0 1 0 1),

(1 1 1 1 1 1)}, $V_3^2$ = {(0 0 0 0 0 0 0), (1 1 0 0 0 0 0), (0 0 0 0 0 1 1), (1 1 0 0 0 1 1)}. It is easily verified V is a semigroup bicode.

We proceed onto prove the following theorem.

**THEOREM 3.1.1:** *Every set repetition bicode is a semigroup bicode.*

*Proof:* Let C = $C_1 \cup C_2$ = $\left\{C_1^1, C_2^1, \cdots, C_{n_1}^1\right\} \cup \left\{C_1^2, C_2^2, \cdots, C_{n_2}^2\right\}$ be a set repetition bicode; clearly each $C_{t_i}^i$ is a repetition code, ie $C_{t_i}^i$ = {(0 0 … 0), (1 1 1 … 1)}; $1 \le t_i \le n_i$ ; i = 1, 2. Clearly each $C_{t_i}^i$ is a semigroup code for every i. Hence C = $C_1 \cup C_2$ is a semigroup bicode.

**THEOREM 3.1.2:** *Every semigroup bicode is a set bicode but a set bicode in general is not a semigroup bicode.*

*Proof:* Let C = $C_1 \cup C_2$ be a semigroup bicode as every semigroup is a set we see C = $C_1 \cup C_2$ is a set bicode. On the other hand every set bicode need not in general be a semigroup bicode. We prove this by a counter example. Consider a set bicode C = $C_1 \cup C_2$ = $\left\{C_1^1, C_2^1, C_3^1\right\} \cup \left\{C_1^2, C_2^2, C_3^2, C_4^2\right\}$ where $C_1^1$ = {(1 1 1 0 0 0 0), (0 0 1 1 0 0 1), (0 0 1 0 1 0 1), (0 1 1 1 0 1 0), (1 0 1 1 1 0 1), (0 0 0 0 0 0 0)}, $C_2^1$ = {(0 0 0 0 0), (1 1 0 0 0), (0 0 0 1 0), (1 1 1 1 1)}, $C_3^1$ = {(0 0 0 0 0 0), (1 1 0 1 1 0), (0 1 1 1 0 1), (0 1 0 1 1 0), (1 0 0 1 0 1)}, $C_1^2$ = {(0 0 0 0 0 0), (1 1 1 1 1 1), (1 1 0 0 1 1), (0 1 1 1 0 0)}, $C_2^2$ = {(0 0 0 0 0 0 0), (1 1 0 0 1 1 0), (1 0 1 1 1 1 0), (0 0 1 1 1 1 0)}, $C_3^2$ = {(0 0 0 0 0 0 0 0), (1 1 0 0 1 1 0 0), (0 0 1 1 1 1 0 0), (1 1 0 0 1 1 1 1)} and $C_4^2$ = {(1 1 1 1 1), (1 1 0 1 1), (1 1 1 0 0), (1 0 1 0 1), (0 0 0 0 0)}.

It is easily verified that the set codes $C_{t_i}^i$ ; $1 \le i \le 2$ are not all semigroup codes. Thus C is a set bicode which is not a semigroup bicode.

Now we proceed on to define (m, m) weight semigroup bicode.

**DEFINITION 3.1.15:** *Let*

$$V = V_1 \cup V_2 = \left\{V_1^1, \ldots, V_{n_1}^1\right\} \cup \left\{V_1^2, V_1^2, \ldots, V_{n_2}^2\right\}$$

*be a semigroup bicode. If every bicode word $X = X_1 \cup X_2$ of V is such that Hamming weight of $X_1$ is equal Hamming weight of $X_2$ equal to m then we call V to be a (m, m) weighted semigroup bicode.*

We illustrate this by an example.

***Example 3.1.24:*** Let

$$V = V_1 \cup V_2 = \left\{V_1^1, V_2^1, V_3^1\right\} \cup \left\{V_1^2, V_2^2\right\}$$

be a semigroup bicode were

$$V_1^1 = \{(0\ 0\ 0\ 0\ 0\ 0), (0\ 0\ 0\ 1\ 1\ 1)\},$$

$$V_2^1 = \{(0\ 0\ 0\ 0\ 0), (0\ 1\ 1\ 1\ 0)\},$$

$$V_3^1 = \{(0\ 0\ 0\ 0\ 0\ 0\ 0), (1\ 1\ 1\ 0\ 0\ 0\ 0)\},$$

$$V_1^2 = \{(0\ 0\ 0\ 0\ 0\ 0), (1\ 0\ 1\ 0\ 1\ 0)\}$$

and

$$V_2^2 = \{(0\ 0\ 0\ 0\ 0\ 0\ 0), (0\ 1\ 0\ 1\ 0\ 1\ 0)\}.$$

Clearly V is a (3, 3) weighted semigroup bicode.

We give yet another example.

***Example 3.1.25:*** Let

$$V = V_1 \cup V_2 = \left\{V_1^1, V_2^1\right\} \cup \left\{V_1^2, V_2^2, V_3^2\right\}$$

be a semigroup bicode over the set S = {0, 1} where $V_1^1$ = {(0 0 0 0 0 0), (1 1 0 0 1 1), (0 0 1 1 0 1 1), (1 1 1 1 0 0)}, $V_2^1$ = {(0 0 0 0 0 0 0 0), (1 1 0 1 1 0 0 0), (0 0 0 1 1 1 1 0), (1 1 0 0 0 1 1 0)}, $V_1^2$ = {(0 0 0 0 0), (1 1 1 1 0)}, $V_2^2$ = {(0 0 0 0 0 0 0), (1 1 0 1 1 0 0), (0 0 0 1 1 1 1), (1 1 0 0 0 1 1)}, $V_3^2$ = {(0 0 0 0 0 0),

(0 0 1 1 1 1), (1 1 0 0 1 1), (1 1 1 1 0 0)}. V is a (4, 4) weighted semigroup bicode.

We have the following interesting theorem.

**THEOREM 3.1.3:** *Let*

$$V = V_1 \cup V_2 = \left\{V_1^1, \cdots, V_{n_1}^1\right\} \cup \left\{V_1^2, V_1^2, \cdots, V_{n_2}^2\right\}$$

*be a repetition semigroup bicode, V is not a (m, m) weight semigroup bicode.*

*Proof:* Clearly each code in $V_1 = \left\{V_1^1, V_2^1, \cdots, V_{n_1}^1\right\}$ is a repetition code and each must of different length. Similarly each code in $V_2$ is a repetition code of distinct length. Since in any repetition code the weight of every non zero code word is the same as its length.

We see no two codes in $V_i$ have same length i = 1, 2. Hence a repetition semigroup bicode is never a (m, m), weight semigroup bicode.

We have the converse also to be true for we have a (m, m), weight semigroup bicode is never a repetition semigroup bicode.

The proof is left as an exercise for the reader to prove.

We just at this juncture give an example of a repetition semigroup bicode.

***Example 3.1.26:*** Let

$$V = V_1 \cup V_2 = \left\{V_1^1, V_2^1, V_3^1, V_4^1\right\} \cup \left\{V_1^2, V_2^2, V_3^2, V_4^2, V_5^2\right\}$$

be a repetition semigroup bicode where $V_1^1$ = {(1 1 1 1 1), (0 0 0 0 0)}, $V_2^1$ = {(1 1 1 1 1 1), (0 0 0 0 0 0)}, $V_3^1$ = {(0 0 0 0 0 0 0), (1 1 1 1 1 1 1)}, $V_4^1$ = {(1 1 1 1), (0 0 0 0)}, $V_1^2$ = (1 1 1 1 1 1 1 1), (0 0 0 0 0 0 0 0)}, $V_2^2$ = {(1 1 1 1 1), (0 0 0 0 0)}, $V_3^2$ = {(0 0 0), (1 1 1)}, $V_4^2$ = {(1 1 1 1 1 1 1 1 1), (0 0 0 0 0 0 0 0 0)} and $V_5^2$ = {(1 1 1 1 1 1), (0 0 0 0 0 0)}. We see the weights of $V_1^1$ = 5, $V_2^1$ = 6, $V_3^1$ = 7, $V_4^1$ = 4, $V_1^2$ = 8, $V_2^2$ = 5, $V_4^2$ = 9, and $V_5^2 = 6$. Thus it can never be a (m, m), weighted semigroup bicode.

Now we proceed on to define the notion of (m, n) (m≠n), weighted semigroup bicode.

**DEFINITION 3.1.16:** *Let*

$$C = C_1 \cup C_2 = \left\{C_1^1, C_2^1, \cdots, C_{n_1}^1\right\} \cup \left\{C_1^2, C_2^2, \cdots, C_{n_2}^2\right\}$$

*be a semigroup bicode. If every code word in the code $C_j^1$ in $C_1$ is of same weight m for every j = 1, 2, …, $n_1$ and every code word in the code $C_k^2$ in $C_2$ is of the same weight n for every k = 1, 2, …, $n_1$, then we call C to be a (m, n) weight semigroup bicode.*

We illustrate this by some simple examples.

***Example 3.1.27:*** Let

$$C = C_1 \cup C_2 = \left\{C_1^1, C_2^1, C_3^1, C_4^1\right\} \cup \left\{C_1^2, C_2^2, C_3^2\right\}$$

where $C_1^1$ = {(0 0 0 0 0 0 0 0), (1 1 1 0 0 0)}, $C_2^1$ = {(0 0 0 0 0 0 0 0), (1 1 0 0 0 0 0 1)}, $C_3^1$ = {(0 0 0 0), (1 1 1 0)}, $C_4^1$ = {(0 0 0 0 0 0 0), (1 1 1 0 0 0 0)}, $C_1^2$ = (1 1 0 0 1 1), (0 0 0 0 0 0), (1 1 1 1 0 0), (0 0 1 1 1 1)}, $C_2^2$ = {(0 0 0 0 0 0 0 0), (1 1 1 1 0 0 0 0), (0 0 1 1 0 0 1 1), (1 1 0 0 0 0 1 1)} and $C_3^2$ = {(1 1 1 1 0 0 0 0), (0 0 0 0 0 0 0 0), (0 0 1 1 1 1 0 0), (1 1 0 0 1 1 0 0)} is a (3, 4) weighted semigroup bicode.

These bicodes can be used in cryptography as it can mislead the intruder. Secondly these bicodes are such that error is detected easily and error correction is also possible as it holds the condition of being a semigroup. Infact the semigroup bicodes C = $C_1 \cup C_2$ can also be defined or called as bisemigroup codes for C = $C_1 \cup C_2$ is clearly a bisemigroup.

Now we proceed on to define the new notion of orthogonal semigroup bicode.

**DEFINITION 3.1.17:** *Let*

$$C = C_1 \cup C_2 = \left\{C_1^1, C_2^1, \cdots, C_{n_1}^1\right\} \cup \left\{C_1^2, C_2^2, \cdots, C_{n_2}^2\right\}$$

*be a semigroup bicode over the set {0, 1}. The orthogonal semigroup bicode of dual semigroup bicode of C denoted by*

$$\{(C_1 \cup C_2)^{\perp}\} = C_1^{\perp} \cup C_2^{\perp}$$

$$= \left(\left(C_1^1\right)^{\perp`}, \left(C_2^1\right)^{\perp}, \ldots, \left(C_{n_1}^1\right)^{\perp}\right) \cup \left(\left(C_1^2\right)^{\perp`}, \left(C_2^2\right)^{\perp}, \ldots, \left(C_{n_2}^2\right)^{\perp}\right)$$

*is the dual code of each* $C_{t_j}^i$ *where* $\left(C_{t_j}^i\right)^{\perp}$ *= {x | x.y = 0 for all y* $\in$ $C_{t_j}^i$*, such that the collection forms a semigroup under addition}, true for* $1 \le t_j \le n_i$*, i = 1, 2.*

We illustrate this by a simple example.

***Example 3.1.28:*** Let

$$C = C_1 \cup C_2 = \left\{C_1^1, C_2^1, C_3^1\right\} \cup \left\{C_1^2, C_2^2, C_3^2, C_4^2\right\}$$

where $C_1^1$ = {(0 0 0 0 0), (1 1 1 0 0), (1 1 1 1 1), (0 0 0 1 1)}, $C_2^1$ = {(0 0 0 0 0 0), (1 1 0 0 0 0), (1 1 1 1 0 0), (0 0 1 1 0 0), (0 0 0 0 1 1), (1 1 1 1 1 1)}, $C_3^1$ = {(0 0 0 0 0 0 0), (1 1 1 0 0 0 0), (0 0 0 1 1 0 0), (0 0 0 0 0 1 1), (1 1 1 0 0 1 1), (0 0 0 1 1 1 1), (1 1 1 1 1 0 0)}, $C_1^2$ = {(0 0 0 0 0), (1 0 1 0 1), (1 1 1 0 0), (0 0 0 1 1), (0 1 0 0 1), (1 0 1 1 0), (1 1 1 1 1), (0 1 0 1 0)}, $C_2^2$ = {(0 0 0 0 0 0), (1 1 1 1 1 1)}, $C_3^2$ = {(0 0 0 0 0 0 0), (1 1 1 1 0 0 0), (1 1 1 1 1 1 1), (0 0 0 0 1 1 1)} and $C_4^2$ = {(0 0 0 0 0 0 0 0), (1 1 0 0 1 1 0 0), (0 0 1 1 0 0 1 1), (1 1 1 1 0 0 0 0), (1 1 1 1 1 1 1 1), (0 0 1 1 1 1 0 0), (0 0 0 0 1 1 1 1)} be a semigroup bicode.

$C^{\perp}$ = $C_1^{\perp} \cup C_2^{\perp}$ = {{(0 0 0 0 0), (0 1 1 0 0), (0 0 0 1 1), (0 1 1 1 1)} = $\left(C_1^1\right)^{\perp}$ }, $\left(C_2^1\right)^{\perp}$ = {(0 0 0 0 0 0), (1 1 0 0 0 0), (1 1 1 1 0

0), (0 0 1 1 0 0)}, $\left(C_3^1\right)^\perp$ = {(0 0 0 0 0 0 0), (1 1 0 0 0 0 0), (0 0 0 0 0 1 1), (1 1 0 0 0 1 1)}} ∪ {$\left(C_1^2\right)^\perp$ = {(0 0 0 0 0), (1 0 1 0 0)}, $\left(C_2^2\right)^\perp$ = {(0 0 0 0 0 0), (1 1 0 0 0 0), (0 0 1 1 0 0), (1 1 1 1 0 0), (1 1 0 0 1 1), (0 0 1 1 1 1), (1 1 1 1 1 1)}, $\left(C_3^2\right)^\perp$ = {(0 0 0 0 0 0 0), (1 1 0 0 0 0 0), (0 0 1 1 0 0 0), (1 1 1 1 0 0 0)}, $\left(C_4^2\right)^\perp$ = {(0 0 0 0 0 0 0 0), (1 1 1 1 1 1 1 1), (1 1 1 1 0 0 0 0), (0 0 0 0 1 1 1 1)}}. $C^\perp = C_1^\perp \cup C_2^\perp$ is the dual semigroup bicode.

We see we can always restrain ourselves to obtain dual bicode of a dual semigroup bicode to be a semigroup bicode.

Now we just see any message X in C = $C_1 \cup C_2$ would be of the form

$$X = \left(X_1^1, X_2^1, \ldots, X_{n_1}^1\right) \cup \left(X_1^2, X_2^2, \ldots, X_{n_2}^2\right)$$

where $X_{j_i}^i \in C_{j_i}^i$ ; $1 \le j_i \le n_i$; i = 1, 2.

We illustrate this by the following example.

***Example 3.1.29:*** Let

$$C = C_1 \cup C_2 = \left\{C_1^1, C_2^1, C_3^1\right\} \cup \left\{C_1^2, C_2^2, C_3^2, C_4^2\right\}$$

be a semigroup bicode where

$C_1^1$ = {(0 0 0 0 0), (0 0 1 1 0), (1 1 0 0 1)},

$C_2^1$ = {(0 0 0 0 0 0), (0 0 0 1 1 0), (1 1 0 1 1 0), (1 1 0 0 0 0)},

$C_3^1$ = {(0 0 0 0 0 0 0), (1 1 0 1 0 0 0), (0 0 1 0 1 1 1), (1 1 1 1 1 1 1)},

$C_1^2$ = {(0 0 0 0 0 0 0 0), (1 1 1 1 0 0 0 0), (0 0 0 0 1 1 1 1), (1 1 0 0 0 0 1 1), (0 0 1 1 0 0 1 1), (0 0 1 1 1 1 0 0), (1 1 1 1 1 1 1 1)},

$C_2^2$ = {(0 0 0 0 0 0), (1 1 1 1 1 1), (1 1 1 0 0 0), (0 0 0 1 1 1)},

$C_3^2$ = {(0 0 0 0 0), (1 1 0 0 0), (0 0 1 1 1)} and

$C_4^2$ = {(0 0 0 0 0 0 0 0 0 0), (1 1 1 1 1 0 0 0 0 0), (0 0 0 0 0 1 1 1 0 0), (1 1 1 1 1 1 1 1 0 0) (0 0 0 0 0 0 0 0 1 1), (1 1 1 1 1 1 1 1 1 1), (0 0 0 0 0 1 1 1 1 1)}

be a bisemigroup code (or a semigroup bicode) over the set S = {0, 1}. Any element X = $X_1 \cup X_2$ = {(1 1 0 0 1), (1 1 0 0 0 0), (1 1 1 1 1 1 1)} $\cup$ {(0 0 1 1 1 1 0 0), (0 0 0 1 1 1), (1 1 0 0 0), (0 0 0 0 1 1 1 1 1)} $\in C_1 \cup C_2$.

That is they can send at a time by a single transmission 7 set of messages of length (5, 6, 7) $\cup$ (8, 6, 5, 1 0).

If due to some condition some of the receivers are not available then they can send message Y = $Y_1 \cup Y_2$ = {(1 1 0 0 1), (1 1 0 1 1 0)} $\cup$ {(1 1 0 0 0 0 1 1), (0 0 1 1 1), (1 1 1 1 1 1 1 1 0 0)} i.e., 5 set of messages of lengths (5, 6) $\cup$ (8, 5, 1 0) .

The sender can also send only 2 messages of Z = $Z_1 \cup Z_2$ = {(1 1 0 1 1 0)} $\cup$ {(1 1 1 1 1 1 1 1 1 1)} or even a single message as {(1 1 1 1 1 1 1)} $\cup$ $\phi$ or $\phi$ $\cup$ {(1 1 1 1 1 1 1 1 0 0)}.

Thus the flexibility of the use of channels and the number of messages sent at time makes this semigroup bicodes advantageous over other codes. This is true even in case of set bicodes.

We illustrate the typical communication system.

The channel would be a multichannel which can receive and transmit at a time a maximum of $n_1 + n_2$ simultaneously i.e., a $n_1 + n_2$ channel.

According to need a few channel need not transmit remain in the off state. Such sort of transmission coding and decoding is possible due to the advent of supercomputer and proper programming.

Thus these multichannel in time of need can also function as a single channel. The advantage of these codes is they cannot be hacked easily while transmitting secret messages they will transmit $n_1 + n_2$ code words but only a few of them will be really carrying the messages rest will be misleading codes. So at the receiving end the receiver will only decode the codes which carry the messages and ignore the misleading codes totally. By this method of transmission it is impossible for the intruder to easily hack the message or even guess which codes really carry the messages!

Apart from cryptography these codes are best suited for multichannel transmission that too when the lengths of the messages are different.

When we use the notion of (m, n), weighted semigroup bicodes the error detection and correction can also be carried out easily.

These codes can be used in computers, televisions and also in counter when several different sets of messages are to be transmitted simultaneously.

Now we proceed on to describe the new notion of semigroup cyclic bicode.

**DEFINITION 3.1.18:** *Let*

$$C = \left(C_1^1, ..., C_{n_1}^1\right) \cup \left(C_1^2, ..., C_{n_2}^2\right)$$

*be a semigroup bicode, if each* $C_{t_j}^i$ *is a cyclic code for* $1 \leq t_j \leq n_i$; $1 \leq i \leq 2$ *then we call S to be a cyclic semigroup bicode or semigroup cyclic bicode.*

**THEOREM 3.1.4:** *Every repetition semigroup bicode is a cyclic semigroup bicode.*

*Proof:* Given

$$C = C_1 \cup C_2 = \left(C_1^1, ..., C_{n_1}^1\right) \cup \left(C_1^2, ..., C_{n_2}^2\right)$$

is a repetition semigroup bicode. Thus each $C_{j_i}^i$ = {(0 0 … 0), (1 1 … 1)}, $1 \leq j_i \leq n_i$, $1 \leq i \leq 2$; is a cyclic code, hence C is a cyclic semigroup bicode.

Further the class of set bicodes contains the class of semigroup bicodes.

We can as in case of semigroup codes use the method of approximation described in chapter one to find the correct message. We define group bicode and illustrate it by examples.

**DEFINITION 3.1.19:** *Let*

$$C = C_1 \cup C_2 = \left\{C_1^1, C_2^1, \cdots C_{n_1}^1\right\} \cup \left\{C_1^2, C_2^2, \cdots C_{n_2}^2\right\}$$

*where each* $C_{t_i}^i$ *is the collection of code words which forms a group under addition;* $1 \leq t_j \leq n_i$, $i = 1, 2$; *we call C to be a group bicode or special bigroup code, i.e., C is a group bicode if each* $C_i$ *is a group code.*

We illustrate this by some simple examples.

***Example 3.1.30:*** Let

$$C = C_1 \cup C_2 = \left\{C_1^1, C_2^1, C_3^1, C_4^1, C_5^1\right\} \cup \left\{C_1^2, C_2^2, C_3^2\right\}$$

where

$C_1^1$ = {(0 0 0 0 0) (1 1 1 1 1) (1 0 0 1 1) (0 1 1 0 0) (0 0 1 1 1) (1 1 0 0 0) (0 1 0 1 1) (1 0 0 1 0)}

$C_2^1$ = {(0 0 0 0) (1 0 1 1) (0 1 0 1) (1 1 1 0)}

$C_3^1$ = {(0 0 0 0 0 0) (0 0 1 0 0 1) (1 0 0 1 0 0) (0 1 0 0 1 0) (0 1 1 0 1 1) (1 1 1 1 1 1) (1 1 0 1 1 0) (1 0 1 1 0 1)}

$C_4^1$ = {(0 0 0 0 0 0 0) (1 1 0 1 0 0 0) (0 1 1 0 1 0 0) (0 0 1 1 0 1 0) (0 0 0 1 1 0 1) (1 0 1 1 1 0 0) (1 1 1 0 0 1 0) (0 1 0 1 1 1 0) (1 0 0 0 1 1 0) (0 1 0 0 0 1 1) (1 0 0 1 0 1 1) (1 0 1 0 1 0 1) (1 1 1 1 0 0 1) (0 0 1 0 1 1 1) (1 1 0 0 1 0 1) (0 1 1 1 0 0 1)}

$C_5^1$ = {(1 0 0 0 1 0 1) (0 1 0 0 1 1 1) (0 0 1 0 1 1 0) (0 0 0 1 0 1 1) (1 1 0 0 0 1 0) (0 1 1 0 0 0 1) (0 0 1 1 1 0 1) (1 0 0 1 1 1 0) (1 0 1 0 0 1 1) (0 1 0 1 1 0 0) (1 1 1 0 1 0 0) (1 1 0 1 0 1 1) (1 0 1 1 0 0 0) (0 1 1 1 0 1 0) (0 0 0 0 0 0 0)(1 1 1 1 1 1 1)}

$C_1^2$ = {(0 0 0 0 0 0) (1 0 0 1 0 0) (0 1 0 0 1 0) (0 0 1 0 0 1) (1 1 0 1 1 0) (0 1 1 0 1 1) (1 1 1 1 1 1) (1 0 1 1 0 1)}

$C_2^2$ = {(0 0 0 0 0 0 0) (1 1 1 0 1 0 0) (0 1 1 1 0 1 0) (0 0 1 1 1 0 1) (1 0 0 1 1 1 0) (1 1 0 1 0 0 1) (0 1 0 0 1 1 1) (1 0 1 0 0 1 1)} and

$C_3^2$ = {(0 0 0 0 0) (1 1 1 0 0) (0 0 1 1 0) (1 1 1 1 1 1) (1 1 0 1 0) (0 0 1 0 1) (0 0 0 1 1) (1 1 0 0 1)}

is a group bicode. We see each code $C^i_{t_i}$ is a usual code for $1 \le t_i \le n_i$; i = 1, 2. Now we make the following interesting observation.

1. All group bicodes are set bicodes as well as semigroup bicodes.
2. A set bicode in general is not a group bicode.
3. A semigroup bicode need not be a group bicode.
4. When we use group bicode the main advantage being that both error detection and error correction can be done in a easy way.

The advantage of set bicode over group bicode is that we can in the set bicode

$$C = C_1 \cup C_2 = \left\{C^1_1, C^1_2, \ldots, C^1_{n_1}\right\} \cup \left\{C^2_1, C^2_2, \ldots, C^2_{n_2}\right\},$$

take one of the $C^i_{t_i}$ 's to be a usual code (general linear binary code) and use these codes in cryptography or in defence department for the sender and the receiver knows very well that the code which carries the messages in $C = C_1 \cup C_2$ and can code and decode the message for all other codes sent are just to mislead the intruder.

We call these bicodes in which one or more codes $C^i_{t_i}$ are usual codes as set bicodes only.

We shall illustrate this by a simple example, how these codes can preserve confidentiality.

***Example 3.1.31:*** Let

$$C = C_1 \cup C_2 = \left(C^1_1, C^1_2, C^1_3\right) \cup \left(C^2_1, C^2_2, C^2_3, C^2_4\right)$$

where

$C^1_1$ = {(0 0 0 0 0 0) (1 1 0 0 1 0 0) (1 0 1 1 1 1) (1 0 1 0 1 0)}

$C^1_2$ = {(0 0 0 0 0) (1 1 1 1 1) (1 1 0 0 1) (1 1 0 0 1) (0 1 1 0 1)}

$C^1_3$ = {(0 0 0 0 0 0 0) (1 1 1 0 1 0 0) (0 1 1 1 0 1 0) (0 0 1 1 1 0 1) (1 0 0 1 1 1 0) (1 1 0 1 0 0 1) (0 1 0 0 1 1 1) (1 0 1 0 0 1 1)}

$C_1^2$ = {(0 0 0 0 0 0 0) (0 1 1 1 0 1 0) (1 1 1 0 1 0 0) (0 1 0 0 1 1 1) (1 1 1 1 1 1 1)}

$C_2^2$ = {(0 0 0 0 0 0) (0 0 1 0 0 1) (0 1 0 0 1 0) (0 1 1 0 1 1) (1 0 0 1 0 1) (1 0 1 1 0 1) (1 1 0 1 1 0) (1 1 1 1 1 1)}

$C_3^2$ = {(0 0 0 0 0) (1 1 0 0 0) (1 0 0 0 1) (1 1 1 1 1)} and

$C_4^2$ = {(0 0 0 0 0 0 0 0) (1 1 0 0 0 1 1 0) (1 0 1 0 1 0 1 1) (1 1 1 1 0 0 1 1) (1 1 0 0 1 1 1 1)}

be a set bicode in which the messages are carried only by $C_3^1$ and $C_2^2$ and all other codes are just to mislead the intruder.

Any message X sent across would be like X = {(1 1 0 1 0 0), (1 1 1 1 1), (0 0 1 1 1 0 1)} ∪ {(1 1 1 0 1 0 0), (1 1 0 1 1 0), (1 1 1 1 1), (1 1 0 0 1 1 1 1)}. Clearly both the sender and the receiver knows only the code words (0 0 1 1 1 0 1) and (1 1 0 1 1 0 0) carry the messages and all other symbols are just to mislead the intruder. So with the usual coding technique he would only work with (0 0 1 1 1 0 1) ∪ (1 1 0 1 1 0) and would ignore all other messages.

The advantages of these codes are

(1) When these bicodes are used it is impossible for the intruder to guess which codes carry the real message their by the security while transmission is preserved so best suited for cryptography.

(2) Further these codes can be used in defence department so that only one or two defence personals alone knows where the secret informations are stored hence no traitor can easily guess the secrets.
Even if invaded by the enemy nation it may not be easy for them to get the true information. Thus these set bicodes can be used as storage codes in defence departments.

(3) These codes are suited to the present situation, i.e., the computer world were mode of transmission is not very difficult. These codes can be used in communication in computers were secrecy of the identity is to be maintained.

It is pertinent to mention that even only one code in the set bicode be a usual code i.e., a group or a subspace of the vector space of dimension $2^n$ and all the other codes are just sets which are used only for misleading purposes.By making the number $n_1 + n_2$ of the set bicode $C = C_1 \cup C_2 = \left(C_1^1, C_2^1, \ldots, C_{n_1}^1\right) \cup \left(C_1^2, C_2^2, \ldots, C_{n_2}^2\right)$ arbitrarily large we see it is difficult for any intruder to break the code or hack the message; for only one code in these $n_1 + n_2$ codes will carry the message the rest of the $n_1 + n_2 - 1$ code words are just misleading one which may be even identical with the code word which carries the message.

We illustrate this by a simple example.

***Example 3.1.32:*** Let

$$V = V_1 \cup V_2 = \left(V_1^1, V_2^1, V_3^1, V_4^1, V_5^1\right) \cup \left(V_1^2, V_2^2, V_3^2, V_4^2\right)$$

where $V_3^1$ alone is a code which carries the messages all other code words in the set bicode $V = V_1 \cup V_2$ are only misleading code words or false code words. $V = V_1 \cup V_2 =$

$V_1^1$ = {(1 1 1 1 1 1) (0 0 0 0 0 0) (1 1 0 0 1 1) (0 0 1 1 0 0) (1 1 1 0 1 0)}

$V_2^1$ = {(0 0 0 0 0) (1 0 0 1 1 0) (0 1 1 1 0) (1 1 1 1 1)}

$V_3^1$ = {(0 0 0 0 0 0 0) (1 0 0 0 1 0 1) (0 1 0 0 1 1 1) (0 0 1 0 1 1 0) (0 0 0 1 0 1 1) (1 1 0 0 0 1 0) (0 1 1 0 0 0 1) (0 0 1 1 1 0 1) (1 0 0 1 1 1 0) (1 0 1 0 0 1 1) (0 1 0 1 1 0 0) (1 1 1 0 1 0 0) (0 1 1 1 0 1 0) (1 1 0 1 0 0 1) (1 0 1 1 0 0 0) (1 1 1 1 1 1 1)}

$V_4^1$ = {(1 0 0 0 1 1 0 1) (0 0 0 0 0 0 0 0) (1 1 0 0 0 1 1 0) (1 0 1 0 1 0 1 0 0)}

$V_5^1$ = {(0 0 0 0 0 0 0 0 0) (1 1 1 1 0 0 0 0 0) (0 0 0 1 1 1 0 0 0) (1 1 1 1 1 1 0 0 0) (0 0 1 1 1 1 1 0 0)}

$V_1^2$ = {(1 1 1 1 1 1 1) (0 0 0 0 0 0 0) (0 1 1 1 0 1 0) (1 1 1 1 1 1 0) (0 0 1 1 1 1 1) (1 0 1 1 0 0 0)}

$V_2^2$ = {(1 1 0 0 0) (0 0 0 0 0) (0 1 0 0 1) (0 1 1 0 0 0)}

$V_3^2$ = {(0 0 0 0 0 0) (1 1 1 0 1 0) (0 0 1 1 0 0) (1 1 0 0 1 1)}

$V_4^2$ = {(0 0 0 0 0 0 0 0) (1 1 1 1 1 1 1 1)(1 1 1 1 0 0 0 0) (1 1 0 0 1 1 1 1) (1 0 1 0 1 0 1 0) (0 1 0 1 0 1 0 1) (1 1 1 0 0 0 1 1)}

is a set bicode in which only the code $V_3^1$ is the (7, 4) code which is cyclic. Thus if X = $X_1 \cup X_2 \in V = V_1 \cup V_2$ = {(1 1 1 0 1 0) (0 1 1 1 0) (0 1 0 0 1 1 1) (1 0 0 0 1 1 0 1) (1 1 1 1 1 1 0 0 0)} ∪ {(0 1 1 1 0 1 0) (0 1 0 0 1) (1 1 1 0 1 0) (1 1 0 0 1 1 1 1)}.

The receiver knows the messages are only in the code (0 1 0 0 1 1 1) and its other code words got by cyclic shift of (0 1 0 0 1 1 1) i.e {(1 0 1 0 0 1 1), (1 1 0 1 0 0 1), (1 1 1 0 1 0 0), (0 1 1 1 0 1 0), (0 0 1 1 1 0 1) and (1 0 0 1 1 1 0)} thus only these 7 code words carry the messages. So the sender instead of sending these 7 code words sends the bicode word $X_1 \cup X_2$ and (0 1 0 0 1 1 1) carries the hint for the receiver to know the messages as the key is already known to him.

Such set bicodes when used in cryptography, it is near impossibility for the intruder in the first place to find the code words which carries the message from these $n_1 + n_2$ code words and secondly to guess it is the totally of all cyclic shift's of that word which carry the message. These codes can be used in defence department with complete confidence.

We give yet another unbreakable set bicodes.

***Example 3.1.33:*** Let

$$C = C_1 \cup C_2 = \left(C_1^1, C_2^1, C_3^1, C_4^1, C_5^1\right) \cup \left(C_1^2, C_2^2, C_3^2, C_4^2\right)$$

where

$C_1^1$ = {(0 0 0 0 0 0) (1 1 0 0 1 0) (1 1 1 1 1 1) (0 1 1 0 1 0) (1 0 1 0 1 0)}

$C_2^1$ = {(0 0 0 0 0 0 0) (1 0 1 0 1 0 1) (1 1 1 0 0 0 1) (1 0 1 1 0 1 1) (0 0 1 1 1 1 0)}

$C_3^1$ = {(0 0 0 0 0) (1 1 1 1 0) (1 1 1 0 0) (0 1 1 1 1) (0 1 0 1 0) (1 1 0 1 1)}

$C_4^1$ = {(1 1 1 1 0 0 0 0) (1 1 0 1 1 0 0 0) (0 0 0 0 0 0 0 0) (1 1 0 1 1 1 1 1)}

$C_5^1$ = {(1 1 1 1) (0 0 0 0) (0 1 0 1) (1 0 1 0)}

$C_1^2$ = {(1 1 1 0 0 0 0 1) (0 0 0 0 0 0 0 0) (0 1 1 1 0 0 0 1) (0 0 1 1 1 1 1 0) (1 1 0 0 0 1 1 0)}

$C_2^2$ = {(0 0 0 0 0 0 0) (1 1 0 1 0 0 0) (0 1 1 0 1 0 0) (0 0 1 1 0 1 0) (0 0 0 1 1 0 1) (1 0 1 1 1 0 0) (0 1 0 1 1 1 0) (0 0 1 0 1 1 1) (1 1 1 0 0 1 0) (0 1 1 1 0 0 1) (1 1 0 0 1 0 1) (1 0 0 0 1 1 0)(0 1 0 0 0 1 1) (1 0 1 0 0 0 1) (1 1 1 1 1 1 1) (1 0 0 1 0 1 1)}

$C_3^2$ = {(0 0 0 0 0 0) (1 1 0 1 1 0) (0 0 1 1 1 0) (1 0 1 0 1 0) (0 1 0 1 0 1)}

$C_4^2$ = {(1 1 1 1 1) (0 0 0 0 0) (1 1 0 1 0) (0 0 1 0 1)} .

The receiver knows that only the message is carried by the code $C_2^2$ in $C_2$. Thus if X = $X_1 \cup X_2$ = {(1 1 1 1 1 1) (1 0 1 0 1 0 1) (1 1 1 0) (1 1 0 1 1 1 1 1) (0 1 0 1)} ∪ {(1 1 1 0 0 0 0 1) (0 0 1 1 0 1 0) (1 0 1 0 1 0) (0 0 1 0 1)} is sent as a message the receiver knows only those code words orthogonal to the code word (0 0 1 1 0 1 0) carries the information and all other code words are just to mislead everyone, thus he works with all code words y of length 7 which are orthogonal to (0 0 1 1 0 1 0); i.e (0 0 1 1 0 1 0)y = 0. Clearly y = (0 0 0 0 0 0 0) but this is ignored. y = {(1 1 0 0 1 0 1), (1 1 0 0 0 0 0), (1 1 1 1 0 0 0), (1 0 0 0 0 0 1), (1 0 0 0 1 0 1), (0 0 1 1 1 0 1) and (0 0 1 0 0 1 0)} and so on will carry the message.

To this end we define the notion of codes orthogonal to a given code words.

**DEFINITION 3.1.20:** *Let C be a (n, k) code. Suppose x = ($x_1$, $x_2$, …, $x_n$) where $x_i \in \{0, 1\}$; $1 \le i \le n$. Now the codes orthogonal to x is given by $x^\perp$ = {y ∈ {S x S x S x S x S x S x S} where S = {0, 1} such that x.y = (0)}. $x^\perp$ is a set which is always non empty for (0 0 0 … 0) ∈ $x^\perp$.*

Now we shall illustrate by a simple example.

***Example 3.1.34:*** Let C = {(0 1 1 1 1) (1 0 0 0 0) (0 1 0 1 0) (0 0 0 0 0)} be a set code. Take x = (0 1 0 1 0) ∈ C. $x^\perp$ = {(0 0 0 0 0) (0 1 0 1 0) (1 0 1 0 0) (1 0 1 0 1) (0 0 1 0 1) (1 0 0 0 0) (0 0 1

0 0) (0 0 0 0 1) (1 1 1 1 0) (1 0 0 0 1) (1 1 1 1 1) (0 1 1 1 1) (1 1 0 1 1) (0 1 1 1 0) (1 1 0 1 0) (0 1 0 1 1)}

Certainly the elements of $x^{\perp}$ are not orthogonal with every element of given set code C.

Now we will indicate how the working can carried out by a cryptologist in cryptology. When a set bicode

$$X = X_1 \cup X_2 = \left(X_1^1, X_2^1, \ldots, X_{n_1}^1\right) \cup \left(X_1^2, \ldots, X_{n_2}^2\right)$$

is transmitted the sender knows out of these $n_1 + n_2$ code words the particular $X_{t_i}^i$, $1 \leq t_i \leq n_i$, i = 1, 2 which carries the message. That message may not be direct, once the receiver spots that code he has to supply a password or confidential code to the sender only on sending the correct one the message would be sent as the message is in one or more of the orthogonal codes to X which he has to get correctly by identifying the weight or by the finding the syndrome of it for already it is said the sender and receiver will make one of the $X_{t_j}^i$ to be a normal code with the supplied parity check matrix or the orthogonal parity check matrix. So the receiver decodes the message using the parity check matrix. This sort of set bicodes when used it is not easy for the intruder to hack the information easily.

We can prove the following.

**THEOREM 3.1.5:** *All set repetition bicodes are group bicodes.*

Proof is left as an exercise for the reader.

However we give an example of it.

***Example 3.1.35:*** Let

$$V = V_1 \cup V_2 = \left(V_1^1, V_2^1, V_3^1\right) \cup \left(V_1^2, V_2^2, V_3^2, V_4^2\right)$$

where $V_1^1$ = {(0 0 0 0 0) (1 1 1 1 1)}, $V_2^1$ = {(1 1 1 1 1 1 1) (0 0 0 0 0 0 0)}, $V_3^1$ ={(0 0 0 0 0 0 0 0) (1 1 1 1 1 1 1 1)}, $V_1^2$ = {(1 1 1 1 1 1) (0 0 0 0 0 0)}, $V_2^2$ = {(1 1 1 1) (0 0 0 0)}, $V_3^2$ = {(1 1 1 1 1 1 1 1 1) (0 0 0 0 0 0 0 0 0) and $V_4^2$ = {(1 1 1 1 1 1 1) (0 0 0 0 0 0)} is a repetition set bicode which is also a group bicode as every code in $V_i$, i = 1, 2 is a repetition code.

However a set bicode in general is not a group bicode. Now we proceed on to define the notion of group parity check bicode and group weak parity check bicode.

**DEFINITION 3.1.21:** *Let*

$$C = C_1 \cup C_2 = \left(C_1^1, C_2^1, \ldots, C_{n_1}^1\right) \cup \left(C_1^2, \ldots, C_{n_2}^2\right)$$

*be a group bicode if each of the codes* $C_{t_i}^i$ *,* $1 \le t_i \le n_i$*, i = 1, 2 is such that it is associated with a parity check matrix* $H_{j_i}^i$ *of the form (1 1 … 1),* $1 \le j_i \le n_i$*, i = 1, 2 then we call C to be a group parity check bicode. If the group bicode*

$$C = C_1 \cup C_2 = \left(C_1^1, C_2^1, \ldots, C_{n_1}^1\right) \cup \left(C_1^2, \ldots, C_{n_2}^2\right)$$

*where some of the* $C_{j_1}^1$ *and* $C_{j_2}^2$ *are parity check code then we call C to be a group weak parity check bicode,* $1 \le j_1 \le n_1$ *and* $1 \le j_1 \le n_2$*.That is not all the codes* $C_{j_i}^i$ *are parity check codes only a few of them are parity check codes.*

We will illustrate these definitions by some examples.

***Example 3.1.36:*** Let

$$\mathrm{V} = \mathrm{V}_1 \cup \mathrm{V}_2 = \left(\mathrm{V}_1^1, \mathrm{V}_2^1, \ldots, \mathrm{V}_{n_1}^1\right) \cup \left(\mathrm{V}_1^2, \ldots, \mathrm{V}_{n_2}^2\right)$$

where $n_1 = 4$ and $n_2 = 3$ with

$V_1^1$ = {(0 0 0 0 0 0) (1 1 0 0 0 0) (0 0 1 1 1 1) (1 1 1 1 1 1)}

$V_2^1$ = {(0 0 0 0 0) (1 1 0 0 0) (0 0 1 1 0) (1 0 1 0 0) (1 1 1 1 0) (0 0 0 1 1) (1 0 0 0 1) (0 1 1 0 0)}

$V_3^1$ = {(0 0 0 0 0 0 0) (1 1 1 1 0 0 0) (0 0 0 1 1 1 1) (1 1 0 0 0 1 1) (1 1 0 1 1 0) (0 1 1 1 1 0 0) (0 0 0 1 1 0 0) (1 1 1 1 1 0) (0 1 1 1 1 1 1) (1 1 1 0 0 1 0)}

$V_4^1$ = {(1 1 1 1) (0 0 0 0) (1 1 0 0) (0 0 1 1) (1 0 0 1) (0 1 1 0) (1 0 1 0)}

$V_2^1$ = {(1 1 1 1 1 1) (0 1 0 1 0 0) (1 0 1 0 1 1) (1 1 1 1 0 0) (0 0 1 1 0 0) (0 0 1 1 1 1)}

$V_2^2$ = {(1 1 1 1 1 1 1 1) (1 1 0 0 1 1 0 0) (0 0 1 1 0 0 1 1) (1 1 1 1 0 0 0 0) (0 0 0 0 1 1 1 1) (1 1 1 1 1 1 0 0) (0 0 1 1 1 1 1 1)}

$V_3^2$ = {(0 0 0 0 0) (1 1 1 1 0) (0 1 1 1 1) (0 1 0 1 0) (1 0 1 0 0) (1 1 0 1 1) (1 1 1 0 1)}

with the group parity check matrices H = $H_1 \cup H_2$ = { $H_1^1$ = {(1 1 1 1 1 1), $H_2^1$ = {(1 1 1 1 1), $H_3^1$ = (1 1 1 1 1 1 1), $H_4^1$ = (1 1 1 1)} $\cup H_1^2$ = {(1 1 1 1 1 1), $H_2^2$ = {(1 1 1 1 1 1 1 1), $H_3^2$ = (1 1 1 1 1)}. Clearly V is a group parity check bicode.

Now we proceed on to give examples of weak group parity check bicode.

***Example 3.1.37:*** Let

$$C = C_1 \cup C_2 = \left\{C_1^1, C_2^1, C_3^1, C_4^1\right\} \cup \left\{C_1^2, C_2^2, C_3^2, C_4^2, C_5^2\right\}$$

where $C_1^1$ is generated by the generator matrix

$$G_1^1 = \begin{pmatrix} 0 & 1 & 1 & 1 & 1 \\ 1 & 0 & 0 & 1 & 0 \end{pmatrix},$$

$C_2^2$ is generated by the matrix

$$G_2^1 = \begin{pmatrix} 1 & 0 & 0 & 0 & 1 & 1 \\ 0 & 1 & 0 & 1 & 0 & 1 \\ 0 & 0 & 1 & 1 & 1 & 0 \end{pmatrix}.$$

$C_3^1$ associated with the parity check matrix

$$H_3^1 = (1\ 1\ 1\ 1\ 1\ 1\ 1),$$

$C_4^1$ associated with the parity check matrix

$$H_4^1 = (1\ 1\ 1\ 1\ 1\ 1\ 1\ 1),$$

$C_1^2$ generated by

$$G_1^2 = \begin{pmatrix} 1 & 1 & 1 & 0 & 0 \\ 0 & 0 & 1 & 1 & 0 \\ 1 & 1 & 1 & 1 & 1 \end{pmatrix}$$

$$G_2^2 = \begin{pmatrix} 1 & 0 & 0 & 1 & 0 & 0 \\ 0 & 1 & 0 & 0 & 1 & 0 \\ 0 & 0 & 1 & 0 & 0 & 1 \end{pmatrix}$$

generates the code $C_1^2$ and $C_2^2$, $C_3^2$ has the associated parity check matrix

$$H_4^2 = (1\ 1\ 1\ 1\ 1),$$

$C_4^2$ is generated by

$$G_4^2 = \begin{pmatrix} 1 & 1 & 0 & 1 & 0 & 0 & 0 \\ 0 & 1 & 1 & 0 & 1 & 0 & 0 \\ 0 & 0 & 1 & 1 & 0 & 1 & 0 \\ 0 & 0 & 0 & 1 & 1 & 0 & 1 \end{pmatrix},$$

and $C_5^2$ has the associated parity check matrix

$$H_5^2 = (1\ 1\ 1\ 1\ 1\ 1\ 1\ 1).$$

Thus $C = C_1 \cup C_2$ is a weak parity check group bicode.

Now we proceed on to define the notion of Hamming group bicode.

**DEFINITION 3.1.22:** *Let*

$$C = C_1 \cup C_2 = \left(C_1^1, C_2^1, \ldots, C_{n_1}^1\right) \cup \left(C_1^2, C_2^2, \ldots, C_{n_2}^2\right)$$

*be a group bicode where each of the codes $C_{j_i}^i$ is a Hamming code for $i \le j_i \le n_i$; i = 1, 2. We say $C = C_1 \cup C_2$ to be a group Hamming bicode. If in the group bicode $C = C_1 \cup C_2 = \left(C_1^1, C_2^1, \cdots, C_{n_1}^1\right) \cup \left(C_1^2, C_2^2, \cdots, C_{n_2}^2\right)$ atleast some of the $C_{j_i}^i$ 's are Hamming codes then we call C to be a weak Hamming group bicode or weak group Hamming bicode.*

We will illustrate these by simple examples.

***Example 3.1.38:*** Let

$$C = C_1 \cup C_2 = \left\{C_1^1, C_2^1, C_3^1, C_4^1, C_5^1\right\} \cup \left\{C_1^2, C_2^2, C_3^2, C_4^2\right\}$$

where

$C_1^1$ = {(0 0 0 0 0 0) (1 1 1 0 0 0) (0 0 0 1 1 1) (1 1 1 1 1 1)},

$C_2^1$ is generated by the matrix

$$\begin{pmatrix} 1 & 0 & 0 & 0 & 1 & 0 & 1 \\ 0 & 1 & 0 & 0 & 1 & 1 & 1 \\ 0 & 0 & 1 & 0 & 1 & 1 & 0 \\ 0 & 0 & 0 & 1 & 0 & 1 & 1 \end{pmatrix} = G_2^1,$$

$C_3^1$ is associated with the parity check matrix

$$H_3^1 = \begin{pmatrix} 1 & 1 & 0 & 0 & 0 & 0 & 0 \\ 1 & 0 & 1 & 0 & 0 & 0 & 0 \\ 1 & 0 & 0 & 1 & 0 & 0 & 0 \\ 1 & 0 & 0 & 0 & 1 & 0 & 0 \\ 1 & 0 & 0 & 0 & 0 & 1 & 0 \\ 1 & 0 & 0 & 0 & 0 & 0 & 1 \end{pmatrix};$$

the generator matrix associated with $C_4^1$ is

$$G_4^1 = \begin{pmatrix} 1 & 1 & 1 & 0 & 0 \\ 0 & 0 & 1 & 1 & 0 \\ 1 & 1 & 1 & 1 & 1 \end{pmatrix}$$

and

$$H_5^1 = \begin{pmatrix} 1 & 1 & 0 & 0 & 0 & 0 & 0 & 0 \\ 1 & 0 & 1 & 0 & 0 & 0 & 0 & 0 \\ 1 & 0 & 0 & 1 & 0 & 0 & 0 & 0 \\ 1 & 0 & 0 & 0 & 1 & 0 & 0 & 0 \\ 1 & 0 & 0 & 0 & 0 & 1 & 0 & 0 \\ 1 & 0 & 0 & 0 & 0 & 0 & 1 & 0 \\ 1 & 0 & 0 & 0 & 0 & 0 & 0 & 1 \end{pmatrix}$$

is the parity check matrix associated with the code $C_5^1$.

$$C_1^2 = \{(1\ 1\ 1\ 1\ 1\ 1\ 1\ 1)\ (0\ 0\ 0\ 0\ 0\ 0\ 0\ 0)\}.$$

$C_2^2$ generated by the generator matrix

$$G_2^2 = \begin{pmatrix} 1 & 0 & 0 & 1 & 0 & 0 \\ 0 & 1 & 0 & 0 & 1 & 0 \\ 0 & 0 & 1 & 0 & 0 & 1 \end{pmatrix}.$$

$C_3^2$ is associated with the parity check matrix

$$H_3^2 = \begin{pmatrix} 1 & 1 & 0 & 0 & 0 & 0 & 0 \\ 1 & 0 & 1 & 0 & 0 & 0 & 0 \\ 1 & 0 & 0 & 1 & 0 & 0 & 0 \\ 1 & 0 & 0 & 0 & 1 & 0 & 0 \\ 1 & 0 & 0 & 0 & 0 & 1 & 0 \\ 1 & 0 & 0 & 0 & 0 & 0 & 1 \end{pmatrix}$$

and $C_4^2$ is associated with the parity check matrix

$$
H_4^2 = \begin{pmatrix} 1 & 1 & 0 & 0 & 0 & 0 & 0 & 0 \\ 1 & 0 & 1 & 0 & 0 & 0 & 0 & 0 \\ 1 & 0 & 0 & 1 & 0 & 0 & 0 & 0 \\ 1 & 0 & 0 & 0 & 1 & 0 & 0 & 0 \\ 1 & 0 & 0 & 0 & 0 & 1 & 0 & 0 \\ 1 & 0 & 0 & 0 & 0 & 0 & 1 & 0 \\ 1 & 0 & 0 & 0 & 0 & 0 & 0 & 1 \end{pmatrix}.
$$

Clearly this a group bicode which is a weak group Hamming bicode (weak group Hamming bicode).

**DEFINITION 3.1.23:** *Let*

$$
C = C_1 \cup C_2 = \left(C_1^1, C_2^1, \ldots, C_{n_1}^1\right) \cup \left(C_1^2, C_2^2, \ldots, C_{n_2}^2\right)
$$

*be a group bicode, if each code* $C_{j_t}^i$, $1 \leq j_t \leq n_i$ ; $i = 1, 2$ *are cyclic codes then we call C to be a group cyclic bicode.*

*In the group bicode atleast a few of them are cyclic codes and the rest codes which form a group under addition then we call C to be a weak group cyclic bicode.*

We shall illustrate these definitions by some examples.

***Example 3.1.39:*** Let

$$
C = C_1 \cup C_2 = \left\{C_1^1, C_2^1, C_3^1\right\} \cup \left\{C_1^2, C_2^2, C_3^2, C_4^2\right\}
$$

where $C_1^1$ = {(0 0 0 0 0) (1 1 0 0 0), (0 0 1 1 1), (1 1 1 1 1)}, $C_2^1$ = {(0 0 0 0 0 0) (1 1 0 1 0 1) (1 1 1 0 1 0) (0 1 1 1 0 1) (1 0 1 1 1 0) (0 1 0 1 1 1) (1 0 1 0 1 1) (1 1 0 1 0 1)}. Clearly $C_2^1$ is a cyclic code, $C_3^1$ = {(0 0 0 0 0 0 0) (1 1 0 0 0 0 0) (0 0 1 1 0 0 0) (0 0 0 0 1 1 1), (1 1 0 0 1 1 1) (0 0 1 1 1 1 1)} , $C_1^2$ = {(0 0 0 0 0 0) (1 1 0 0 0 0) (1 1 0 0 1 1) (0 0 0 0 1 1) (0 0 1 1 0 0) (1 1 1 1 0 0) (0 0 1 1 1 1) (1 1 1 1 1 1) and $C_2^2$ = {(1 1 0 1 0 1 1), (0 0 0 0 0 0 0) (1 1 1 0 1 0 1) (1 1 1 1 0 1 0) (0 1 1 1 1 0 1), (1 0 1 1 1 1 0) (0 1 0 1 1 1 1) (1 0 1 0 1 1 1)} is a cyclic code.

Let $C_3^2$ be generated by the generator matrix

$$G_3^2 = \begin{pmatrix} 1 & 0 & 0 & 0 & 1 & 0 & 0 & 0 \\ 1 & 0 & 0 & 0 & 0 & 1 & 0 & 0 \\ 1 & 0 & 0 & 0 & 0 & 0 & 1 & 0 \\ 1 & 0 & 0 & 0 & 0 & 0 & 0 & 1 \end{pmatrix}.$$

$C_4^2$ be the cyclic code generated by the generator matrix

$$G_4^2 = \begin{pmatrix} 1 & 1 & 1 & 0 & 1 & 0 & 0 \\ 0 & 1 & 1 & 1 & 0 & 1 & 0 \\ 0 & 0 & 1 & 1 & 1 & 0 & 1 \end{pmatrix}.$$

Clearly $C = C_1 \cup C_2$ is a weak cyclic group bicode.

Now we proceed onto give an example of a group cyclic bicode.

***Example 3.1.40:*** Let

$$C = C_1 \cup C_2 = \left\{C_1^1, C_2^1, C_3^1, C_4^1\right\} \cup \left\{C_1^2, C_2^2, C_3^2\right\}$$

where each $C_{t_j}^i$ is a cyclic code; $1 \leq t_j \leq 4$ or 3; i = 1, 2. $C_1^1$ is generated by the generator matrix.

$$G_1^1 = \begin{pmatrix} 1 & 0 & 0 & 1 & 0 & 0 \\ 0 & 1 & 0 & 0 & 1 & 0 \\ 0 & 0 & 1 & 0 & 0 & 1 \end{pmatrix},$$

$$G_2^1 = \begin{pmatrix} 1 & 1 & 1 & 0 & 1 & 0 & 0 \\ 0 & 1 & 1 & 1 & 0 & 1 & 0 \\ 0 & 0 & 1 & 1 & 1 & 0 & 1 \end{pmatrix}$$

is the generator matrix of the code $C_2^1$.

$$G_3^1 = \begin{pmatrix} 1 & 0 & 1 & 1 & 1 \\ 0 & 1 & 0 & 1 & 0 \end{pmatrix}$$

is the generator matrix of the $C_3^1$ code. $C_4^1$ = {(1 1 1 0 0 0 1 1 0) (0 1 1 1 0 0 0 1 1) (1 0 1 1 1 0 0 0 1) (1 1 0 1 1 1 0 0 0) (0 1 1 0 1 1 1 0 0) (0 0 1 1 0 1 1 1 0) (0 0 0 1 1 0 1 1 1) (1 0 0 0 1 1 0 1 1) (1 1 0 0 0 1 1 0 1) (0 0 0 0 0 0 0 0 0)}.

$$\begin{pmatrix} 1 & 0 & 0 & 0 & 1 & 1 & 0 \\ 0 & 1 & 0 & 0 & 0 & 1 & 1 \\ 0 & 0 & 1 & 0 & 0 & 0 & 1 \\ 0 & 0 & 0 & 1 & 1 & 1 & 1 \end{pmatrix} = H_1^2$$

is the parity check matrix of the code $C_1^2$. $C_2^2$ is a code associated with the parity check matrix

$$H_2^2 = \begin{pmatrix} 1 & 1 & 0 & 0 & 0 & 0 & 0 & 0 & 0 \\ 1 & 0 & 1 & 0 & 0 & 0 & 0 & 0 & 0 \\ 1 & 0 & 0 & 1 & 0 & 0 & 0 & 0 & 0 \\ 1 & 0 & 0 & 0 & 1 & 0 & 0 & 0 & 0 \\ 1 & 0 & 0 & 0 & 0 & 1 & 0 & 0 & 0 \\ 1 & 0 & 0 & 0 & 0 & 0 & 1 & 0 & 0 \\ 1 & 0 & 0 & 0 & 0 & 0 & 0 & 1 & 0 \\ 1 & 0 & 0 & 0 & 0 & 0 & 0 & 0 & 1 \end{pmatrix}.$$

$C_3^2$ is the cyclic group given by {(0 0 0 0 0 0 0 0) (1 0 1 1 1 0 1) (1 1 1 0 1 1 1 0) (0 1 1 1 0 1 1 1) (1 0 1 1 1 0 1 1) (1 1 0 1 1 1 0 1) (1 1 1 0 1 1 1 0)}. $C_4^2$ is the code generated by the matrix

$$G_4^2 = \begin{pmatrix} 1 & 0 & 0 & 1 & 0 & 0 \\ 0 & 1 & 0 & 0 & 1 & 0 \\ 0 & 0 & 1 & 0 & 0 & 1 \end{pmatrix}.$$

Clearly the group bicode is a weak group cyclic bicode.

***Example 3.1.41:*** Let

$$C = C_1 \cup C_2 = \{C_1^1, C_2^1, C_3^1, C_4^1\} \cup \{C_1^2, C_2^2, C_3^2\}$$

be a group bicode where $C_1^1$ is the cyclic code generated by the generator matrix

$$G_1^1 = \begin{pmatrix} 1 & 1 & 0 & 1 & 0 & 0 & 0 \\ 0 & 1 & 1 & 0 & 1 & 0 & 0 \\ 0 & 0 & 1 & 1 & 0 & 1 & 0 \\ 0 & 0 & 0 & 1 & 1 & 0 & 1 \end{pmatrix}.$$

$C_2^1$ is the cyclic code generated by the generator matrix

$$G_2^1 = \begin{pmatrix} 1 & 0 & 0 & 1 & 0 & 0 \\ 0 & 1 & 0 & 0 & 1 & 0 \\ 0 & 0 & 1 & 0 & 0 & 1 \end{pmatrix}.$$

$C_3^1$ is the cyclic code generated by the generator matrix

$$G_3^1 = \begin{pmatrix} 1 & 1 & 1 & 0 & 1 & 0 & 0 \\ 0 & 1 & 1 & 1 & 0 & 1 & 0 \\ 0 & 0 & 1 & 1 & 1 & 0 & 1 \end{pmatrix}$$

and $C_4^1$ is the cyclic code generated by the generator matrix

$$G_4^1 = \begin{pmatrix} 1 & 0 & 0 & 0 & 1 & 0 & 1 \\ 0 & 1 & 0 & 0 & 1 & 1 & 1 \\ 0 & 0 & 1 & 0 & 1 & 1 & 0 \\ 0 & 0 & 0 & 1 & 0 & 1 & 1 \end{pmatrix}.$$

Now $C_1^2$ = {(1 1 0 1 0 1) (1 1 1 0 1 0) (0 1 1 1 0 1) (1 0 1 1 1 0) (0 1 0 1 1 1) (1 0 1 0 1 1) (0 0 0 0 0 0)} is a cyclic code. $C_2^2$ be a cyclic code generated by the matrix

$$G_2^2 = \begin{pmatrix} 1 & 0 & 0 & 0 & 1 & 0 & 1 \\ 0 & 1 & 0 & 0 & 1 & 1 & 1 \\ 0 & 0 & 1 & 0 & 1 & 1 & 0 \\ 0 & 0 & 0 & 1 & 0 & 1 & 1 \end{pmatrix}.$$

$C_3^2$ is a cyclic code generated by

$$G_3^2 = \begin{pmatrix} 1 & 0 & 0 & 1 & 0 & 0 \\ 0 & 1 & 0 & 0 & 1 & 0 \\ 0 & 0 & 1 & 0 & 0 & 1 \end{pmatrix}.$$

Thus C = $C_1 \cup C_2$ is a group cyclic bicode.

Now we proceed on to define the orthogonal code of a group code.

**DEFINITION 3.1.24:** *Let*

$$C = C_1 \cup C_2 = \left(C_1^1, C_2^1, \ldots, C_{n_1}^1\right) \cup \left(C_1^2, C_2^2, \ldots, C_{n_2}^2\right)$$

*be a group bicode, the dual group bicode of C denoted by*

$$C^\perp = \left(C_1^\perp \cup C_2^\perp\right) =$$

$$\left(\left(C_1^1\right)^\perp, \left(C_2^1\right)^\perp, \ldots, \left(C_{n_1}^1\right)^\perp\right) \cup \left(\left(C_1^2\right)^\perp, \left(C_2^2\right)^\perp, \ldots, \left(C_{n_2}^2\right)^\perp\right);$$

*i.e., for each $C_{j_t}^i$ in C, $C^\perp$ contains $\left(C_{j_t}^i\right)^\perp$ for $1 \le j_t \le n_i$; i = 1,2.*
*If for $C = C_1 \cup C_2$ we have*

$$C' = \left(\left(C_1^1\right)^\perp, C_2^1, \ldots, C_{n_1}^1\right) \cup \left(C_1^2, \left(C_2^2\right)^\perp, \ldots, \left(C_{n_2}^2\right)^\perp\right)$$

*with only some of the* $\left(C^{i}_{j_t}\right)^{\perp}$ *are in C' and other codes are just* $C^{i}_{j_s}$ *; 1 ≤ $j_t$, $j_s$ ≤ $n_i$; i = 1, 2; then we call C' to be the weak dual group code of C.*

We will illustrate this by some examples.

***Example 3.1.42:*** Let

$$C = C_1 \cup C_2 = \left\{C^1_1, C^1_2, C^1_3\right\} \cup \left\{C^2_1, C^2_2, C^2_3, C^2_4\right\}$$

where $C^1_1$ = {(0 0 0 0 0 0) (1 1 0 0 1 1) (0 0 1 1 0 0) (1 1 1 1 1 1)}, $C^1_2$ = {(0 0 0 0 0 0 0) (1 0 1 0 1 0 1) (0 1 0 1 0 1 0) (1 1 1 1 1 1 1)}, $C^1_3$ = {(0 0 0 0 0 0) (0 0 1 0 0 1) (0 1 0 0 1 0) (0 1 1 0 1 1), (1 0 0 1 0 0) (1 0 1 1 0 1) (1 1 0 1 1 0), (1 1 1 1 1 1)}, $C^2_1$ = {(0 0 0 0) (1 0 1 1) (1 1 1 0) (0 1 0 1)}, $C^2_2$ = {(0 0 0 0 0 0 0) (1 1 1 0 1 0 0) (0 1 1 1 0 1 0) (0 0 1 1 1 0 1) (1 0 0 1 1 1 0) (1 1 0 1 0 0 1) (0 1 0 0 1 1 1) (1 0 1 0 0 1 1)}, $C^2_3$ = {(0 0 0 0 0 0) (1 0 0 1 0 0) (0 1 1 0 1 1) (0 1 0 0 1 0) (1 0 1 1 0 1) (1 1 0 1 1 0) (1 1 1 1 1 1)} and $C^2_4$ = {(0 0 0 0 0 0) (1 0 0 0 1 0 1) (0 1 0 0 1 1 1) (0 0 1 0 1 1 0) (0 0 0 1 0 1 1) (1 1 0 0 0 1 0) (0 1 1 0 0 0 1) (0 0 1 1 1 0 1) (1 0 1 0 0 1 1) (0 1 0 1 1 0 0) (1 0 0 1 1 1 0) (1 1 1 0 1 0 0) (0 1 1 1 0 1 0) (1 1 0 1 0 0 0) (1 0 1 1 0 0 0) (1 1 1 1 1 1 1)} is the group bicode.

$C^{\perp} = C_1^{\perp} \cup C_2^{\perp}$ = { $\left(C^1_1\right)^{\perp}$ = (0 0 0 0 0 0) (1 1 0 0 1 1) (0 0 1 1 0 0) (1 1 1 1 1 1) (1 1 1 1 0 0 0 (0 0 0 0 1 1 0, (0 0 1 1 1 1)}, $\left(C^1_2\right)^{\perp}$ = {(0 0 0 0 0 0 0) (1 0 1 0 1 0 1), (1 0 1 0 0 0 0) (0 0 0 1 0 1)}, $\left(C^1_3\right)^{\perp}$ = {(0 0 0 0 0 0) (0 0 1 0 0 1) (1 0 0 1 0 0) (1 1 0 1 1 0) (0 1 0 0 1 0) (0 1 1 0 1 1)}} ∪ { $\left(C^2_1\right)^{\perp}$ = {(0 0 0 0) (1 0 1 0) (0 1 1 1) (1 1 0 1)}, $\left(C^2_2\right)^{\perp}$ = {(0 0 0 0 0 0 0) (1 1 1 0 1 0 0)}, $\left(C^2_3\right)^{\perp}$ = {(0 0 0 0 0 0) (1 0 0 1 0 0) (0 1 1 0 1 1) (1 1 1 1 1 1)} $\left(C^2_4\right)^{\perp}$ = {(0 0 0 0 0 0 0)}} is the dual group bicode of C.

Now we proceed on to give an example of a weak dual group bicode.

***Example 3.1.43:*** Let C = $C_1 \cup C_2$ = $\left(C_1^1, C_2^1, C_3^1\right) \cup \left(C_1^2, C_2^2\right)$ where $C_1^1$ = {(0 0 0 0 0 0), (0 0 1 0 0 1) (0 1 0 0 1 0) (0 1 1 0 1 1) (1 0 0 1 0 0) (1 0 1 1 0 1) (1 1 0 1 1 0) (1 1 1 1 1 1)}, $C_2^1$ = {(0 0 0 0 0 0 0) (1 1 1 0 1 0 0) (0 1 1 1 0 1 0) (0 0 1 1 1 0 1) (1 0 0 1 1 1 0) (0 1 0 0 1 1 1) (1 1 0 1 0 0 1) (1 0 1 0 0 1 1)}, $C_3^1$ = {(0 0 0 0 0) (0 1 1 1 1) (1 0 0 1 0) (1 1 1 0 1)}, $C_1^2$ = {(1 1 1 0 0) (0 0 1 1 0) (1 1 1 1 1) (1 1 0 1 0) (0 0 1 0 1) (1 1 0 0 1) (0 0 0 1 1) (0 0 0 0 0)} and $C_2^2$ = {(1 1 1 0) (0 1 1 0) (0 0 1 1) (1 0 0 0) (0 1 0 1) (1 1 0 1) (1 0 1 1) (0 0 0 0)}. Now the weak dual group bicode of C denoted by

$$C' = \left(\left(C_1^1\right)^{\perp}, C_2^1, \left(C_3^1\right)^{\perp}\right) \cup \left(\left(C_1^2\right)^{\perp}, C_2^2\right)$$

= { $\left(C_1^1\right)^{\perp}$ = {(0 0 0 0 0 0) (1 1 1 1 1 1) (1 0 0 1 0 0) (0 1 0 0 1 0) (0 0 1 0 0 1) (0 1 1 0 1 1) (1 0 1 1 0 1) (1 1 0 1 1 0)}, $C_2^1$, $\left(C_3^1\right)^{\perp}$ = {(0 0 0 0 0), (0 1 0 0 1), (0 1 1 0 0), (0 0 1 0 1)}} ∪ { $\left(C_1^2\right)^{\perp}$ = {(1 1 0 0 0) (0 0 0 0 0)}, $C_2^2$ }. C' is a weak dual group bicode. One can have several dual group bicodes for a given group bicode.

Now we proceed on to define the new notion of whole group bicode.

**DEFINITION 3.1.25:** *Let*

$$C = C_1 \cup C_1^{\perp} = \left(C_1^1, C_2^1, \ldots, C_n^1\right) \cup \left(\left(C_1^1\right)^{\perp}, \left(C_2^1\right)^{\perp}, \ldots, \left(C_n^1\right)^{\perp}\right)$$

*be a group bicode such that not all* $\left(C_i^1\right)^{\perp}$ = $\left(C_i^1\right)$ *for i = 1, 2, …, n. Then we call C to be a whole group bicode.*

We illustrate this by some examples.

***Example 3.1.44:*** Let

$$\left(C_1^1, C_2^1, C_3^1, C_4^1\right) \cup \left(\left(C_1^1\right)^{\perp}, \left(C_2^1\right)^{\perp}, \left(C_3^1\right)^{\perp}, \left(C_4^1\right)^{\perp}\right)$$

where $C_1^1$ = {(0 0 0 0) (1 0 1 1) (0 1 0 1) (1 1 1 0)}, $\left(C_1^1\right)^{\perp}$ = {(0 0 0 0) (1 0 1 0) (1 1 0 1) (0 1 1 1)}, $C_2^1$ = {(0 0 0 0 0 0) (0 0 1 0 0 1) (0 1 0 0 1 0 0 (0 1 1 0 1 1) (1 0 0 1 0 0) (1 0 1 1 0 1) (1 1 0 1 1 0) (1 1 1 1 1 1)}, $\left(C_2^1\right)^{\perp}$ = {(0 0 0 0 0 0) (0 0 1 0 0 1) (0 1 0 0 1 0) (0 1 1 0 1 1) (1 0 0 1 0 0) (1 0 1 1 0 1) (1 1 0 1 1 0) (1 1 1 1 1 1)}, $C_3^1$ = {(0 0 0 0 0) (0 1 1 1 1) (1 0 0 1 0) (1 1 1 0 1)}, $\left(C_3^1\right)^{\perp}$ = {(0 0 0 0 0) (0 1 1 0 0 0) (0 1 0 0 1) (0 0 1 0 1 0 (1 0 1 1 0) (1 0 0 1 1) (1 1 1 1 1)}, $C_4^1$ = {(0 0 0 0 0 0 0) (0 0 0 1 1 1 1) (0 1 1 0 0 1 1) (1 0 1 0 1 0 1) (0 1 1 1 1 0 0) (1 0 1 1 0 1 0) (1 1 0 0 1 1 0) (1 1 0 1 0 0 1)} and $\left(C_4^1\right)^{\perp}$ = {(0 0 0 0 0 0 0) (0 0 0 1 1 1 1) (0 1 1 0 0 1 1) (1 0 1 0 1 0 1) (0 1 1 1 1 0 0) (1 0 1 1 0 1 0) (1 1 0 0 1 1 0) (1 1 0 1 0 0 1 0 (1 1 1 1 1 1 1 0, (1 1 1 0 0 0 0), (1 0 0 0 0 1 1) (0 1 0 0 1 0 1) (1 0 0 1 1 0 0) (0 1 0 1 0 1 0) (0 0 1 1 0 0 1) (0 0 1 0 1 1 0)}. Thus $C = C_1 \cup C_1^{\perp}$ is a whole group bicode.

Now we proceed on to define the new notion of weighted group bicode.

**DEFINITION 3.1.26:** *Let*

$$C = C_1 \cup C_2 = \left(C_1^1, C_2^1, \ldots, C_{n_1}^1\right) \cup \left(C_1^2, C_2^2, \ldots, C_{n_2}^2\right)$$

*be a group bicode such that each $C_{j_t}^i$ is of weight m, $1 \le j_t \le n_i$, i = 1, 2; then we call the group bicode $C = C_1 \cup C_2$ to be (m, m) weighted group bicode.*

***Example 3.1.45:*** Let $C = C_1 \cup C_2 = \left(C_1^1, C_2^1, C_3^1\right) \cup \left(C_1^2, C_2^2, C_3^2\right)$ where $C_1^1$ = {(0 0 0 0), (1 1 1 1)}, $C_2^1$ = {(0 0 0 0 0 0) (1 1 1 1 0 0) (1 1 0 0 1 1) (0 0 1 1 1 1)}, $C_3^1$ = {(0 0 0 0 0 0 0 0) (1 1 0 0 0

0 0 1 1) (1 1 1 1 0 0 0 0) (0 0 1 1 0 0 1 1)}, $C_1^2$ = {(0 0 0 0 0) (1 1 1 1 0)}, $C_2^2$ = {(0 1 1 1 1 0) (0 0 0 0 0 0) (1 0 1 1 0 1) (1 1 0 0 1 1)} and $C_3^2$ = {(1 0 1 0 1 0 1) (0 0 0 0 0 0 0)} is a (4, 4) weighted group bicode.

The main advantage of this bicode is that both error detection and error correction is easy. Also in channels were retransmission is not very easy these codes can be used.

Now we proceed on to define the new notion of mixed weighted group bicode.

**DEFINITION 3.1.27:** *Let*

$$C = C_1 \cup C_2 = \left(C_1^1, C_2^1, \ldots, C_{n_1}^1\right) \cup \left(C_1^2, C_2^2, \ldots, C_{n_2}^2\right)$$

*be a group bicode if each code $C_j^1$ in $C_1$ is of weight n, $1 \leq j \leq n_1$ and each code $C_k^2$ in $C_2$ is of weight m, $1 \leq k \leq n_2$ then we call the group bicode to be a (m, n) weighted group bicode.*

We illustrate this by an example.

***Example 3.1.46:*** Let

$$C = C_1 \cup C_2 = \left\{C_1^1, C_2^1, C_3^1\right\} \cup \left\{C_1^2, C_2^2, C_3^2, C_4^2\right\}$$

be a group bicode where $C_1^1$ = {(0 0 0 0 0 0) (1 1 1 1 0 0) (0 0 1 1 1 1) (1 1 0 0 1 1)}, $C_2^1$ = {(1 1 1 1 0 0 0 0) (0 0 1 1 1 1 0 0) (0 0 0 0 0 0 0 0) (1 1 0 0 1 1 0 0)}, $C_3^1$ = {(0 0 0 0) (1 1 1 1)}; thus $C_1 = \left\{C_1^1, C_2^1, C_3^1\right\}$ is a weighted m group code here m = 4.

Take $C_2 = \left\{C_1^2, C_2^2, C_3^2, C_4^2\right\}$ where $C_1^2$ = {(1 1 1 1 1 1) (0 0 0 0 0 0)}, $C_2^2$ = {(1 1 1 1 1 1 0 0) (0 0 0 0 0 0 0 0)}, $C_3^2$ = {(1 1 1 1 1 1 0 0 0 0) (0 0 0 0 0 0 0 0 0 0) (0 0 0 1 1 1 1 1 1 0) (1 1 1 0 0 0 1 1 1 0)}, $C_4^2$ = {(0 0 0 0 0 0 1 1 1 1 1 1) (0 0 0 0 0 0 0 0 0 0 0 0) (0 0 0 1 1 1 1 1 1 0 0 0) (0 0 0 1 1 1 0 0 0 1 1 1)}.

Clearly $C_2$ is a group code which is a 6 weighted thus C = $C_1 \cup C_2$ is a (4, 6) weighted group bicode.

## 3.2 Set n-codes and their applications

Now having seen the notions of several types of bicodes we now proceed on to define different types of n-codes n ≥ 3 when n = 3 we can also call them as tricodes.

**DEFINITION 3.2.1:** *Let C = $C_1 \cup C_2 \cup \ldots \cup C_n$ (n ≥ 3) be such that each $C_i$ = $\left(C_1^i, C_2^i, \ldots, C_{n_i}^i\right)$ is a set code i = 1, 2, …, n. Then we call C to be a set n-code; when n = 3 we call C to be a set tri code.*

We will illustrate this by the following example.

***Example 3.2.1:*** Let

$$C = C_1 \cup C_2 \cup C_3 \cup C_4 \cup C_5$$
$$= \left\{C_1^1, C_2^1, C_3^1\right\} \cup \left\{C_1^2, C_2^2\right\} \cup \left\{C_1^3, C_2^3, C_3^3, C_4^3\right\}$$
$$\cup \left\{C_1^4, C_2^4, C_3^4\right\} \cup \left\{C_1^5, C_2^5\right\}$$

where $C_1^1$ = {(1 1 1 1 1) (0 0 0 0 0) (1 1 0 0 0) (1 0 1 0 1) (0 1 0 1 0)}, $C_2^1$ = {(0 0 0 0 0 0) (1 1 1 1 0 0) (0 1 1 1 1 0) (0 0 1 1 1 1) (1 1 0 0 1 1)}, $C_3^1$ = {(0 0 0 0 0 0 0) (1 1 1 0 0 1 1) (0 1 1 0 1 1 0) (0 1 0 1 0 1 0)}, $C_1^2$ = {(0 0 0 0 0 0) (1 0 1 0 1 0) (0 1 0 1 0 1) (1 1 1 0 0 0) (0 0 0 1 1 1) (1 1 0 0 1 1) (1 1 1 0 1 1)}, $C_2^2$ = {(0 0 0 0 0 0 0 0) (1 1 0 0 1 1 0 0) (1 1 1 1 0 0 0 0) (0 0 0 0 1 1 1 1) (1 1 0 0 0 0 1 1) (1 1 1 0 0 1 1 1) (0 1 1 0 1 1 1 0) (1 1 1 1 1 0 0 0) (0 0 0 1 1 1 1 1) (0 0 1 1 1 1 0 0)}, $C_1^3$ = {(0 1 0 0) (0 0 0 0) (1 0 1 0) (1 1 0 1)}, $C_2^3$ = {(0 0 0 0 0 0) (1 0 1 0 1 1) (1 0 1 1 1 0) (1 1 1 0 1 1) (1 0 1 1 1 1) (1 1 0 0 1 0)}, $C_3^3$ = {(0 0 0 0 0) (1 1 0 0 1) (1 1 1 0 0) (0 1 1 1 0) (0 0 1 1 0) (0 0 1 1 1) (1 1 0 0 0)}, $C_4^3$ = {(1 1 1 1 1 1 1) (0 0 0 0 0 0 0) (1 1 0 0 0 1 1) (0 0 1 1 1 0 0) (1 0 1 0 1 1 0) (0 1 0 1 1 1 1)}, $C_1^4$ = {(0 0 0 0 0 0 0 0) (1 1 1 1 0 0 1 1) (0 0 1 1 0 0 1 1) (0 0 0 0 1 1 1 1) (1 1 1 1 1 1 0 0) (0 0 1 1 1 1 1 1) (1 1 1 0 0 1 1 1)}, $C_2^4$ = {(1 1 1 1 1 1 1 1 1) (0 0 0 0 0 0 0 0 0)}, $C_3^4$ = {(1 1 1 1 1) (0 0 0 0 0)}, $C_1^5$ = {(0 0 0 0

0 0) (1 0 0 0 1 0) (0 1 1 0 0 1) (0 1 1 1 0 1) (1 0 0 1 1 0)} and $C_2^5$ = {(0 0 0 0 0 0 0 0) (1 1 0 0 0 1 1 1) (0 0 1 1 1 0 0 0) (1 1 0 1 1 0 0 1) (0 0 1 0 0 1 1 0)} is a set n-code n = 5.

***Example 3.2.2:*** Let

$$C = C_1 \cup C_2 \cup C_3 =$$
$$\left\{C_1^1, C_2^1, C_3^1\right\} \cup \left\{C_1^2, C_2^2, C_3^2, C_4^2\right\} \cup \left\{C_1^3, C_2^3\right\}$$

where

$C_1^1$ = {(1 1 1 0 0 0) (0 0 0 1 1 1) (0 0 0 1 1 1) (0 0 0 0 0 0) (1 0 1 0 1 0) (1 0 1 1 0 0) (0 1 0 1 1 0) (0 1 1 0 1 0)},

$C_2^1$ = {(0 0 0 0 0) (1 1 0 0 0) (0 0 1 1 1) (1 1 1 1 0) (0 1 1 1 1)},

$C_3^1$ = {(0 0 0 0 0 0 0) (1 1 1 0 0 0 1) (0 0 0 1 1 1 0) (1 1 1 1 1 0 0) (1 0 1 1 0 1 1) (1 1 1 1 1 1 1)},

$C_1^2$ = {(1 1 1 1 1 1) (0 0 0 0 0 0)},

$C_2^2$ = {(1 1 1 1) (0 1 0 1) (1 0 1 1) (0 0 0 0)},

$C_3^2$ = {(0 0 0 0 0 0 0 0) (1 1 1 0 0 0 0 0) (1 1 0 1 1 1 0 0) (1 1 1 1 1 1 1 1) (1 1 0 1 1 0 1 1) (1 1 1 0 1 1 1 0)},

$C_4^2$ = {(1 1 1 1 0 0 0 0 1) (0 0 0 0 0 0 0 0 0) (0 0 0 0 1 1 1 1 0) (1 1 1 0 0 0 1 1 1) (0 0 0 1 1 1 0 0 0) (1 1 1 1 1 0 0 0 1) (0 1 1 1 0 1 1 1 0) (1 1 1 1 1 1 0 0 0)},

$C_1^3$ = {(1 1 1 1 1 1 1) (0 0 0 0 0 0 0) (1 1 0 0 0 1 1) (0 0 1 1 1 0 0) (1 0 1 0 1 0 1)} and

$C_2^3$ = {(1 1 1 1 1 1) (0 1 1 1 0 1) (0 1 1 0 1 1) (1 0 1 0 1 0) (0 1 0 1 0 1) (0 0 0 0 0 0)}.

Thus C is a set tri code.
Now we proceed on to define repetition set n-code.

**DEFINITION 3.2.2:** *Let $C = C_1 \cup C_2 \cup \ldots \cup C_n$ where each $C_i$ is a repetition set code for i = 1, 2, …, n then we call C to be a set repetition n-code ($n \geq 3$).*

We illustrate this by a few examples.

***Example 3.2.3:*** Let $C = C_1 \cup C_2 \cup C_3 \cup C_4 \cup C_5 \cup C_6$ where $C_1 = \{ C_1^1 = \{(0\ 0\ 0\ 0\ 0)\ (1\ 1\ 1\ 1\ 1)\},\ C_2^1 = \{(0\ 0\ 0\ 0)\ (1\ 1\ 1\ 1)\},\ C_3^1 = \{(0\ 0\ 0\ 0\ 0\ 0)\ (1\ 1\ 1\ 1\ 1\ 1)\},\ C_4^1 = \{(0\ 0\ 0\ 0\ 0\ 0\ 0)\ (1\ 1\ 1\ 1\ 1\ 1\ 1)\}\}$

$C_2 = \{ C_1^2 = \{(1\ 1\ 1\ 1\ 1\ 1)\ (0\ 0\ 0\ 0\ 0\ 0)\},\ C_2^2 = \{(1\ 1\ 1\ 1\ 1\ 1\ 1\ 1)\ (0\ 0\ 0\ 0\ 0\ 0\ 0\ 0)\},\ C_3^2 = \{(0\ 0\ 0\ 0\ 0)\ (1\ 1\ 1\ 1\ 1)\}\}$

$C_3 = \{ C_1^3 = \{(1\ 1\ 1\ 1)\ (0\ 0\ 0\ 0)\},\ C_2^3 = \{(1\ 1\ 1\ 1\ 1)\ (0\ 0\ 0\ 0\ 0)\},\ C_3^3 = \{(0\ 0\ 0\ 0\ 0\ 0\ 0\ 0\ 0)\ (1\ 1\ 1\ 1\ 1\ 1\ 1\ 1\ 1)\}\}$

$C_4 = \{ C_1^4 = \{(1\ 1\ 1\ 1\ 1)\ (0\ 0\ 0\ 0\ 0)\},\ C_2^4 = \{(0\ 0\ 0\ 0\ 0\ 0\ 0\ 0\ 0\ 0)\ (1\ 1\ 1\ 1\ 1\ 1\ 1\ 1\ 1\ 1)\},\ C_3^4 = \{(1\ 1\ 1\ 1\ 1\ 1)\ (0\ 0\ 0\ 0\ 0\ 0)\}\}$

$C_5 = \{ C_1^5 = \{(0\ 0\ 0\ 0\ 0\ 0\ 0\ 0\ 0)\ (1\ 1\ 1\ 1\ 1\ 1\ 1\ 1\ 1)\},\ C_2^5 = \{(1\ 1\ 1\ 1)\ (0\ 0\ 0\ 0\ 0)\},\ C_3^5 = \{(1\ 1\ 1\ 1\ 1\ 1)\ (0\ 0\ 0\ 0\ 0\ 0)\}\}$

$C_6 = \{ C_1^6 = \{(1\ 1\ 1\ 1\ 1\ 1)\ (0\ 0\ 0\ 0\ 0\ 0)\},\ C_2^6 = \{(0\ 0\ 0\ 0\ 0\ 0\ 0\ 0)\ (1\ 1\ 1\ 1\ 1\ 1\ 1\ 1)\},\ C_3^6 = \{(0\ 0\ 0\ 0\ 0\ 0\ 0\ 0\ 0)\ (1\ 1\ 1\ 1\ 1\ 1\ 1\ 1\ 1)\},\ C_4^6 = \{(0\ 0\ 0\ 0\ 0\ 0\ 0)\ (1\ 1\ 1\ 1\ 1\ 1\ 1)\},\ C_5^6 = \{(0\ 0\ 0\ 0\ 0\ 0\ 0\ 0\ 0\ 0)\ (1\ 1\ 1\ 1\ 1\ 1\ 1\ 1\ 1\ 1)\}\}$

is a set repetition 6-code.

The main advantage in using these set repetition n-code is that both error detection and error correction is very easy. Further these n-codes can be used in cryptography by using several of the codes as misleading codes and only the receiver and the sender alone know the code words which carry the message so it is very difficult for the intruder to break the key.

The main use of set n-codes is that when we want to mislead the intruder it is best suited. Out of these set of n-codes $C_1, \ldots, C_n$ we can use one or two codes say some $C_{j_t}^i$ for some i in {1, 2, …, n} to be a message transmitter the rest being to mislead for in this case we do not want the codes to be semigroups or groups it is sufficient if they are a proper subset of a code which forms a group under addition. It is left as an experiment for the interested to use these codes in cryptography. We will just illustrate this situation by an example.

***Example 3.2.4:*** Let

$$C = C_1 \cup C_2 \cup C_3 \cup C_4$$
$$= \left\{C_1^1, C_2^1, C_3^1\right\} \cup \left\{C_1^2, C_2^2, C_3^2, C_4^2\right\} \cup \left\{C_1^3, C_2^3\right\} \cup \left\{C_1^4, C_2^4, C_3^4\right\}$$

be a set 4-code where

$C_1^1$ = {(0 0 0 0 0) (1 1 1 1 1) (1 1 0 0 0) (1 0 1 0 1)}

$C_2^1$ = {(0 0 0 0 0 0) (1 1 1 1 0 0) (1 1 0 0 1 1) (0 0 1 1 1 1)}

$C_3^1$ = {(1 1 1 0) (0 1 1 1) (0 1 0 1) (0 0 0 0)}

$C_1^2$ = {(0 0 0 0 0 0) (1 0 1 0 1 0) (0 1 0 1 0 1) (1 1 0 0 0 1) (1 1 1 0 0 1)}

$C_2^2$ = {(0 0 0 0 0 0 0 0) (1 1 0 0 1 1 0 1) (0 1 1 0 1 1 0 1) (1 1 1 1 0 0 0 1) (1 0 0 0 1 1 1 0)}

$C_3^2$ = {(0 0 0 0 0 0 0) (1 1 1 0 1 0 0) (0 1 1 1 0 1 0) (0 0 1 1 1 0 1) (1 0 0 1 1 1 0) (1 1 0 0 1 1 1) (0 1 0 1 0 1 1) (1 0 1 0 0 1 1)}

$C_4^2$ = {(1 1 1 1 0 0 0 0 0) (1 1 1 0 0 0 1 1 1) (0 0 0 1 1 1 0 0 1) (0 0 0 0 0 0 0 0 0) (1 0 1 0 1 0 1 0 1)}

$C_1^3$ = {(0 0 0 0 0 0 0) (1 1 0 1 0 0 0) (0 1 1 0 1 0 0) (0 0 1 1 0 1 0) (0 0 0 1 1 0 1) (1 0 1 1 1 0 0) (0 1 0 1 1 1 0) (0 0 1 0 1 1 1) (1 1 0 0 1 0 1) (1 1 1 0 0 1 0) (0 1 1 1 0 0 1) (1 0 0 0 1 1 0) (0 1 0 0 0 1 1) (1 0 1 0 0 0 1) (1 1 1 1 1 1 1) (1 0 0 1 0 1 1)}

$C_2^3$ = {(0 0 0 0 0 0) (1 1 0 1 0 0) (1 0 0 1 1 0) (1 0 1 0 1 0) (0 1 0 1 0 1) (0 1 1 0 1 0), (0 1 1 1 0 0)}

$C_1^4$ = {(1 1 1 1 1 1) (0 0 0 0 0 0) (1 1 1 1 0 0) (1 1 0 0 1 1) (1 0 1 1 0 0) (0 1 1 1 1 0) (0 0 1 1 1 1) (1 1 1 0 0 1)}

$C_2^4$ = {(0 0 0 0 0 0) (0 0 1 0 0 1) (1 0 0 1 0 0) (0 1 0 0 1 0) (0 1 1 0 1 1) (1 1 1 1 1 1) (1 1 0 1 1 0) (1 0 1 1 0 1)} and

$C_3^4$ = {(1 1 1 1 1 1 1 1) (0 0 0 0 0 0 0 0) (1 1 0 1 1 0 1 1) (0 1 1 1 1 0 0 1) (0 1 1 0 0 0 1 1) (0 0 1 1 1 1 1 1) (1 0 1 0 1 0 1 0) (0 1 0 1 0 1 0 1)}.

Now $C = C_1 \cup C_2 \cup C_3 \cup C_4$ is a set 4-code we see in the codes $C_3^2$, $C_1^3$ and $C_1^4$ are usual codes that is they are closed with respect to addition and form a group. Further they are subspaces of a vector space.

Now any message x ∈ C = $C_1 \cup C_2 \cup C_3 \cup C_4$ would be of the form x = {(1 1 0 0 0) (1 1 1 1 0 0) (1 1 1 0)} ∪ {(1 1 0 0 0 1) (1 1 1 1 0 0 0 1) (0 1 1 1 0 1 0) (1 0 1 0 1 0 1 0 1)} ∪ {(0 1 1 0 1 0 0) (0 1 1 0 1 0)} ∪ {(1 1 0 0 1 1) (0 0 1 0 0 1) (1 1 1 1 1 1 1 1)} only the code words (0 0 1 1 1 0 1), (0 1 1 0 1 0 0) and (1 1 0 0 1 1) alone need to be decoded for all other code words are just to mislead. The receiver and the sender know these facts. But a general intruder will not know these facts and it is a remote impossiblility he guesses this.

Now if we take set codes of same weight we call them to m-weighted set codes. We will define this situation and then illustrate it by an example.

The main use of such weighted set n-codes is that it is easy to detect errors. These codes can be thought of an easy error detectable codes.

**DEFINITION 3.2.3:** *Let*

$$C = C_1 \cup C_2 \cup \ldots C_n$$
$$= \left\{C_1^1, C_2^1, \ldots C_{r_1}^1\right\} \cup \left\{C_1^2, C_2^2 \ldots C_{r_2}^2\right\} \cup \ldots \cup \left\{C_1^n, C_2^n \ldots C_{r_n}^n\right\}$$

*be a set n-code if the weight of each code word in each* $C_{j_t}^i$ *in* $C_i$ *is of weight m for* $1 \leq j_i \leq r_i$, *i = 1, 2, …, n then we call C to be a m-weighted set n-code.*

***Example 3.2.5:*** Let C = $C_1 \cup C_2 \cup C_3 \cup C_4 \cup C_5$ where

$$C_1 = \left\{C_1^1, C_2^1, C_3^1\right\} \text{ where}$$

$C_1^1$ = {(1 1 1 0 0) (0 1 1 1 0) (0 0 0 0 0) (1 0 1 0 1) (1 1 0 0 1) (1 0 1 1 0)}

$C_2^1$ = {(1 1 1 0 0 0) (0 1 1 1 0 0) (1 0 1 1 0 0) (0 1 0 1 0 1) (1 0 1 0 1 0) (0 0 0 0 0 0)}

$C_3^1$ = {(0 0 0 1 0 0 1 1) (1 1 1 0 0 0 0 0) (0 0 0 0 0 0 0 0) (1 0 1 0 1 0 0 0) (1 0 0 1 0 0 0 1) (0 0 1 1 1 0 0 0) (0 1 0 1 1 0 0 0) (0 0 0 0 0 1 1 1)}.

$$C_2 = \left\{C_1^2, C_2^2\right\}$$

$C_1^2$ = {(1 1 0 1) (1 0 1 1) (1 1 1 0) (0 1 1 1) (0 0 0 0)}

$C_2^2$ = {1 1 0 0 0 1) (1 1 0 0 1 0) (1 1 0 1 0 0) (0 1 1 0 1 0) (0 1 1 0 0 1) (0 0 1 1 0 1) (0 0 0 0 0 0)}

$$C_3 = \left\{C_1^3, C_2^3, C_3^3, C_4^3\right\}$$

where

$C_1^3$ = {(1 0 1 0 1 0) (0 1 0 1 0 1) (0 0 0 0 0 0), (1 0 1 1 0 0) (0 1 0 1 1 0) (1 1 0 1 0 0)}

$C_2^3$ = {(0 0 0 0 0 0 0) (1 1 1 0 0 0 0) (0 1 1 1 0 0 0) (0 0 1 1 1 0 0) (0 0 0 1 1 1 0) (0 0 0 0 1 1 1)}

$C_3^3$ = {(0 0 0 0 0 0 0 0) (1 1 1 0 0 0 0 0) (0 0 0 0 1 1 1 0) (0 0 0 0 0 1 1 1)} and

$C_4^3$ = {(0 0 0 0 0 0 1 1 1) (0 0 0 0 0 0 0 0 0) (0 0 0 0 0 1 1 1 0) (0 0 0 0 1 1 1 0 0) (0 0 0 1 1 1 0 0 0) (0 0 1 1 1 0 0 0 0) (0 1 1 1 0 0 0 0 0) (1 1 1 0 0 0 0 0 0)}.

$$C_4 = \left\{C_1^4, C_2^4, C_3^4\right\}$$

where

$C_1^4$ = {(1 1 0 0 0 0 1) (1 1 0 0 0 1 0) (0 0 0 0 0 0 0) (1 1 0 0 0 1 0 0) (1 1 0 1 0 0 0) (0 1 1 0 1 0 0) (0 0 1 1 0 1 0) (0 0 0 1 1 0 1) (1 0 0 0 0 1 1) (0 1 0 0 0 1 1) (0 0 0 1 0 1 1) (1 0 0 0 1 1 0) (1 0 1 1 0 0 0)}

$C_2^4$ = {(1 1 0 0 0 0 0 1) (1 0 1 1 0 0 0 0) (1 1 0 0 0 0 1 0) (1 1 0 0 0 1 0 0) (1 1 0 0 1 0 0 0) (1 1 0 1 0 0 0 0) (1 0 0 1 1 0 0 0) (1 0 0 0 1 1 0 0) (1 0 0 0 0 1 1 0) (1 0 0 0 0 0 1 1) (0 1 0 0 0 0 1 1) (0 0 1 0 0 0 1 1) (0 0 0 1 0 0 1 1) (0 0 0 0 1 0 1 1) (0 1 1 0 1 0 0 0) (0 1 1 0 0 1 0 0) (0 1 1 0 0 0 1 0) (0 1 1 0 0 0 0 1) (0 0 0 0 0 0 0 0)} and

$C_3^4$ = {(1 1 1 0 0 0) (0 0 0 0 0 0) (0 1 1 1 0 0) (0 0 1 1 1 0) (0 0 0 1 1 1) (1 0 1 0 1 0) (0 1 0 1 0 1)}.

$$C_5 = \left\{C_1^5, C_2^5\right\}$$

where

$C_1^5$ = {(0 0 0 0 0 1 1 1) (1 0 1 0 1 0 0 0) (0 1 0 1 0 1 0 0) (0 0 1 0 1 0 0 1) (0 0 0 1 0 1 0 1) (0 0 0 0 0 0 0 0) (1 1 1 0 0 0 0 0) (0 1 1 1 0 0 0 0) (0 0 1 1 1 0 0 0) (0 0 0 1 1 1 0 0) (0 0 0 0 1 1 1 0)} and

$C_2^5$ = {(0 0 0 0 0 0 0) (1 0 1 0 1 0 0) (0 1 0 1 0 1 0) (0 0 1 0 1 0 1) (1 0 0 1 0 1 0) (0 1 0 0 1 0 1) (1 0 1 0 0 1 0) (0 1 0 1 0 0 1)}.

Clearly C is 3 weighted 5-code. This code can also be used in cryptography.

These codes are best suited in networking in computer as well as fragmenting for the set n-code, $C = C_1 \cup C_2 \cup ... \cup C_n$ can be fragmented at any $C_i$ according to need. Also these codes are of different lengths it will find its use in networking. Thus it is left for the computer scientists to use these set n-codes. Also they can use these set n-codes in places were the intruder should not easily hack the privacy.

Now we proceed on to define the notion of group n-codes for when we make use of the binary symbols i.e., {0, 1} the notion of semigroup and group coincide under addition.

We proceed on to define group n-codes.

**DEFINITION 3.2.4:** *Let*

$$C = C_1 \cup C_2 \cup ... C_n$$
$$= \left\{C_1^1, C_2^1, \ldots, C_{r_1}^1\right\} \cup ... \cup \left\{C_1^n, C_2^n ... C_{r_n}^n\right\}$$

*($n \geq 3$) where each is a group code, $1 \leq j_i \leq r_i$; $i = 1, 2, \ldots, n$. We call C the set group n-code.*

*When n = 2 we get the group bicode. When n = 3 we get the group tricode.*

We will illustrate this by some examples.

***Example 3.2.6:*** Let

$$C = C_1 \cup C_2 \cup C_3 \cup C_4 \cup C_5$$
$$= \left\{C_1^1, C_2^1, C_3^1\right\} \cup \left\{C_1^2, C_2^2, C_3^2, C_4^2\right\} \cup \left\{C_1^3, C_2^3\right\} \cup$$
$$\left\{C_1^4, C_2^4, C_3^4\right\} \cup \left\{C_1^5, C_2^5, C_3^5, C_4^5\right\}$$

where

$C_1^1$ = {(0 0 0 0 0) (1 1 1 1 1)}

$C_2^1$ = {(0 0 0 0 0 0 0) (1 0 0 0 1 0 1) (0 1 0 0 1 1 1) (0 0 1 0 1 1 0) (0 0 0 1 0 1 1), (1 1 0 0 0 1 0) (1 0 1 0 0 1 1) (1 0 0 1 1 1 0) (0 1 1 0 0 0 1) (0 1 0 1 1 0 0) (0 0 1 1 1 0 1) (1 1 1 0 1 0 0) (0 1 1 1 0 1 0) (1 1 0 1 0 0 0) (1 0 1 1 0 0 0) (0 1 1 1 0 1 0) (1 1 1 1 1 1 1)}

$C_3^1$ = {(0 0 0 0) (1 0 1 1) (0 1 0 1) (1 1 1 0)}

$C_1^2$ = {(0 0 0 0 0 0) (0 0 1 0 0 1) (0 1 0 0 1 0) (0 1 1 0 1 1) (1 1 1 1 1 1) (1 1 0 1 1 0) (1 0 1 1 0 1) (1 0 0 1 0 0)}

$C_2^2$ = {(0 0 0 0 0 0 0 0) (1 1 1 1 1 1 1)}

$C_3^2$ = {(0 0 0 0 0) (0 1 1 1 1) (1 0 0 1 0) (1 1 1 0 1)}

$C_4^2$ = {(0 0 0 1 1 1 1) (0 1 1 0 0 1 1) (1 0 1 0 1 0 1) (0 0 0 0 0 0 0) (0 1 1 1 1 0 0) (1 1 0 0 1 1 0) (1 0 1 1 0 1 0) (1 1 0 1 0 0 1)}

$C_1^3$ = {(0 0 0 0 0 0 0 0 0) (1 1 1 1 1 1 1 1)}

$C_2^3$ = {(0 0 0 0) (1 1 1 0) (0 1 1 0) (0 0 1 1) (1 0 0 0) (0 1 0 1) (1 1 0 1) (1 0 1 1)}

$C_1^4$ = {(1 1 1 1 1 1) (0 0 0 0 0 0)}

$C_2^4$ = {(0 0 0 0) (1 0 1 1) (0 1 1 1) (1 0 0 1) (1 1 0 0) (0 0 1 0) (1 1 1 0) (0 1 0 1)}

$C_3^4$ = {(0 0 0 0 0) (0 1 0 0 1) (0 0 1 0 1) (1 0 0 1 1) (0 1 1 0 0) (1 0 1 1 0) (1 1 0 1 0) (1 1 1 1 1)}

$C_1^5$ = {(0 0 0 0 0 0 0) (1 1 1 1 1 1 1)}

$C_2^5$ = {(0 0 0 0 0), (1 1 1 1 1)}

$C_3^5$ = {(0 0 0 1 1 1 1), (0 1 1 0 0 1 1) (1 0 1 0 1 0 1) (0 0 0 0 0 0 0), (0 1 1 1 1 0 0) (1 1 0 0 1 1 0) (1 0 1 1 0 1 0) (1 1 0 1 0 0 1)}

$C_4^5$ = {(0 0 0 0) (1 0 1 1) (0 1 0 1) (1 1 1 0)}

is a set group 5-code any typical group 5 code word X $\in$ C = X = $X_1 \cup X_2 \cup X_3 \cup X_4 \cup X_5$
= {(1 1 1 1 1) (1 0 0 0 1 0 1) (1 0 1 1)} $\cup$ {(0 0 1 0 0 1) (1 1 1 1 1 1 1 1) (1 1 1 0 1) (0 0 0 1 1 1 1)} $\cup$ {(1 1 1 1 1 1 1 1), (1 1 1

0)} ∪ {(0 0 0  0 0 0) (1 0 1 1) (0 1 0 0 1)} ∪ {(1 1 1 1 1 1 1) (1 1 1 1), (0 0 0 1 1 1 1) (0 1 0 1)}.

These codes can be used in computer networking when fragmenting of codes is to be carried out these codes can be best suited as at the union the codes can be fragmented further these codes are of varied lengths which is also an added advantage to the user.

Let us now define different types of group n-codes.

**DEFINITION 3.2.5:** *Let*

$$C = C_1 \cup C_2 \cup \ldots C_n$$

$$= \left\{C_1^1, C_2^1, \ldots, C_{r_1}^1\right\} \cup \left\{C_1^2, C_2^2, \ldots, C_{r_2}^2\right\} \cup \ldots \cup \left\{C_1^n, C_2^n, \ldots, C_{r_n}^n\right\}$$

*be a group n-code (n ≥ 3) if each of the codes* $C_{j_i}^i$ *are repetition codes,* $1 \leq j_i \leq r_i$*, i = 1, 2, …, n; then we call C to be a repetition group n-code.*

We will illustrate this situation by some examples.

***Example 3.2.7:*** Let

$$C = C_1 \cup C_2 \cup C_3 \cup C_4 =$$

$$\left\{C_1^1, C_2^1, C_3^1\right\} \cup \left\{C_1^2, C_2^2\right\} \cup \left\{C_1^3, C_2^3, C_3^3, C_4^3\right\} \cup \left\{C_1^4, C_2^4, C_3^4\right\}$$

where $C_1^1$ = {(0 0 0 0 0) (1 1 1 1 1)}, $C_2^1$ = {(1 1 1 1) (0 0 0 0)}, $C_3^1$ = {(0 0 0 0 0 0 0 0 0) (1 1 1 1 1 1 1 1 1}, $C_1^2$ = {(1 1 1 1 1 1) (0 0 0 0 0 0)}, $C_2^2$ = {(0 0 0 0) (1 1 1 1)}, $C_1^3$ = {(0 0 0 0 0 0) (1 1 1 1 1 1)}, $C_2^3$ = {(0 0 0 0 0 0 0) (1 1 1 1 1 1 1)}, $C_3^3$ = {(1 1 1 1 1) (0 0 0 0 0)}, $C_4^3$ = {(1 1 1 1 1 1 1 1 1 1) (0 0 0 0 0 0 0 0 0 0)}, $C_1^4$ = {(1 1 1 1) (0 0 0 0)}, $C_2^4$ = {(0 0 0 0 0 0 0 0) (1 1 1 1 1 1 1 1)} and $C_3^4$ = {(1 1 1 1 1 1 1) (0 0 0 0 0 0 0)}.

It is easily verified that C is a repetition group 4 code.

We give yet another example of a repetition group tricode.

***Example 3.2.8:*** Let C = $C_1 \cup C_2 \cup C_3$

$$= \left\{C_1^1, C_2^1, C_3^1, C_4^1\right\} \cup \left\{C_1^2, C_2^2\right\} \cup \left\{C_1^3, C_2^3, C_3^3\right\}$$

where $C_1^1$ = { (0 0 0 0) (1 1 1 1)}, $C_2^1$ = {(0 0 0 0 0 0) (1 1 1 1 1 1)}, $C_3^1$ = {(1 1 1 1 1 1 1 1) (0 0 0 0 0 0 0 0)}, $C_4^1$ = {(0 0 0 0 0) (1 1 1 1 1)}, $C_1^2$ = {(0 0 0 0 0 0) (1 1 1 1 1 1)}, $C_2^2$ = {(1 1 1 1 1 1 1 1 1) (0 0 0 0 0 0 0 0 0)}, $C_1^3$ = {(0 0 0 0 0 0) (1 1 1 1 1 1)} and $C_2^3$ = {(1 1 1 1 1 1 1 1 1 1) (0 0 0 0 0 0 0 0 0 0)}. C is a repetition group tricode.

These codes can be used when it is not possible to retransmit message or the cost is high; also these codes can carry at a time several messages unlike the repetition code or the repetition bicode. With the advent of computers these group n-codes can be transmitted with space crafts or satellites.

Now we proceed on to define the group parity check n-matrix and the group generator n-matrix associated with a group n-code and illustrate them by examples.

**DEFINITION 3.2.6:** *Let*

$$C = C_1 \cup C_2 \cup \ldots \cup C_n$$

$$= \left\{C_1^1, C_2^1, \ldots, C_{r_1}^1\right\} \cup \left\{C_1^2, C_2^2, \ldots, C_{r_2}^2\right\} \cup \ldots \cup \left\{C_1^n, C_2^n, \ldots, C_{r_n}^n\right\}$$

*be a group n-code. We call*

$$H = \left\{H_1^1, H_2^1, \ldots H_{r_1}^1\right\} \cup \left\{H_1^2, H_2^2, \ldots H_{r_2}^2\right\} \cup \ldots \cup \left\{H_1^n, H_2^n, \ldots, H_{r_n}^n\right\}$$

*where each $H_{j_i}^i$ is the parity check matrix of the code $C_{j_i}^i$ ; $1 \le j_i \le r_i$; i = 1, 2, …, n to be the group parity check n-code. If*

$$G = \left\{G_1^1, G_2^1, \ldots G_{r_1}^1\right\} \cup \left\{G_1^2, G_2^2, \ldots, G_{r_2}^2\right\} \cup \ldots \cup \left\{G_1^n, G_2^n, \ldots, G_{r_n}^n\right\}$$

*is such that $G_{j_i}^i$ is the generator matrix of the code $C_{j_i}^i$ ;$1 \le j_i \le r_i$; i = 1, 2, …, n then we call G to be the group generator n-matrix of the group n-code C.*

Let us illustrate this by an example.

***Example 3.2.9:*** Let

$$C = C_1 \cup C_2 \cup C_3 \cup C_4$$

$$= \left\{C_1^1, C_2^1, C_3^1\right\} \cup \left\{C_1^2, C_2^2\right\} \cup \left\{C_1^3, C_2^3, C_3^3, C_4^3\right\} \cup \left\{C_1^4, C_2^4, C_3^4\right\}$$

be a group 4 code associated with the group parity check 4-matrix

$$H = H_1 \cup H_2 \cup H_3 \cup H_4 =$$

$$\left\{\begin{pmatrix} 1 & 1 & 0 & 0 & 0 \\ 1 & 0 & 1 & 0 & 0 \\ 1 & 0 & 0 & 1 & 0 \\ 1 & 0 & 0 & 0 & 1 \end{pmatrix} = H_1^1,\right.$$

$$H_2^1 = \begin{pmatrix} 0 & 1 & 1 & 1 & 0 & 0 \\ 1 & 0 & 1 & 0 & 1 & 0 \\ 1 & 1 & 0 & 0 & 0 & 1 \end{pmatrix},$$

$$\left. H_3^1 = \begin{pmatrix} 0 & 0 & 1 & 0 & 1 & 1 & 1 \\ 0 & 1 & 0 & 1 & 1 & 1 & 0 \\ 1 & 0 & 1 & 1 & 1 & 0 & 0 \end{pmatrix}\right\} \cup$$

$$\left\{\begin{pmatrix} 1 & 1 & 0 & 0 \\ 1 & 0 & 1 & 0 \\ 1 & 0 & 0 & 1 \end{pmatrix} = H_1^2,\right.$$

$$\left. H_2^2 = \begin{pmatrix} 1 & 1 & 1 & 0 & 1 & 0 & 0 \\ 0 & 1 & 1 & 1 & 0 & 1 & 0 \\ 1 & 1 & 0 & 1 & 0 & 0 & 1 \end{pmatrix}\right\} \cup$$

$$\{ H_1^3 = \begin{pmatrix} 1 & 1 & 1 & 0 & 1 & 0 & 0 \\ 0 & 1 & 1 & 1 & 0 & 1 & 0 \\ 1 & 1 & 0 & 1 & 0 & 0 & 1 \end{pmatrix},$$

$$H_2^3 = \begin{pmatrix} 1 & 1 & 0 & 0 & 0 & 0 \\ 1 & 0 & 1 & 0 & 0 & 0 \\ 1 & 0 & 0 & 1 & 0 & 0 \\ 1 & 0 & 0 & 0 & 1 & 0 \\ 1 & 0 & 0 & 0 & 0 & 1 \end{pmatrix},$$

$$H_3^3 = \begin{pmatrix} 1 & 0 & 0 & 1 & 0 & 0 \\ 0 & 1 & 0 & 0 & 1 & 0 \\ 0 & 0 & 1 & 0 & 0 & 1 \end{pmatrix},$$

$$H_4^3 = \left.\begin{pmatrix} 1 & 1 & 0 & 0 & 0 & 0 & 0 & 0 \\ 1 & 0 & 1 & 0 & 0 & 0 & 0 & 0 \\ 1 & 0 & 0 & 1 & 0 & 0 & 0 & 0 \\ 1 & 0 & 0 & 0 & 1 & 0 & 0 & 0 \\ 1 & 0 & 0 & 0 & 0 & 1 & 0 & 0 \\ 1 & 0 & 0 & 0 & 0 & 0 & 1 & 0 \\ 1 & 0 & 0 & 0 & 0 & 0 & 0 & 1 \end{pmatrix}\right\} \cup$$

$$\{ H_1^4 = \begin{pmatrix} 1 & 1 & 0 & 0 & 0 & 0 & 0 \\ 1 & 0 & 1 & 0 & 0 & 0 & 0 \\ 1 & 0 & 0 & 1 & 0 & 0 & 0 \\ 1 & 0 & 0 & 0 & 1 & 0 & 0 \\ 1 & 0 & 0 & 0 & 0 & 1 & 0 \\ 1 & 0 & 0 & 0 & 0 & 0 & 1 \end{pmatrix},$$

$$H_2^4 = \{(1\ 1\ 1\ 1\ 1\ 1\ 1\ 1),$$

$$H_3^4 = \left.\begin{pmatrix} 0 & 0 & 1 & 0 & 1 & 1 & 1 \\ 0 & 1 & 0 & 1 & 1 & 1 & 0 \\ 1 & 0 & 1 & 1 & 1 & 0 & 0 \end{pmatrix}\right\}$$

is the group parity check n-matrix.

The elements of a group n-code will be known as the group n-code word. Suppose

$$x = \left(x_1^1, x_2^1, x_3^1\right) \cup \left(x_1^2, x_2^2\right) \cup \left(x_1^3, x_2^3, x_3^3, x_4^3\right) \cup \left(x_1^4, x_2^4, x_3^4\right)$$

is a group n-code word. Then we say x ∈ C if and only if

$$Hx^t = \left(H_1^1\left(x_1^1\right)^t, H_2^1\left(x_2^1\right)^t, H_3^1\left(x_3^1\right)^t,\right) \cup \left(H_1^2\left(x_1^2\right)^t, H_2^2\left(x_2^2\right)^t\right)$$
$$\cup \left(H_1^3\left(x_1^3\right)^t, H_2^3\left(x_2^3\right)^t, H_3^3\left(x_3^3\right)^t, H_4^3\left(x_4^3\right)^t\right)$$
$$\cup \left(H_1^4\left(x_1^4\right)^t, H_2^4\left(x_2^4\right)^t, H_3^4\left(x_3^4\right)^t\right)$$

= {(0 0 0 0) (0 0 0) (0 0 0)} ∪ {(0 0 0) (0 0 0)} ∪ {(0 0 0) (0 0 0 0 0) (0 0 0) (0 0 0 0 0 0 0)} ∪ {(0 0 0 0 0 0) (0) (0 0 0)},

which will be known as the group zero n-code word.

Now we use the properties of group parity check n-matrix for error detection, for which we define the notion of group n-syndrome of a group n-code.

**DEFINITION 3.2.7:** *Let*

$$C = C_1 \cup C_2 \cup \ldots \cup C_n$$
$$= \left\{C_1^1, C_2^1, \ldots, C_{r_1}^1\right\} \cup \left\{C_1^2, C_2^2, \ldots, C_{r_2}^2\right\} \cup \ldots \cup \left\{C_1^n, C_2^n, \ldots, C_{r_n}^n\right\}$$

*be a group n-code with the associated group parity check n-matrix,*

$$H = H_1 \cup H_2 \cup \ldots \cup H_n$$
$$= \left\{H_1^1, H_2^1, \ldots, H_{r_1}^1\right\} \cup \left\{H_1^2, H_2^2, \ldots, H_{r_2}^2\right\} \cup \ldots \cup \left\{H_1^n, H_2^n, \ldots, H_{r_n}^n\right\}.$$

*For any y ∈ C = $C_1 \cup C_2 \cup \ldots \cup C_n$ we define the group n-syndrome of y to be*

$$GS(y) = \left\{S\left(y_1^1\right), S\left(y_2^1\right), \ldots, S\left(y_{r_1}^1\right)\right\} \cup$$
$$\left\{S\left(y_1^2\right), S\left(y_2^2\right), \ldots, S\left(y_{r_2}^2\right)\right\} \cup \ldots \cup \left\{S\left(y_1^n\right), S\left(y_2^n\right), \cdots S\left(y_{r_n}^n\right)\right\}.$$

*If GS (y) = 0 $\cup \ldots \cup$ 0 then we accept y $\in$ C if GS(y) $\neq$ (0 … 0) then we say y $\notin$ C we have detected error and we have to correct it.*

As in case of usual codes we make use of the coset leader properties to correct the error or in other words find the approximately correct n-code word or the most likely transmitted n-code word.

Now for this we recall the notion of group n-vector space

$$V = V_1 \cup V_2 \cup \ldots \cup V_n$$
$$= \left\{V_1^1, V_2^1, \ldots, V_{r_1}^1\right\} \cup \left\{V_1^2, V_2^2, \ldots, V_{r_2}^2\right\} \cup \ldots \cup \left\{V_1^n, V_2^n, \ldots, V_{r_n}^n\right\}$$

where each $V_{j_i}^i$ is a vector space over the group $F_2 = \{0, 1\}$ under addition modulo 2; $i \leq j_i \leq r_i$ and i = 1, 2, …, n.

For a group n vector space $V = V_1 \cup V_2 \cup \ldots \cup V_n$ let $C_1 \cup C_2 \cup \ldots \cup C_n$ be a proper n-subset of V ie. each $C_i$ is a proper subset of $V_i$, $1 \leq i \leq n$ i.e., each $C_{j_i}^i$ is a proper subset of $V_{j_i}^i$ and further if each $C_{j_i}^i$ is a subspace of $V_{j_i}^i$; $1 \leq j_i \leq r_i$; i = 1, 2, …, n then we call C to be a group n-subspace of V.

The n factor space V/C consists of n cosets a + c =

$$a_1 + c_1 \cup a_2 + c_2 \cup \ldots \cup a_n + c_n$$
$$= \left(a_1^1 + c_1^1, a_2^1 + c_2^1 \cup \cdots \cup a_{r_1}^1 + c_{r_1}^1\right) \cup$$
$$\left(a_1^2 + c_1^2, a_2^2 + c_2^2 \cup \cdots \cup a_{r_2}^2 + c_{r_2}^2\right) \cup \ldots \cup$$
$$\left(a_1^n + c_1^n, a_2^n + c_2^n \cup \cdots \cup a_{r_n}^n + c_{r_n}^n\right);$$

i.e., for each $c^i_{j_i}$ $a^i_{j_i}$ + $c^i_{j_i}$ = $\left\{a^i_{j_i} + x^i_{j_i} \middle| x^i_{j_i} \in c^i_{j_i}\right\}$ for arbitrary $a^i_{j_i} \in F_2^{p_i}$ ; $1 \le j_i \le r_i$, i = 1, 2, …, n. where $p_i$ is the length of the code words in $c^i_{j_i}$ .

This is a partition of

$$\left(F^{p_1^1}, F^{p_2^1}, \ldots, F^{p_{r_1}^1}\right) \cup \left(F^{p_1^2}, F^{p_2^2}, \ldots, F^{p_{r_2}^2}\right)$$

$$\cup \ldots \cup \left(F^{p_1^n}, F^{p_2^n}, \ldots, F^{p_{r_n}^n}\right)$$

of the form $F_1 \cup F_2 \cup \ldots \cup F_n$

$$= \left\{c_1^1 \cup \left(a_1^1 + c_1^1\right) \cup \ldots \cup \left(a_1^{t_1} + c_1^1\right), c_2^1 \cup \left(a_1^2 + c_2^2\right) \cup \ldots \cup\right.$$

$$\left(a_1^{t_2} + c_2^1\right), \ldots, c^1 \cup \left(a_1^2 + c_2^2\right) \ldots$$

$$\cup \left(a_1^{t_2^1} + c_2^1\right), \ldots, c_{r_1}^1 \cup \left(a_{r_1}^1 + c_{r_1}^1\right) \cup \ldots.$$

$$\cup \left.\left(a_1^{t_{r_1}^1} + c_{r_1}^1\right) \mid t_i^1 = 2^{n_i^1 - k_i^1} - 1\right\} \cup \ldots$$

$$\cup \left\{c_1^n \cup \left(a_1^n + c_1^n\right) \cup \ldots \cup \left(a_1^{t_1^n} + c_1^n\right),\right.$$

$$\left\{c_2^n \cup \left(a_2^n + c_2^n\right) \cup \ldots \cup \left(a_2^{t_2^n} + c_2^n\right),\right.$$

$$\left.\ldots, c_{r_n}^n \cup \left(a_{r_n}^n + c_{r_n}^n\right) \cup \ldots \cup \left(a_n^{t_2^n} + c_2^n\right) \mid t_i^n = 2^{n_i^1 - k_i^1} - 1\right\}$$

i.e., each code $c^i_{j_i}$ is of length $n^i_{j_i}$ with $k^i_{j_i}$ message symbols of $c^i_{j_i}$ is a ( $n^i_{j_i}$ , $k^i_{j_i}$ ) code true for $1 \le j_i \le r_i$. i = 1, 2, …, n.

Clearly if y is any received group n-code then y must be an element of one of these n-cosets say

$$\left(a_1^{t_1^1} + c_1^1, a_1^{t_2^1} + c_2^1, \ldots, a_{r_1}^{t_{r_1}^1} + c_{r_1}^1\right) \cup$$

$$\left(a_1^{t_1^2} + c_1^2, a_2^{t_2^2} + c_2^2, \ldots, a_{r_2}^{t_{r_2}^1} + c_{r_2}^2\right) \cup$$

$$\ldots \cup \left(a_1^{t_1^n} + c_1^n, a_1^{t_2^n} + c_2^n, \ldots, a_{r_n}^{t_{r_n}^n} + c_{r_n}^n\right).$$

If the code word x has been transmitted then the error n vector e is given as $e = y - x \in I - x = I$ i.e.,

$$\left\{e_1^1, e_2^1, \ldots, e_{r_1}^1\right\} \cup \left\{e_1^2, e_2^2, \ldots, e_{r_2}^2\right\} \cup \ldots \cup \left\{e_1^n, e_2^n, \ldots, e_{r_n}^n\right\} =$$
$$\left\{y_1^1, y_2^1, \ldots, y_{r_1}^1\right\} \cup \left\{y_1^2, y_2^2, \ldots, y_{r_2}^2\right\} \cup \ldots \cup \left\{y_1^n, y_2^n, \ldots, y_{r_n}^n\right\} -$$
$$\left\{x_1^1, x_2^1, \ldots, x_{r_1}^1\right\} \cup \left\{x_1^2, x_2^2, \ldots, x_{r_2}^2\right\} \cup \ldots \cup \left\{x_1^n, x_2^n, \ldots, x_{r_n}^n\right\}$$
$$\in I - x = I.$$

We have the decoding rule for each vector $y_{j_i}^i$ to find the coset to which it belongs. The minimum weigtht in a coset is called the coset leader. If we have several vector we choose one of them as the coset leader thus in case of group n-code we have n sets of coset leaders.

We will give a simple illustration for a code $C_{j_i}^i$; the coset table this can be applied for $1 \le j_i \le n_i$. $i = 1, 2, \ldots, n$.

***Example 3.2.10:*** Let

$$C = C_1 \cup C_2 \cup C_3 \cup C_4$$
$$= \left\{C_1^1, C_2^1\right\} \cup \left\{C_1^2, C_2^2, C_3^2\right\} \cup \left\{C_1^3, C_2^3\right\} \cup \left\{C_1^4, C_2^4, C_3^4\right\}$$

where

$$\{\, H_1^1 = (1\ 1\ 1\ 1\ 1)$$

$$\left. H_2^1 = \begin{pmatrix} 1 & 0 & 0 & 1 & 0 & 0 \\ 0 & 1 & 0 & 0 & 1 & 0 \\ 0 & 0 & 1 & 0 & 0 & 1 \end{pmatrix} \right\} \cup$$

$$\{\, H_1^2 = (1\ 1\ 1\ 1\ 1\ 1\ 1\ 1),$$

$$\begin{pmatrix} 1 & 1 & 1 & 0 & 1 & 0 & 0 \\ 0 & 1 & 1 & 1 & 0 & 1 & 0 \\ 1 & 1 & 0 & 1 & 0 & 0 & 1 \end{pmatrix} = H_2^2,$$

$$H_3^2 = \begin{pmatrix} 1 & 0 & 0 & 1 & 0 & 0 \\ 0 & 1 & 0 & 0 & 1 & 0 \\ 0 & 0 & 1 & 0 & 0 & 1 \end{pmatrix} \cup$$

$$\{(1\ 1\ 1\ 1\ 1\ 1) = H_1^3,$$

$$\begin{pmatrix} 0 & 0 & 1 & 0 & 1 & 1 & 1 \\ 0 & 1 & 0 & 1 & 1 & 1 & 0 \\ 1 & 0 & 1 & 1 & 1 & 0 & 0 \end{pmatrix} = H_2^3 \} \cup$$

$$\{(1\ 1\ 1\ 1\ 1\ 1) = H_1^4,$$

$$\begin{pmatrix} 1 & 1 & 0 & 0 & 0 \\ 1 & 0 & 1 & 0 & 0 \\ 1 & 0 & 0 & 1 & 0 \\ 1 & 0 & 0 & 0 & 1 \end{pmatrix} = H_2^4,$$

$$H_3^4 = \left.\begin{pmatrix} 0 & 1 & 1 & 1 & 0 & 0 \\ 1 & 0 & 1 & 0 & 1 & 0 \\ 1 & 1 & 0 & 0 & 0 & 1 \end{pmatrix}\right\}$$

is the group parity 4 check matrix associated with C.

Suppose

$y = y_1 \cup y_2 \cup y_3 \cup y_4$

$= \{(1\ 1\ 0\ 1\ 1)\ (1\ 1\ 0\ 1\ 0\ 0)\}$

$\cup \{(1\ 1\ 0\ 1\ 1\ 0\ 1\ 1)\ (1\ 1\ 1\ 0\ 1\ 0\ 1)\ (1\ 1\ 1\ 1\ 1\ 0)\}$

$\cup \{(1\ 1\ 0\ 0\ 0\ 0)\ (1\ 1\ 0\ 1\ 1\ 1\ 0)\}$

$\cup \{(1\ 0\ 1\ 0\ 1\ 1), (1\ 1\ 1\ 1\ 0)\ (1\ 1\ 1\ 1\ 0\ 1)\}$

be a received group 4 code word.

What is the procedure used to detect the error first?

We calculate the group 4 syndrome

$$GS(y) = \left\{ \begin{pmatrix} 1 & 1 & 1 & 1 & 1 \end{pmatrix} \begin{bmatrix} 1 \\ 1 \\ 0 \\ 1 \\ 1 \end{bmatrix}, \begin{pmatrix} 1 & 0 & 0 & 1 & 0 & 0 \\ 0 & 1 & 0 & 0 & 1 & 0 \\ 0 & 0 & 1 & 0 & 0 & 1 \end{pmatrix} \begin{bmatrix} 1 \\ 1 \\ 0 \\ 1 \\ 0 \\ 0 \end{bmatrix} \right\}$$

$$\cup \left\{ \begin{pmatrix} 1 & 1 & 1 & 1 & 1 & 1 & 1 & 1 \end{pmatrix} \begin{bmatrix} 1 \\ 1 \\ 0 \\ 1 \\ 1 \\ 0 \\ 1 \\ 1 \end{bmatrix}, \right.$$

$$\left. \begin{pmatrix} 1 & 1 & 1 & 0 & 1 & 0 & 0 \\ 0 & 1 & 1 & 1 & 0 & 1 & 0 \\ 1 & 1 & 0 & 1 & 0 & 0 & 1 \end{pmatrix} \begin{bmatrix} 1 \\ 1 \\ 1 \\ 0 \\ 1 \\ 0 \\ 1 \end{bmatrix}, \begin{pmatrix} 1 & 0 & 0 & 1 & 0 & 0 \\ 0 & 1 & 0 & 0 & 1 & 0 \\ 0 & 0 & 1 & 0 & 0 & 1 \end{pmatrix} \begin{bmatrix} 1 \\ 1 \\ 1 \\ 1 \\ 1 \\ 0 \end{bmatrix} \right\}$$

$$\cup \left\{ \begin{pmatrix}1 & 1 & 1 & 1 & 1 & 1\end{pmatrix} \begin{bmatrix}1\\1\\0\\0\\0\\0\end{bmatrix}, \begin{pmatrix}0 & 0 & 1 & 0 & 1 & 1 & 1\\0 & 1 & 0 & 1 & 1 & 1 & 0\\1 & 0 & 1 & 1 & 1 & 0 & 0\end{pmatrix} \begin{bmatrix}1\\1\\0\\1\\1\\1\\0\end{bmatrix} \right\}$$

$$\cup \left\{ \begin{pmatrix}1 & 1 & 1 & 1 & 1 & 1\end{pmatrix} \begin{bmatrix}1\\0\\1\\0\\1\\1\end{bmatrix}, \right.$$

$$\left. \begin{pmatrix}1 & 1 & 0 & 0 & 0\\1 & 0 & 1 & 0 & 0\\1 & 0 & 0 & 1 & 0\\1 & 0 & 0 & 0 & 1\end{pmatrix} \begin{bmatrix}1\\1\\1\\1\\0\end{bmatrix}, \begin{pmatrix}0 & 1 & 1 & 1 & 0 & 0\\1 & 0 & 1 & 0 & 1 & 0\\1 & 1 & 0 & 0 & 0 & 1\end{pmatrix} \begin{bmatrix}1\\1\\1\\1\\0\\1\end{bmatrix} \right\}$$

= {(0) (0 1 0)} $\cup$ {(0) (0 0 1) (0 0 1)} $\cup$ {(0) (0 0 1)} $\cup$ {(0) (0 0 0 1) (1 0 1)} $\neq$ {(0) (0 0 0)} $\cup$ {(0) (0 0 0) (0 0 0)} $\cup$ {(0) (0 0 0)} $\cup$ {(0) (0 0 0 0) (0 0 0)}. So the received message has an error.

Now we use n-coset leader method to correct the codes. First we see the error during transmission. We see error has occurred in the code $C_2^1, C_2^2 C_3^2, C_2^3, C_2^4$ and $C_3^4$.

However $C_2^4$ can be corrected as (1 1 1 1 1) for a repetition code has (0 0 0 0 0) or (1 1 1 1 1). Thus we have to correct the code words from $C_2^1, C_1^2, C_3^2, C_2^3$ and $C_3^4$ only.

We give the procedure how it is corrected using coset leaders.

| Message words | 0 0 0 | 0 1 0 | 0 0 1 | 1 0 0 |
|---|---|---|---|---|
| Code words | 0 0 0 0 0 0 | 0 1 0 0 1 0 | 0 0 1 0 0 1 | 1 0 0 1 0 0 |
| Other cosets | 1 0 0 0 0 0<br>0 1 0 0 0 0<br>0 0 1 0 0 0<br>1 1 0 0 0 0<br>0 1 1 0 0 0<br>1 0 0 0 0 1<br>1 1 1 0 0 0 | 1 1 0 0 1 0<br>0 0 0 0 1 0<br>0 1 1 0 1 0<br>1 0 0 0 1 0<br>0 0 1 0 1 0<br>1 1 0 0 1 1<br>1 0 1 0 1 0 | 1 0 1 0 0 1<br>0 1 1 0 0 1<br>0 0 0 0 0 1<br>1 1 1 0 0 1<br>0 1 0 0 0 1<br>1 0 1 1 0 1<br>1 1 0 0 0 1 | 0 0 0 1 0 0<br>1 1 0 1 0 0<br>1 0 1 1 0 0<br>0 1 0 1 0 0<br>1 1 1 1 0 0<br>0 0 0 1 0 1<br>0 1 1 1 0 0 |
| Message words | 1 1 0 | 0 1 1 | 1 0 1 | 1 1 1 |
| Code words | 1 1 0 1 1 0 | 0 1 1 0 1 1 | 1 0 1 1 0 1 | 1 1 1 1 1 1 |
| Other cosets | 0 1 0 1 1 0<br>1 0 0 1 1 0<br>1 1 1 0 1 1<br>0 0 0 1 1 0<br>1 0 1 1 1 0<br>0 1 0 1 1 1<br>0 0 1 1 1 0 | 1 1 1 0 1 1<br>0 0 1 0 1 1<br>0 1 0 0 1 1<br>1 0 1 0 1 1<br>0 0 0 0 1 1<br>1 1 1 0 1 0<br>1 0 0 0 1 1 | 0 0 1 1 0 1<br>1 1 1 1 0 1<br>1 0 0 1 0 1<br>0 1 1 1 0 1<br>1 1 0 1 0 1<br>0 0 1 1 0 0<br>0 1 0 1 0 1 | 0 1 1 1 1 1<br>1 0 1 1 1 1<br>1 1 0 1 1 1<br>0 0 1 1 1 1<br>1 0 0 1 1 1<br>0 1 1 1 1 0<br>0 0 0 1 1 1 |

The correct message is (1 0 0 1 0 0) from the code $C_2^1$ and not (1 1 0 1 0 0). The message (1 1 1 0 1 0 1) is not a correct one, to the correct once we again use the coset leader method, which is as follows:

| | | | | | | | | | | | | | | |
|---|---|---|---|---|---|---|---|---|---|---|---|---|---|---|
| message | 0 | 0 | 0 | 0 | | | | 1 | 0 | 0 | 0 | 1 | | |
| code words | 0 | 0 | 0 | 0 | 0 | 0 | 0 | 1 | 0 | 0 | 0 | 1 | 0 | 1 |
| | 1 | 0 | 0 | 0 | 0 | 0 | 0 | 0 | 0 | 0 | 0 | 1 | 0 | 1 |
| | 0 | 1 | 0 | 0 | 0 | 0 | 0 | 1 | 1 | 0 | 0 | 1 | 0 | 1 |
| | 0 | 0 | 1 | 0 | 0 | 0 | 0 | 1 | 0 | 1 | 0 | 1 | 0 | 1 |
| | 0 | 0 | 0 | 1 | 0 | 0 | 0 | 1 | 0 | 0 | 1 | 1 | 0 | 1 |
| | 0 | 0 | 0 | 0 | 1 | 0 | 0 | 1 | 0 | 0 | 0 | 0 | 0 | 1 |
| | 0 | 0 | 0 | 0 | 0 | 1 | 0 | 1 | 0 | 0 | 0 | 1 | 1 | 1 |
| | 0 | 0 | 0 | 0 | 0 | 0 | 1 | 1 | 0 | 0 | 0 | 1 | 0 | 0 |

| | | | | | | | | | | | | | |
|---|---|---|---|---|---|---|---|---|---|---|---|---|---|
| 0 | 1 | 0 | 0 | | | | 0 | 0 | 1 | 0 | | | |
| 0 | 1 | 0 | 0 | 1 | 1 | 1 | 0 | 0 | 1 | 0 | 1 | 1 | 0 |
| 1 | 1 | 0 | 0 | 1 | 1 | 1 | 1 | 0 | 1 | 0 | 1 | 1 | 0 |
| 0 | 0 | 0 | 0 | 1 | 1 | 1 | 0 | 1 | 1 | 0 | 1 | 1 | 0 |
| 0 | 1 | 1 | 0 | 1 | 1 | 1 | 0 | 1 | 0 | 0 | 1 | 1 | 0 |
| 0 | 1 | 0 | 1 | 1 | 1 | 1 | 0 | 0 | 1 | 1 | 1 | 1 | 0 |
| 0 | 1 | 0 | 0 | 0 | 1 | 1 | 0 | 0 | 1 | 0 | 0 | 1 | 0 |
| 0 | 1 | 0 | 0 | 1 | 0 | 1 | 0 | 0 | 1 | 0 | 1 | 0 | 0 |
| 0 | 1 | 0 | 0 | 1 | 1 | 0 | 0 | 0 | 1 | 0 | 1 | 1 | 1 |
| 0 | 0 | 0 | 1 | | | | 1 | 1 | 0 | 0 | | | |
| 0 | 0 | 0 | 1 | 0 | 1 | 1 | 1 | 1 | 0 | 0 | 0 | 1 | 0 |
| 1 | 0 | 0 | 1 | 0 | 1 | 1 | 0 | 1 | 0 | 0 | 0 | 1 | 0 |
| 0 | 1 | 0 | 1 | 0 | 1 | 1 | 1 | 0 | 0 | 0 | 0 | 1 | 0 |
| 0 | 0 | 1 | 1 | 0 | 1 | 1 | 1 | 1 | 1 | 0 | 0 | 1 | 0 |
| 0 | 0 | 0 | 0 | 0 | 1 | 1 | 1 | 1 | 0 | 1 | 0 | 1 | 0 |
| 0 | 0 | 0 | 1 | 1 | 1 | 1 | 1 | 1 | 0 | 0 | 1 | 1 | 0 |
| 0 | 0 | 0 | 1 | 0 | 0 | 1 | 1 | 1 | 0 | 0 | 0 | 0 | 0 |
| 0 | 0 | 0 | 1 | 0 | 1 | 0 | 1 | 1 | 0 | 0 | 0 | 1 | 1 |

| 1 | 0 | 1 | 0 | | | |
|---|---|---|---|---|---|---|
| 1 | 0 | 1 | 0 | 0 | 1 | 1 |
| 0 | 0 | 1 | 0 | 0 | 1 | 1 |
| 1 | 1 | 1 | 0 | 0 | 1 | 1 |
| 1 | 0 | 0 | 0 | 0 | 1 | 1 |
| 1 | 0 | 1 | 1 | 0 | 1 | 1 |
| 1 | 0 | 1 | 1 | 0 | 1 | 1 |
| 1 | 0 | 1 | 0 | 0 | 0 | 1 |
| 1 | 0 | 1 | 0 | 0 | 1 | 0 |

| 1 | 0 | 0 | 1 | | | |
|---|---|---|---|---|---|---|
| 1 | 0 | 0 | 1 | 1 | 1 | 0 |
| 0 | 0 | 0 | 1 | 1 | 1 | 0 |
| 1 | 1 | 0 | 1 | 1 | 1 | 0 |
| 1 | 0 | 1 | 1 | 1 | 1 | 0 |
| 1 | 0 | 0 | 0 | 1 | 1 | 0 |
| 1 | 0 | 0 | 1 | 0 | 1 | 0 |
| 1 | 0 | 0 | 1 | 1 | 0 | 0 |
| 1 | 0 | 0 | 1 | 1 | 1 | 1 |

| 0 | 1 | 0 | 1 | | | |
|---|---|---|---|---|---|---|
| 0 | 1 | 0 | 1 | 1 | 0 | 1 |
| 1 | 1 | 0 | 1 | 1 | 0 | 1 |
| 0 | 0 | 0 | 1 | 1 | 0 | 1 |
| 0 | 1 | 1 | 1 | 1 | 0 | 1 |
| 0 | 1 | 0 | 0 | 1 | 0 | 1 |
| 0 | 1 | 0 | 1 | 0 | 0 | 1 |
| 0 | 1 | 0 | 1 | 1 | 1 | 1 |
| 0 | 1 | 0 | 1 | 1 | 0 | 0 |

| 0 | 1 | 1 | 0 | | | |
|---|---|---|---|---|---|---|
| 0 | 1 | 1 | 0 | 0 | 0 | 1 |
| 1 | 1 | 1 | 0 | 0 | 0 | 1 |
| 0 | 0 | 1 | 0 | 0 | 0 | 1 |
| 0 | 1 | 0 | 0 | 0 | 0 | 1 |
| 0 | 1 | 1 | 1 | 0 | 0 | 1 |
| 0 | 1 | 1 | 0 | 1 | 0 | 1 |
| 0 | 1 | 1 | 0 | 0 | 1 | 1 |
| 0 | 1 | 1 | 0 | 0 | 0 | 0 |

| 0 | 0 | 1 | 1 | | | |
|---|---|---|---|---|---|---|
| 0 | 0 | 1 | 1 | 1 | 0 | 1 |
| 1 | 0 | 1 | 1 | 1 | 0 | 1 |
| 0 | 1 | 1 | 1 | 1 | 0 | 1 |
| 0 | 0 | 0 | 1 | 1 | 0 | 1 |
| 0 | 0 | 1 | 0 | 1 | 0 | 1 |
| 0 | 0 | 1 | 1 | 0 | 0 | 1 |
| 0 | 0 | 1 | 1 | 1 | 1 | 1 |
| 0 | 0 | 1 | 1 | 1 | 0 | 0 |

| 1 | 1 | 1 | 0 | | | |
|---|---|---|---|---|---|---|
| 1 | 1 | 1 | 0 | 1 | 0 | 0 |
| 0 | 1 | 1 | 0 | 1 | 0 | 0 |
| 1 | 0 | 1 | 0 | 1 | 0 | 0 |
| 1 | 1 | 0 | 0 | 1 | 0 | 0 |
| 1 | 1 | 1 | 1 | 1 | 0 | 0 |
| 1 | 1 | 1 | 0 | 0 | 0 | 0 |
| 1 | 1 | 1 | 0 | 1 | 1 | 0 |
| 1 | 1 | 1 | 0 | 1 | 0 | 1 |

| 1 | 1 | 0 | 1 | | | |
|---|---|---|---|---|---|---|
| 1 | 1 | 0 | 1 | 1 | 0 | 1 |
| 0 | 1 | 0 | 1 | 1 | 0 | 1 |
| 1 | 0 | 0 | 1 | 1 | 0 | 1 |
| 1 | 1 | 1 | 1 | 1 | 0 | 1 |
| 1 | 1 | 0 | 0 | 1 | 1 | 0 |
| 1 | 1 | 0 | 1 | 0 | 0 | 1 |
| 1 | 1 | 0 | 1 | 1 | 1 | 1 |
| 1 | 1 | 0 | 1 | 1 | 0 | 0 |

| | | |
|---|---|---|
| 1 0 1 1 | 0 1 1 1 | 1 1 1 1 |
| 1 0 1 1 0 0 0 | 0 1 1 1 0 1 0 | 1 1 1 1 1 1 1 |
| 0 0 1 1 0 0 0 | 1 1 1 1 0 1 0 | 0 1 1 1 1 1 1 |
| 1 1 1 1 0 0 0 | 0 0 1 1 0 1 0 | 1 0 1 1 1 1 1 |
| 1 0 0 1 0 0 0 | 0 1 0 1 0 1 0 | 1 1 0 1 1 1 1 |
| 1 0 1 0 0 0 0 | 0 1 1 0 0 1 0 | 1 1 1 0 1 1 1 |
| 1 0 1 1 1 0 0 | 0 1 1 1 1 1 0 | 1 1 1 1 0 1 1 |
| 1 0 1 1 0 1 0 | 0 1 1 1 0 0 0 | 1 1 1 1 1 0 1 |
| 1 0 1 1 0 0 1 | 0 1 1 1 0 1 1 | 1 1 1 1 1 1 0 |

The received word is (1 1 1 0 1 0 1) the correct set word must be (1 1 1 0 1 0 0) from $C_2^2$.

Now the received word (1 1 1 1 1 0) is not a code word of $C_3^2$ to find the correct word we make use of the coset leader method. For this we make use of the coset representation of the code $C_2^1$. We find this word (1 1 1 1 1 0) occurs with the word (1 1 0 1 1 0) so the correct message is (1 1 0 1 1 0) and this occurs with the coset leader (0 0 1 0 0 0) so (0 0 1 0 0 0) + (1 1 1 1 1 0) = (1 1 0 1 1 0).

Now the received message (1 1 0 1 1 1 0) is not a correct code word. To find the approximately proper word using the method of coset leaders

| | | |
|---|---|---|
| Message | 0 0 0 0 | 1 0 0 0 |
| codewords | 0 0 0 0 0 0 0 | 1 0 0 0 1 1 0 |
| | 1 0 0 0 0 0 0 | 0 0 0 0 1 1 0 |
| | 0 1 0 0 0 0 0 | 1 1 0 0 1 1 0 |
| | 0 0 1 0 0 0 0 | 1 0 1 0 1 1 0 |
| | 0 0 0 1 0 0 0 | 1 0 0 1 1 0 0 |
| | 0 0 0 0 1 0 0 | 1 0 0 0 0 1 0 |
| | 0 0 0 0 0 1 0 | 1 0 0 0 1 0 0 |
| | 0 0 0 0 0 0 1 | 1 0 0 0 1 1 1 |
| 0 1 0 0 | 0 0 1 0 | 0 0 0 1 |
| 0 1 0 0 0 1 1 | 0 0 1 0 1 1 1 | 0 0 0 1 1 0 1 |
| 1 1 0 0 0 1 1 | 1 0 1 0 1 1 1 | 1 0 0 1 1 0 1 |
| 0 0 0 0 0 1 1 | 0 1 1 0 1 1 1 | 0 1 0 1 1 0 1 |
| 0 1 1 0 0 1 1 | 0 0 0 0 1 1 1 | 0 0 1 1 1 0 1 |
| 0 1 0 1 0 1 1 | 0 0 1 1 1 1 1 | 0 0 0 0 1 0 1 |
| 0 1 0 0 1 1 1 | 0 0 1 0 0 1 1 | 0 0 0 1 0 0 1 |
| 0 1 0 0 0 0 1 | 0 0 1 0 1 0 1 | 0 0 0 1 1 1 1 |
| 0 1 0 0 0 1 0 | 0 0 1 0 1 1 0 | 0 0 0 1 1 0 0 |
| 1 1 0 0 | 0 1 1 0 | 0 0 1 1 |
| 1 1 0 0 1 0 1 | 0 1 1 0 1 0 0 | 0 0 1 1 0 1 0 |
| 0 1 0 0 1 0 1 | 1 1 1 0 1 0 0 | 1 0 1 1 0 1 0 |
| 1 0 0 0 1 0 1 | 0 0 1 0 1 0 0 | 0 1 1 1 0 1 0 |
| 1 1 1 0 1 0 1 | 0 1 0 0 1 0 0 | 0 0 0 1 0 1 0 |
| 1 1 0 1 1 0 1 | 0 1 1 1 1 0 0 | 0 0 1 0 0 1 0 |
| 1 1 0 0 0 0 1 | 0 1 1 0 0 0 0 | 0 0 1 1 1 1 0 |
| 1 1 0 0 1 1 1 | 0 1 1 0 1 1 0 | 0 0 1 1 0 0 0 |
| 1 1 0 0 1 0 0 | 0 1 1 0 1 0 1 | 0 0 1 1 0 1 1 |

$$
\begin{array}{ccccccc|ccccccc|ccccccc|}
1 & 0 & 1 & 0 & & & & 1 & 0 & 0 & 1 & & & & 0 & 1 & 0 & 1 & & & \\
\hline
1 & 0 & 1 & 0 & 0 & 0 & 1 & 1 & 0 & 0 & 1 & 0 & 1 & 1 & 0 & 1 & 0 & 1 & 1 & 1 & 0 \\
\hline
0 & 0 & 1 & 0 & 0 & 0 & 1 & 0 & 0 & 0 & 1 & 0 & 1 & 1 & 1 & 1 & 0 & 1 & 1 & 1 & 0 \\
1 & 1 & 1 & 0 & 0 & 0 & 1 & 1 & 1 & 0 & 1 & 0 & 1 & 1 & 0 & 0 & 0 & 1 & 1 & 1 & 0 \\
1 & 0 & 0 & 0 & 0 & 0 & 1 & 1 & 0 & 1 & 1 & 0 & 1 & 1 & 0 & 1 & 1 & 1 & 1 & 1 & 0 \\
1 & 0 & 1 & 1 & 0 & 0 & 1 & 1 & 0 & 0 & 0 & 0 & 1 & 1 & 0 & 1 & 0 & 0 & 1 & 1 & 0 \\
1 & 0 & 1 & 0 & 1 & 0 & 1 & 1 & 0 & 0 & 1 & 1 & 1 & 1 & 0 & 1 & 0 & 1 & 0 & 1 & 0 \\
1 & 0 & 1 & 0 & 0 & 1 & 1 & 1 & 0 & 0 & 1 & 0 & 0 & 1 & 0 & 1 & 0 & 1 & 1 & 0 & 0 \\
1 & 0 & 1 & 0 & 0 & 0 & 0 & 1 & 0 & 0 & 1 & 0 & 1 & 0 & 0 & 1 & 0 & 1 & 1 & 1 & 1
\end{array}
$$

$$
\begin{array}{ccccccc|ccccccc|ccccccc|}
1 & 1 & 1 & 0 & & & & 0 & 1 & 1 & 1 & & & & 1 & 1 & 0 & 1 & & & \\
\hline
1 & 1 & 1 & 0 & 0 & 1 & 0 & 0 & 1 & 1 & 1 & 0 & 0 & 1 & 1 & 1 & 0 & 1 & 0 & 0 & 0 \\
\hline
0 & 1 & 1 & 0 & 0 & 1 & 0 & 1 & 1 & 1 & 1 & 0 & 0 & 1 & 0 & 1 & 0 & 1 & 0 & 0 & 0 \\
1 & 0 & 1 & 0 & 0 & 1 & 0 & 0 & 0 & 1 & 1 & 0 & 0 & 1 & 1 & 0 & 0 & 1 & 0 & 0 & 0 \\
1 & 1 & 0 & 0 & 0 & 1 & 0 & 0 & 1 & 0 & 1 & 0 & 0 & 1 & 1 & 1 & 1 & 1 & 0 & 0 & 0 \\
1 & 1 & 1 & 1 & 0 & 1 & 0 & 0 & 1 & 1 & 0 & 0 & 0 & 1 & 1 & 1 & 0 & 0 & 0 & 0 & 0 \\
1 & 1 & 1 & 0 & 1 & 1 & 0 & 0 & 1 & 1 & 1 & 1 & 0 & 1 & 1 & 1 & 0 & 1 & 1 & 0 & 0 \\
1 & 1 & 1 & 0 & 0 & 0 & 0 & 0 & 1 & 1 & 1 & 0 & 1 & 1 & 1 & 1 & 0 & 1 & 0 & 1 & 0 \\
1 & 1 & 1 & 0 & 0 & 1 & 1 & 0 & 1 & 1 & 1 & 0 & 0 & 0 & 1 & 1 & 0 & 1 & 0 & 0 & 1
\end{array}
$$

$$
\left.\begin{array}{ccccccc|ccccccc}
1 & 0 & 1 & 1 & & & & 1 & 1 & 1 & 1 & & & \\
\hline
1 & 0 & 1 & 1 & 1 & 0 & 0 & 1 & 1 & 1 & 1 & 1 & 1 & 1 \\
\hline
0 & 0 & 1 & 1 & 1 & 0 & 0 & 0 & 1 & 1 & 1 & 1 & 1 & 1 \\
1 & 1 & 1 & 1 & 1 & 0 & 0 & 1 & 0 & 1 & 1 & 1 & 1 & 1 \\
1 & 0 & 0 & 1 & 1 & 0 & 0 & 1 & 1 & 0 & 1 & 1 & 1 & 1 \\
1 & 0 & 1 & 0 & 1 & 0 & 0 & 1 & 1 & 1 & 0 & 1 & 1 & 1 \\
1 & 0 & 1 & 1 & 0 & 0 & 0 & 1 & 1 & 1 & 1 & 0 & 1 & 1 \\
1 & 0 & 1 & 1 & 1 & 1 & 0 & 1 & 1 & 1 & 1 & 1 & 0 & 1 \\
1 & 0 & 1 & 1 & 1 & 0 & 1 & 1 & 1 & 1 & 1 & 1 & 1 & 0
\end{array}\right].
$$

We know (1 1 0 1 1 1 0) is the a received word clearly (1 1 0 1 1 1 0) is not a correct word. The coset leader is (1 0 0 0 0 0 0) the correct word is (1 0 0 0 0 0 0) + (1 1 0 1 1 1 0) = (0 1 0 1 1 1 0).

Clearly (1 1 1 1 0) is a received word which is not a correct word as a code is a repetition code, so (1 1 1 1 1) is the correct word.

Now we consider the received code word (1 1 1 1 0 1); clearly some error as occurred during transmission.

To find the correct message using the technique of coset leaders.

The message

| | | |
|---|---|---|
| 0 0 0 | 1 0 0 | 0 0 1 |
| 0 0 0 0 0 0 | 1 0 0 0 1 1 | 0 0 1 1 1 0 |
| 1 0 0 0 0 0 | 0 0 0 0 1 1 | 1 0 1 1 1 0 |
| 0 1 0 0 0 0 | 1 1 0 0 1 1 | 0 1 1 1 1 0 |
| 0 0 1 0 0 0 | 1 0 1 0 1 1 | 0 0 0 1 1 0 |
| 0 0 0 1 0 0 | 1 0 0 1 1 1 | 0 0 1 0 1 0 |
| 0 0 0 0 1 0 | 1 0 0 0 0 1 | 0 0 1 1 0 0 |
| 0 0 0 0 0 1 | 1 0 0 0 1 0 | 0 0 1 1 1 1 |
| 1 0 0 1 0 0 | 1 0 0 1 1 1 | 1 0 1 0 1 0 |
| 0 1 0 | 1 0 1 | 0 1 1 |
| 0 1 0 1 0 1 | 1 0 1 1 0 1 | 0 1 1 0 1 1 |
| 1 1 0 1 0 1 | 0 0 1 1 0 1 | 1 1 1 0 1 1 |
| 0 0 0 1 0 1 | 1 1 1 1 0 1 | 0 0 1 0 1 1 |
| 0 1 1 1 0 1 | 1 0 0 1 0 1 | 0 1 0 0 1 1 |
| 0 1 0 0 0 1 | 1 0 1 0 0 1 | 0 1 1 1 1 1 |
| 0 1 0 1 1 1 | 1 0 1 1 0 1 | 1 0 0 0 1 1 |
| 0 1 0 1 0 0 | 1 0 1 1 0 0 | 0 1 1 0 1 0 |
| 1 1 0 0 0 1 | 0 0 1 0 0 1 | 1 1 1 1 1 1 |

| | | | | | | | | | | | |
|---|---|---|---|---|---|---|---|---|---|---|---|
| 1 | 1 | 1 | | | | 1 | 1 | 0 | | | |
| 1 | 1 | 1 | 0 | 0 | 0 | 1 | 1 | 0 | 1 | 1 | 0 |
| 0 | 1 | 1 | 0 | 0 | 0 | 0 | 1 | 0 | 1 | 1 | 0 |
| 1 | 0 | 1 | 0 | 0 | 0 | 1 | 0 | 0 | 1 | 1 | 0 |
| 1 | 1 | 0 | 0 | 0 | 0 | 1 | 1 | 1 | 1 | 1 | 0 |
| 1 | 1 | 1 | 1 | 0 | 0 | 1 | 1 | 0 | 0 | 1 | 0 |
| 1 | 1 | 1 | 0 | 1 | 0 | 1 | 1 | 0 | 1 | 0 | 0 |
| 1 | 1 | 1 | 0 | 0 | 1 | 1 | 1 | 0 | 1 | 1 | 1 |
| 0 | 1 | 1 | 1 | 0 | 0 | 0 | 1 | 0 | 0 | 1 | 0 |

The coset leader is (0 1 0 0 0 0) and the received word is (1 1 1 1 0 1). Hence the correct word is (received word) t coset leader = (1 1 1 1 0 1) + (0 1 0 0 0 0) = (1 0 1 1 0 1) is the corrected code word.

Thus we have elaborately described the method of both error detection and error correction using the coset leaders. Do we have any other method of finding the correct n-code code? We give the pseudo best group n approximations.

**DEFINITION 3.2.8:** *Let $G = (G_1 \cup G_2 \cup \ldots \cup G_n)$ $(n \geq 3)$ where $G_i = (c^i_1, \ldots, c^i_{r_i})$, $1 \leq i \leq n$ with $c^i_{t_i}$ codes i.e. $c^i_{t_i}$ is a subspace of a vector space $Z_2^{q_i^{t_i}}$, $q_i^j$ is the length of the code $c^i_{t_i}$ with j number of message symbols, $1 \leq j \leq r_i$; i = 1, 2, …, n, be a group n-code. If we have pseudo inner product on each $c^i_{t_i}$ which we denote by $\langle , \rangle_{t_i}$; $1 \leq t_i \leq r_i$; i = 1, 2, …, n.*

*Thus we have on G, $r_1 + r_2 + \ldots + r_n$ number of pseudo inner products.*

*A group n-space endowed with the $r_1 + \ldots + r_n$ pseudo inner product will be known as a $(r_1 + \ldots + r_n)$-pseudo n-inner product space.*

$$G \subseteq (Z_2^{q_1^1} \ldots Z_2^{q_1^{r_1}}) \cup (Z_2^{q_2^1} \ldots Z_2^{q_2^{r_2}}) \cup \ldots \cup (Z_2^{q_n^1} \ldots Z_2^{q_n^{r_n}})$$

*will also be a $(r_1 + \ldots + r_n)$-pseudo inner product subspace.*

We will illustrate this by some example.

***Example 3.2.11:*** Let

$$V = (V_1 \cup V_2 \cup V_3 \cup V_4)$$
$$= \{V_1^1, V_2^1, V_3^1\} \cup \{V_1^2, V_2^2\} \cup \{V_1^3, V_2^3, V_3^3, V_4^3\} \cup \{V_1^4, V_2^4\}$$

= $\{Z_2 \times Z_2 \times Z_2, Z_2 \times Z_2 \times Z_2 \times Z_2, Z_2 \times Z_2 \times Z_2 \times Z_2 \times Z_2 \times Z_2\} \cup \{Z_2 \times Z_2 \times Z_2 \times Z_2, Z_2 \times Z_2 \times Z_2 \times Z_2 \times Z_2 \times Z_2\} \cup \{Z_2 \times Z_2 \times Z_2, Z_2 \times Z_2, Z_2 \times Z_2 \times Z_2, Z_2 \times Z_2 \times Z_2 \times Z_2 \times Z_2 \times Z_2 \times Z_2 \times Z_2\} \cup \{Z_2 \times Z_2 \times Z_2 \times Z_2 \times Z_2 \times Z_2 \times Z_2, Z_2 \times Z_2 \times Z_2 \times Z_2 \times Z_2\}$ be a group vector 4-space over $Z_2$. Define a group 4 inner product $\langle\ \rangle = \{\langle\ \rangle_1^1, \langle\ \rangle_2^1, \langle\ \rangle_3^1\} \cup \{\langle\ \rangle_1^2, \langle\ \rangle_2^2\} \cup \{\langle\ \rangle_1^3, \langle\ \rangle_2^3, \langle\ \rangle_3^3, \langle\ \rangle_4^3\} \cup \{\langle\ \rangle_1^4, \langle\ \rangle_2^4\}$ on V where each $\langle , \rangle_j^i$ is the standard pseudo inner product on $V_j^i$; $1 \leq i \leq 4$, $1 \leq j \leq 3$ 3 or 2 or 4.

Thus the group vector 4-space is endowed with a standard pseudo 4-inner product.

Now we proceed onto define the notion of group pseudo best n-approximation.

**DEFINITION 3.2.9:** *Let*

$$C = C_1 \cup C_2 \cup \ldots \cup C_n$$
$$= \left(C_1^1, C_2^1, \ldots, C_{r_1}^1\right) \cup (C_1^2, C_2^2, \ldots, C_{r_2}^2) \cup \ldots \cup (C_1^n, C_2^n, \ldots, C_{r_n}^n)$$

*be a group n-code. Clearly each $C_j^i$ is a vector subspace of $Z_2^m$ suppose*

$$y = y_1 \cup y_2 \cup \ldots \cup y_n$$
$$= \{y_1^1, y_2^1, \ldots, y_{r_1}^1) \cup \{y_1^2, y_2^2, \ldots, y_{r_2}^2) \cup \ldots \cup \{y_1^n, y_2^n, \ldots, y_{r_n}^n)$$

*be a received code. Suppose $y \notin C$, then we find*

$$x = \{x_1^1, x_2^1, \ldots, x_{r_1}^1) \cup \{x_1^2, x_2^2, \ldots, x_{r_2}^2) \cup \ldots \cup \{x_1^n, x_2^n, \ldots, x_{r_n}^n)$$

*in C which is such that*

$$|| y - x || \leq || y - \gamma || \qquad (I)$$

*for every $\gamma$ in C i.e. the group best n-approximation to y by group n-vectors in C is a n-vector x in C such that (I) is satisfied.*

$$||y - x || = (|| y_1 - x_1 || \cup || y_2 - x_2 || \cup \ldots \cup || y_n - x_n ||)$$

$$= \left(\| y_1^1 - x_1^1 \|, \| y_2^1 - x_2^1 \|, \ldots, \| y_{r_1}^1 - x_{r_1}^1 \|\right) \cup \left(\| y_1^2 - x_1^2 \|,\right.$$
$$\| y_2^2 - x_2^2 \| , \ldots, \left.\| y_{r_2}^2 - x_{r_2}^2 \|\right) \cup \ldots \cup \left(\| y_1^n - x_1^n \|,\right.$$
$$\left.\| y_2^n - x_2^n \|, \ldots, \| y_{r_n}^n - x_{r_n}^n \|\right) \le \left(\| y_1^1 - \gamma_1^1 \|, \| y_2^1 - \gamma_2^1 \|, \ldots,\right.$$
$$\left.\| y_{r_1}^1 - \gamma_{r_1}^1 \|\right) \cup \left(\| y_1^2 - \gamma_1^2 \|, \| y_2^2 - \gamma_2^2 \|, \ldots, \| y_{r_2}^2 - \gamma_{r_2}^2 \|\right) \cup \ldots \cup$$
$$\left(\| y_1^n - \gamma_1^n \|, \| y_2^n - \gamma_2^n \|, \ldots, \| y_{r_n}^n - \gamma_{r_n}^n \|\right)$$

*i.e.* $\| y_{j_i}^i - x_{j_i}^i \| \le \| y_{j_i}^i - \gamma_{j_i}^i \|$; *$1 \le j_i \le r_i$; $i = 1, 2, \ldots, n$.*

*Thus*

$x = x_1 \cup x_2 \cup \ldots \cup x_n = (x_1^1, x_2^1, \ldots, x_{r_1}^1) \cup \ldots \cup (x_1^n, x_2^n, \ldots, x_{r_n}^n) \in C$

*is the group pseudo best n-approximation to y.*

For more literature please refer [39] we illustrate this by a simple example so that the reader can easily understand how to find the group best n-approximation to y; i.e. if $y = y_1 \cup \ldots \cup y_n$ is the received n-code word and y is not in the group n-code $C = C_1 \cup \ldots \cup C_n \subseteq (Z_2^{q_1^1} \ldots Z_2^{q_1^{r_1}}) \cup \ldots \cup (Z_2^{q_n^1} \ldots Z_2^{q_n^{r_n}})$. The group n- best approximation x to y with respect to group n-code C (which is a group n-subspace of $(Z_1 \cup \ldots \cup Z_n)$ gives the approximately correct group n-code word which is in C. i.e. $x \in C$ and x is such that $\| y - x\| \le \| y - \alpha\|$ for every $\alpha$ in C.

We shall illustrate this in a group 3-code.

***Example 3.2.12:*** Let $V = \{Z_2^4, Z_2^5\} \cup \{Z_2^5, Z_2^6\} \cup \{Z_2^6, Z_2^5\} = V_1 \cup V_2 \cup V_3$ be the group 3-vector space. Let $C = C_1 \cup C_2 \cup C_3 = \{c_1^1, c_2^1\} \cup \{c_2^1, c_2^2\} \cup \{c_3^1, c_3^2\}$ be a group code with the code words given by the group 3 parity check matrices

$$H = \{H_1^1, H_2^1\} \cup \{H_1^2, H_2^2\} \cup \{H_1^3, H_2^3\}$$

$$= \left\{ \begin{pmatrix} 1 & 0 & 1 & 0 \\ 1 & 1 & 0 & 1 \end{pmatrix}, \begin{pmatrix} 1 & 0 & 1 & 0 & 0 \\ 0 & 1 & 0 & 1 & 0 \\ 0 & 1 & 0 & 0 & 1 \end{pmatrix} \right\} \cup$$

$$\left\{\begin{pmatrix}1 & 0 & 0 & 1 & 0\\ 0 & 1 & 1 & 0 & 1\end{pmatrix}, \begin{pmatrix}0 & 1 & 1 & 1 & 0 & 0\\ 1 & 0 & 1 & 0 & 1 & 0\\ 1 & 1 & 0 & 0 & 0 & 1\end{pmatrix}\right\} \cup$$

$$\left\{\begin{pmatrix}1 & 0 & 1 & 1 & 0 & 0\\ 1 & 1 & 1 & 0 & 1 & 0\\ 0 & 1 & 1 & 0 & 0 & 1\end{pmatrix}, \begin{pmatrix}1 & 1 & 1 & 0 & 0\\ 1 & 0 & 0 & 1 & 0\\ 0 & 1 & 0 & 0 & 1\end{pmatrix}\right\}$$

be the 3-group code associated with the 3 parity check matrix H. Let us assume we have the standard inner product defined on the group 3 vector space V.

Now the group 3 code associated with H is given by $C = C_1 \cup C_2 \cup C_3$ = {{(0 0 0 0), (1 0 1 1), (0 1 0 1), (1 1 1 0)},{(0 0 0 0 0), (1 0 1 0 0), (0 1 0 1 1), (1 1 1 1 1)}} ∪ {{(0 0 0 0 0), (1 0 0 1 0) (0 1 0 0 1), (0 0 1 0 1), (1 1 0 1 1) (1 0 1 1 1) (0 1 1 0 0) (1 1 1 1 0)}, {(0 0 0 0 0 0), (1 0 0 0 1 1) (0 1 0 1 0 1), (0 0 1 1 1 0), (1 1 0 1 1 0) (0 1 1 0 1 1), (1 0 1 1 0 1), (1 1 1 0 0 0)}} ∪ {{(0 0 0 0 0 0), (1 0 0 1 1 0), (0 1 0 0 1 1) (0 0 1 1 1 1) (1 1 0 1 0 1) (0 1 1 1 0 0) (1 0 1 0 0 1) (1 1 1 0 1 0)}, {(0 0 0 0 0), (1 0 1 1 0) (0 1 1 0 1) (1 1 0 1 1)}}. Suppose the received group 3-code word is $y = y_1 \cup y_2 \cup y_3$ = {(1 1 1 1), (1 1 0 1 1)} ∪ {(0 1 0 1 0), (1 1 1 0 1 1)} ∪ {(1 1 1 1 1 0), (1 1 1 1 0)}. It is easily verified using the group 3-parity check matrix $y \notin C$. Now to find the group 3-best approximation to y relative to C the vector subspace of V. Clearly $y \in V$ and $y \notin C$. To find $x \in C$ such that $\| y - x \| \leq \|y - \alpha \|$ for all $\alpha \in C$. We use the standard pseudo group n-inner product on V.

Choose a group 3-orthonormal basis B for C. B = {{(1 0 1 1), (0 1 0 1)}, {(1 0 1 0 0), (0 1 0 1 1)} ∪ {{(1 0 0 1 0), (0 1 0 0 1), (0 0 1 0 1)}, {(1 0 0 0 1 1), (0 1 0 1 0 1), (0 0 1 1 1 0)}} ∪ {{(1 0 0 1 1 0), (0 1 0 0 1 1), (0 0 1 1 1 1)}, {(1 0 1 1 0), (0 1 1 0 1)}}.

$$x = \sum_i \sum_K (y \mid \alpha_K^i) \alpha_K^i$$

x = {(1 1 1 1), (1 0 1 1)⟩ (1 0 1 1) + ⟨(1 1 1 1), (0 1 0 1)⟩ (0 1 0 1)}, {⟨(1 1 0 1 1), (1 0 1 0 0)⟩ (1 0 1 0 0) + ⟨(1 1 0 1 1), (0 1 0 1 1)⟩ × (0 1 0 1 1)} ∪ {⟨(0 1 0 1 0), ((1 0 0 1 0)⟩ (1 0 0 1 0) + ⟨(0 1 0 1 0), (0 1 0 0 1)⟩ (0 1 0 0 1) + ⟨(0 1 0 1 0), (0 0 1 0 1)⟩ (0 0 1 0 1)}, {⟨(1 1 1 0 1 1), (1 0 0 0 1 1)⟩ (1 0 0 0 1 1) + ⟨(1 1 1 0 1 1), (0 1 0 1 0 1)⟩ (0 1 0 1 0 1) + ⟨(1 1 1 0 1 1), (0 0 1 1 1 0)⟩ × (0 0 1 1 1 0)}} ∪ {⟨(1 1 1 1 1 0), (1 0 0 1 1 0) × (1 0 0 1 1 0) + ⟨(1 1 1 1 1 0), (0 1 0 0 1 1)⟩ (0 1 0 0 1 1) + ⟨(1 1 1 1 1 0), (0 0 1 1 1 1)⟩ (0 0 1 1 1 1)}, {⟨(1 1 1 1 0), (1 0 1 1 0)⟩ (1 0 1 1 0) + ⟨(1 1 1 1 0), (0 1 1 0 1)⟩ (0 1 1 0 1)} = {(1 0 1 1), (1 1 1 1 1)} ∪ { (1 1 0 1 1), (1 0 0 0 1 1)} ∪ {(1 0 1 0 0 1), (1 0 1 0 0)} ∈ C. Thus x is the best 3-approximated group 3-codeword of y.

We mention a few of its applications.

1. In image compression and in the image coding these codes can do multifold job at a time. These codes would be much more welcome with the use of computers and processing would also be fast and economic.
2. These new codes can also be used in Block truncation coding.
3. Since not much of algebraic operations are used these codes can be readily made use of by cryptologists, computer scientist and electrical scientist.

These codes can be used as a storage device in computers. When sets of n-informations are to passed in channels where all the sets of n-informations are to be passed simultaneously need not be distinct repetition in each of the sets is also permitted. These codes can also be used in networking of computers when in the work place when some m of them work on different sets of $n_i$ datas i = 1, 2, …, m. These group n-codes can be used as storage were security is needed. To keep the secret in tact we use misleading codes to be present, so that the intruder cannot easily find out which, are true messages and which are false. For if

$$C = C_1 \cup C_2 \cup \ldots \cup C_n$$

$$= (C_1^1, C_2^1, \ldots, C_{r_1}^1) \cup (C_1^2, C_2^2, \ldots, C_{r_2}^2) \cup \ldots \cup (C_1^n, C_2^n, \ldots, C_{r_n}^n)$$

is the group n-code. In this we can have 60% of the codes to be real data storage one while the rest are just to mislead the intruder. Certainly it is not an easy task for any one to determine the real code which, carry the stored data and those misleading codes for one cannot find any difference between them. Further only one or two officers know the real codes in which the data is stored. Thus these codes will be very useful in defence departments.

Another way of using these codes is for transmissions. In that case if M is a message from a group n-code

$$C = C_1 \cup C_2 \cup \ldots \cup C_n$$

$$(C_1^1, C_2^1, \ldots, C_{r_1}^1) \cup (C_1^2, C_2^2, \ldots, C_{r_2}^2) \cup \ldots \cup (C_1^n, C_2^n, \ldots, C_{r_n}^n)$$

where

$$M = M_1 \cup M_2 \cup \ldots \cup M_n$$

$$= (M_1^1, M_2^1, \ldots, M_{r_1}^1) \cup (M_1^2, M_2^2, \ldots, M_{r_2}^2) \cup \ldots \cup (M_1^n, M_2^n, \ldots, M_{r_n}^n)$$

is the message to transmitted here only a stipulated percentage of the codewords $M_j^i$ carry the message, the rest are only added to mislead the eavesdropper, $1 \leq j \leq r_i$; $1 \leq i \leq n$. In such cases also the eavesdropper cannot easily get the real message. Thus these group n-codes have high security.

Further these can also be used by a group. It is further important to note every one in the group need not know the complete algorithm. Only the chief of the group knows the complete algorithm and its implementations. Each one in the group will only know the algorithm and the implementation of a few codes from $C_j = (C_1^j, C_2^j, \ldots, C_{r_j}^j)$. Also in these a few of the codes $C_s^j$ will be misleading codes. Thus the major merit of these group codes is that even if one accidentally reveals the secret every one need not change their algorithm. If the affected one changes it; it is sufficient. Thus the key will be made use of only on the codes which carry the messages and not on the misleading codes. The key space will be a group n-space. As in case of other algorithms the security is based on the key n-space. It does not mater even if the intruder (or eavesdropper) knows the algorithm for he cannot know the particular key for the key is different for each group member and it is very much far fetched for anyone to find every key. For if one member in the group is hacked nothing happens to the rest for each one has a different key infact a different type of code.

These group n-codes can be used by a n set of groups with varying set of members. Even each group can have a leader who alone knows the algorithms. Thus these codes happen to be more secure and more appropriate with the computerized world.

# FURTHER READING

# INDEX

**B**

**C**

**D**

**N**

**O**

**P**

**R**

**S**

**V**

**W**

# ABOUT THE AUTHORS

**Dr.W.B.Vasantha Kandasamy** is an Associate Professor in the Department of Mathematics, Indian Institute of Technology Madras, Chennai. In the past decade she has guided 12 Ph.D. scholars in the different fields of non-associative algebras, algebraic coding theory, transportation theory, fuzzy groups, and applications of fuzzy theory of the problems faced in chemical industries and cement industries.

She has to her credit 646 research papers. She has guided over 68 M.Sc. and M.Tech. projects. She has worked in collaboration projects with the Indian Space Research Organization and with the Tamil Nadu State AIDS Control Society. This is her $38^{th}$ book.

On India's 60th Independence Day, Dr.Vasantha was conferred the Kalpana Chawla Award for Courage and Daring Enterprise by the State Government of Tamil Nadu in recognition of her sustained fight for social justice in the Indian Institute of Technology (IIT) Madras and for her contribution to mathematics. (The award, instituted in the memory of Indian-American astronaut Kalpana Chawla who died aboard Space Shuttle Columbia). The award carried a cash prize of five lakh rupees (the highest prize-money for any Indian award) and a gold medal.
She can be contacted at vasanthakandasamy@gmail.com
You can visit her on the web at: http://mat.iitm.ac.in/~wbv

---

**Dr. Florentin Smarandache** is a Professor of Mathematics and Chair of Math & Sciences Department at the University of New Mexico in USA. He published over 75 books and 150 articles and notes in mathematics, physics, philosophy, psychology, rebus, literature.

In mathematics his research is in number theory, non-Euclidean geometry, synthetic geometry, algebraic structures, statistics, neutrosophic logic and set (generalizations of fuzzy logic and set respectively), neutrosophic probability (generalization of classical and imprecise probability). Also, small contributions to nuclear and particle physics, information fusion, neutrosophy (a generalization of dialectics), law of sensations and stimuli, etc. He can be contacted at smarand@unm.edu